\documentclass[11pt]{article}
\usepackage[pdftex]{graphicx}
\usepackage[small,bf,sf,textfont={small,sf}]{caption}
\usepackage{vmargin}
\usepackage{subfig}
\usepackage{amsmath}
\usepackage{amssymb}
\usepackage{fancyheadings}
\usepackage{color}
\usepackage[colorlinks=true,
            linkcolor=blue,
            urlcolor=blue,
            citecolor=blue]{hyperref}
\usepackage[numbers]{natbib}
\usepackage{doi}
\usepackage{booktabs}
\usepackage{float}
\usepackage{lineno}
\usepackage{epstopdf}
\usepackage{bm}
\usepackage{dsfont}
\usepackage{multirow}
\usepackage{pdflscape}

\floatstyle{ruled}
\newfloat{floatbox}{thp}{lob}
\floatname{floatbox}{Algorithm}

\newcommand{\trp}{^{\scriptsize \top}}

\newcommand{\inv}{^{\scriptsize -1}}
\newcommand{\sqr}{^{\scriptsize \textrm{1/2}}}
\newcommand{\invsqr}{^{\scriptsize -\textrm{1/2}}}

\newcommand{\x}{\mathbf{x}}

\newcommand{\m}{\mathbf{m}}

\newcommand{\e}{\mathbf{e}}

\newcommand{\mpr}{\mathbf{m}_{\textrm{pr}}}
\newcommand{\mtrue}{\mathbf{m}_{\textrm{true}}}

\newcommand{\wj}{\mathbf{w}_{j}}
\newcommand{\dsim}{\mathbf{d}}
\newcommand{\dtrue}{\mathbf{d}_{\textrm{true}}}
\newcommand{\dobs}{\mathbf{d}_{\textrm{obs}}}

\newcommand{\C}{\mathbf{C}}

\newcommand{\Cd}{\mathbf{C}_{\mathbf{d}}}

\newcommand{\Ce}{\mathbf{C}_{\mathbf{e}}}

\newcommand{\Cm}{\mathbf{C}_{\mathbf{m}}}

\newcommand{\K}{\mathbf{K}}
\newcommand{\R}{\mathbf{R}}

\newcommand{\DD}{\Delta \mathbf{D}}

\setpapersize{USletter}
\setmarginsrb{1.5in}{1in}{1in}{0.5in}{12pt}{11mm}{12pt}{30pt}

\begin{document}

\title{Data-Space Inversion with Ensemble Smoother}

\author{Mateus M. Lima$^1$, Alexandre A. Emerick$^2$ and Carlos E. P. Ortiz$^3$}

\maketitle

\footnotetext[1]{Petrobras and UENF (\texttt{mateusmartins@yahoo.com.br})}
\footnotetext[2]{Petrobras (\texttt{aemerick@gmail.com})}
\footnotetext[3]{UENF (\texttt{capico.LENEP@gmail.com})}

\section*{Abstract}
\label{Sec:Abstract}

Reservoir engineers use large-scale numerical models to predict the production performance in oil and gas fields. However, these models are constructed based on scarce and often inaccurate data, making their predictions highly uncertain. On the other hand, measurements of pressure and flow rates are constantly collected during the operation of the field. The assimilation of these data into the reservoir models (history matching) helps to mitigate uncertainty and improve their predictive capacity. History matching is a nonlinear inverse problem, which is typically handled using optimization and Monte Carlo methods. In practice, however, generating a set of properly history-matched models that preserve the geological realism is very challenging, especially in cases with complicated prior description, such as models with fractures and complex facies distributions. Recently, a new data-space inversion (DSI) approach was introduced in the literature as an alternative to the model-space inversion used in history matching. The essential idea is to update directly the predictions from a prior ensemble of models to account for the observed production history without updating the corresponding models. The present paper introduces a DSI implementation based on the use of an iterative ensemble smoother and demonstrates with examples that the new implementation is computationally faster and more robust than the earlier method based on principal component analysis. The new DSI is also applied to estimate the production forecast in a real field with long production history and a large number of wells. For this field problem, the new DSI obtained forecasts comparable with a more traditional ensemble-based history matching.

\noindent \textbf{Keywords:} Data-space inversion; uncertainty quantification; ensemble smoother; history matching.


\section{Introduction}
\label{Sec:Intro}

Reservoir characterization from static and dynamic data allows to create numerical models that can be used to simulate the performance of petroleum reservoirs under different operating conditions. These models are essential for efficient exploitation and management of oil and gas fields. The incorporation of static data is typically done using geostatistics while the incorporation of dynamic data is done using history-matching methods.

History matching is usually a very difficult task, which involves the integration of interdisciplinary teams and an intensive use of computational resources. A complete study may require months of work and the results are not always satisfactory. The task is even harder if one needs to provide uncertainty estimates, in which case several alternative history-matched models must be generated. In the last few decades, the advances in assisted (or semi-automatic) history-matching techniques were notorious. Yet, history matching remains one of the most time-consuming steps of a field study because the size and complexity of the models have also increased significantly in the same period. \citet{oliver:11b} present a review of the main history-matching methods proposed in the literature. Among these methods, the ones based on the ensemble Kalman filter (EnKF) \citep{evensen:94,evensen:07bk} have become quite popular, especially because of their ease of implementation and integration with commercial reservoir simulators and the ability to generate multiple models with large number of uncertainty parameters at an affordable computational cost. Despite the relative success in a number of recent field cases reported in the literature; see, for example, \citep{emerick:13c,chen:14a,emerick:16a,maucec:16a,abadpour:18a,evensen:18a,lorentzen:19a}, generating a set of models properly conditioned to all historical data and still preserving the geological realism is very challenging, especially in cases with complicated prior description, such as models with fractures and complex facies distributions.

A new approach known as data-space inversion (DSI) \citep{sun:17b} has drawn attention in the literature as an alternative to the model-space inversion approach used in history matching. The basic idea behind DSI is to update directly the predictions from a prior ensemble of models to account for the observed production history without updating the corresponding models. The upside of this approach is to be able to provide an ensemble of forecasts without going through the time-consuming history-matching step. Because there are no model inversions, there are no concerns about losing the geological realism. The downside is that DSI provides only the forecast estimates for a fixed production strategy. In practice, however, reservoir engineers may also be interested in having the corresponding models to study different drainage strategies. For this reason, data- and model-space inversions should be considered complementary rather than alternative (or even competing) approaches.

The DSI method introduced by \citet{sun:17b} uses principal component analysis (PCA) to reparameterize the predicted data from the prior ensemble into a lower dimensional space and the randomized maximum likelihood (RML) method \citep{oliver:08bk} to generate samples of the posterior distribution of predicted data given the observations. The authors used a data transformation before PCA to improve linearity. \citet{sun:17c} extended the original DSI method introducing a more general data transformation procedure. They tested the method in a model of a complex fractured reservoir and obtained reasonable uncertainty estimates of the production forecast. \citet{jiang:18a} modified the DSI approach to allow changes in the well controls during the forecast period so that the method could be used for life-cycle optimization. They noted that the method required a larger prior ensemble to better represent a wider range of possibilities. They showed that the proposed method combined to a direct search optimization algorithm \citep{audet:06a} was able to improve the expected net-present value of a reservoir model.

Similar ideas of DSI have appeared before in the literature. For example, \citet{krishnamurti:00a} and \citet{pagowski:05a} applied linear regression to combine the ozone forecast of ensembles of models. In the atmospheric literature these methods are referred to as aggregation methods or ensemble forecast \citep{mallet:09a}. \citet{emerick:14a} also used a DSI-type of approach to invert 4D seismic impedance data directly to pressure and water saturation. \citet{scheidt:15a} proposed a method named prediction-focused analysis (PFA) based on projecting the prior predictions into a low-dimensional space and using kernel smoothing to estimate the joint distribution of historical and forecasted data. Using the joint distribution, they could predict the uncertainty estimates in a tracer transport problem. However, PFA method seems applicable only to problems with few data points because it requires the projection to very few dimensions (two or three dimension in the examples presented in the paper). \citet{satija:15a} modified the PFA method using canonical functional component analysis to improve the linearity in the projected data. \citet{satija:17a} used the same approach in a reservoir problem and concluded that the method provided uncertainty estimates of production forecast in reasonable agreement with rejection sampling. More recently, \citet{he:18b} used similar ideas from DSI to estimate the uncertainty reduction in a study to compute the value of information of data-acquisition plans. \citet{jeong:18a} applied machine learning techniques (neural networks and support vector regression) to DSI. They concluded that the method can be a more efficient alternative to the computationally demanding history-matching methods, however the method fails to provide satisfactory forecast if the predictions from the prior ensemble (training set) are too far from the expected true response.

In the present paper, we introduce a new DSI implementation based on the use of an iterative ensemble smoother and demonstrate with examples that the new DSI is computationally faster and more robust than the procedure proposed in \citep{sun:17b,sun:17c}. Moreover, we apply the new DSI to a real field case with long production history and large number of wells and show that the method provides forecasts comparable with a more traditional ensemble-based history-matching process. The rest of the paper is organized as follows: Section~\ref{Sec:Methodology} reviews the DSI method proposed in \citep{sun:17b,sun:17c} and the new DSI method. Section~\ref{Sec:TestCases} presents three reservoir test problems. The first problem is a small synthetic case used to demonstrate that the proposed method provides results similar to the original DSI with a lower computational cost. The second problem is a benchmark history-matching case \citep{avansi:15a,unisim-i-h:13} constructed with data from real reservoir in Campos Basis. This problem is used to compare the methods is a more realistic situation with a large number of data points. The last problem corresponds to a real brown-field case where the proposed method is compared against an ensemble-based history matching. Section~\ref{Sec:Conclusions} summarizes the conclusions of the paper.

\section{Methodology}
\label{Sec:Methodology}

\subsection{Preliminaries}
\label{Sec:Preliminaries}

Let $\m \in \mathds{R}^{N_m}$ denote the vector of uncertain parameters of a reservoir model with a historical production period $t_{\textrm{h}}$. Our goal is to predict the production performance for a period $t_{\textrm{f}}$ after $t_{\textrm{h}}$. Let $\dsim \in \mathds{R}^{N_d}$ denote the vector of predicted production data, which is a nonlinear function of $\m$, that is, $\dsim = \dsim (\m)$. In the applications of interest of this paper, $\dsim$ is the result of a reservoir simulation. This vector contains predicted data from both, history, $\dsim_{\textrm{h}}$, and forecast periods, $\dsim_{\textrm{f}}$, that is,

\begin{equation}
  \dsim = \left[ \begin{array}{c}
                   \dsim_{\textrm{h}} \\
                   \dsim_{\textrm{f}}
                 \end{array}\right].
\end{equation}

Let $\dobs \in \mathds{R}^{N_{d,\textrm{h}}}$ denote the vector of field observations, which is corrupted with an additive random noise, $\e_{\textrm{d}}$, due to measurement errors such that

\begin{equation}
  \dobs = \dtrue + \e_{\textrm{d}},
\end{equation}
where $\dtrue$ is the true (noiseless) data. Moreover, assume that model errors are also additive and

\begin{equation}\label{Eq:dtrue}
   \dtrue = \dsim_{\textrm{h}} \left(\mtrue \right) + \e_{\textrm{m}},
\end{equation}
where $\mtrue$ is the vector containing the ``true'' values for the model parameters and $\e_{\textrm{m}}$ is a vector of model errors. Under these conditions, it is straightforward to show \citep{tarantola:05} that the likelihood of $\m$ is

\begin{equation}
  \mathcal{L} (\m | \dobs) = \textrm{const} \times \exp \left( - \frac{1}{2}\left( \dobs - \dsim_{\textrm{h}}(\m)\right)\trp \Ce\inv  \left( \dobs - \dsim_{\textrm{h}}(\m)\right) \right),
\end{equation}
where
\begin{equation}\label{Eq:Ce}
  \Ce = \textrm{cov} \left[ \e_{\textrm{d}} \right] + \textrm{cov} \left[\e_{\textrm{m}}\right]
\end{equation}
is the total data-error covariance matrix. Note that we use the term ``model error'' to refer to any imperfection in the model used to represent the real reservoir. Examples of sources of model errors include numerical and discretization errors, simplifications of the physics and insufficient parameterization. The assumption in Eq.~\ref{Eq:dtrue} is that all sources of model errors can be aggregated into a random vector with $\textrm{E}\left[\e_{\textrm{m}}\right] = \mathbf{0}$ and known covariance. Even though these are strong assumptions which are unlike to hold in reality, it is important to note that including $\textrm{cov} \left[\e_{\textrm{m}}\right]$ in the matrix $\Ce$ reduces the weights attributed to data helping to partially compensate for deficiencies in the models; see, for example, \citep{sun:17a,oliver:18b}.

If the prior model follows a multivariate Gaussian distribution, $\m \sim \mathcal{N} (\mpr, \Cm)$, then the posterior probability density function (PDF) of $\m$ given $\dobs$ has the form

\begin{eqnarray}
  \nonumber p(\m | \dobs) & = &  \textrm{const} \times \mathcal{L} (\m | \dobs) p(\m) \\
  & = &  \textrm{const} \times \exp \left( - \mathcal{O}(\m) \right),
\end{eqnarray}
where

\begin{equation}\label{Eq:O}
  \mathcal{O}(\m) = \frac{1}{2}\left( \dobs - \dsim_{\textrm{h}}(\m)\right)\trp \Ce\inv  \left( \dobs - \dsim_{\textrm{h}}(\m)\right) + \frac{1}{2}\left( \m - \mpr \right)\trp \Cm\inv  \left( \m - \mpr \right).
\end{equation}

The model $\m$ that minimizes $\mathcal{O}(\m)$ corresponds to the maximum a posteriori \citep{tarantola:05}. In practice, however, we are interested in sampling the posterior PDF to quantify uncertainty. In this case, one alternative is the RML method \citep{oliver:08bk}, which provides an approximate sampling of $p(\m | \dobs)$. Each RML sample is obtained by minimizing a modified version of Eq.~\ref{Eq:O} given by

\begin{equation}\label{Eq:Orml}
  \mathcal{O}_{\textrm{rml}}(\m) = \frac{1}{2}\left( \dobs^\ast - \dsim_{\textrm{h}}(\m)\right)\trp \Ce\inv  \left( \dobs^\ast - \dsim_{\textrm{h}}(\m)\right) + \frac{1}{2}\left( \m - \m^\ast \right)\trp \Cm\inv  \left( \m - \m^\ast \right),
\end{equation}
where $\dobs^\ast \sim \mathcal{N}(\dobs, \Ce)$ and $\m^\ast \sim \mathcal{N}(\mpr, \Cm)$.

\subsection{Data-Space Inversion}
\label{Sec:DSI}

In this section, we review the data-space inversion (DSI) procedure as proposed in \citep{sun:17b} and later improved in \citep{sun:17c}. The main idea behind the method is to use PCA to write the vector of predicted data as

\begin{equation}\label{Eq:PCA}
  \dsim_{\textrm{pca}} = \overline{\dsim} + \Cd\sqr \x,
\end{equation}
where $\overline{\dsim}$ and $\Cd$ are the mean and covariance of $\dsim$, respectively. Both are computed using a prior ensemble, that is,

\begin{equation}
  \overline{\dsim} = \frac{1}{N_e}\sum_{j=1}^{N_e} \dsim_j
\end{equation}
and
\begin{eqnarray}
  \nonumber \Cd & = & \frac{1}{N_e - 1}\sum_{j=1}^{N_e} \left(\dsim_j -  \overline{\dsim}\right)\left(\dsim_j -  \overline{\dsim}\right)\trp \\
  & = & \DD \DD\trp,
\end{eqnarray}
where
\begin{equation}\label{Eq:DD}
  \DD = \frac{1}{\sqrt{N_e - 1}}\left[\dsim_1 -  \overline{\dsim}, \ldots, \dsim_{N_e} -  \overline{\dsim} \right].
\end{equation}
The square root of $\Cd$ in Eq.~\ref{Eq:PCA} is computed using the singular value decomposition (SVD) of $\DD$

\begin{equation}\label{Eq:SVD}
  \DD = \mathbf{U} \bm{\Sigma} \mathbf{V}\trp,
\end{equation}
where $\mathbf{U}$ is a $N_d \times N_d$ orthogonal matrix containing the left singular values of $\DD$, which are equivalent to the left eigenvectors of $\Cd$; $\bm{\Sigma}$ is a $N_d \times N_e$ matrix containing as non-zero elements the singular values of $\DD$, or, equivalently, the square root of the eigenvalues of $\Cd$. The matrix $\mathbf{V}$ contains the right singular values of $\DD$. The square-root of $\Cd$ becomes

\begin{equation}
  \Cd\sqr =  \mathbf{U}\bm{\Sigma}.
\end{equation}

The vector $\x$ in Eq.~\ref{Eq:PCA} is a sample from a standard normal distribution, that is, $\x \sim \mathcal{N}(\mathbf{0}, \mathbf{I})$. In practice, we truncate small singular values using an energy criterium, which means that we consider only the $N_r \leq \max\{N_d, N_e -1\}$ largest singular values such that

\begin{equation}\label{Eq:SVDEnergy}
   \frac{\sum_{i=1}^{N_r} \sigma_i}{\sum_{i=1}^{\min\{N_d, N_e - 1\}} \sigma_i} \geq \xi
\end{equation}
where $\xi \leq 1$ is the energy threshold, typically selected between 0.9 and 0.99 and $\sigma_i$ is the $i$th singular value of $\DD$. Truncating small singular values has two positive side effects: it reduces the dimension of $\x$ and introduces regularization in the inversion. Both are important because the DSI method uses RML for sampling, which requires solving several minimization problems.

\begin{figure}
\centering
	\includegraphics[width=0.4\linewidth]{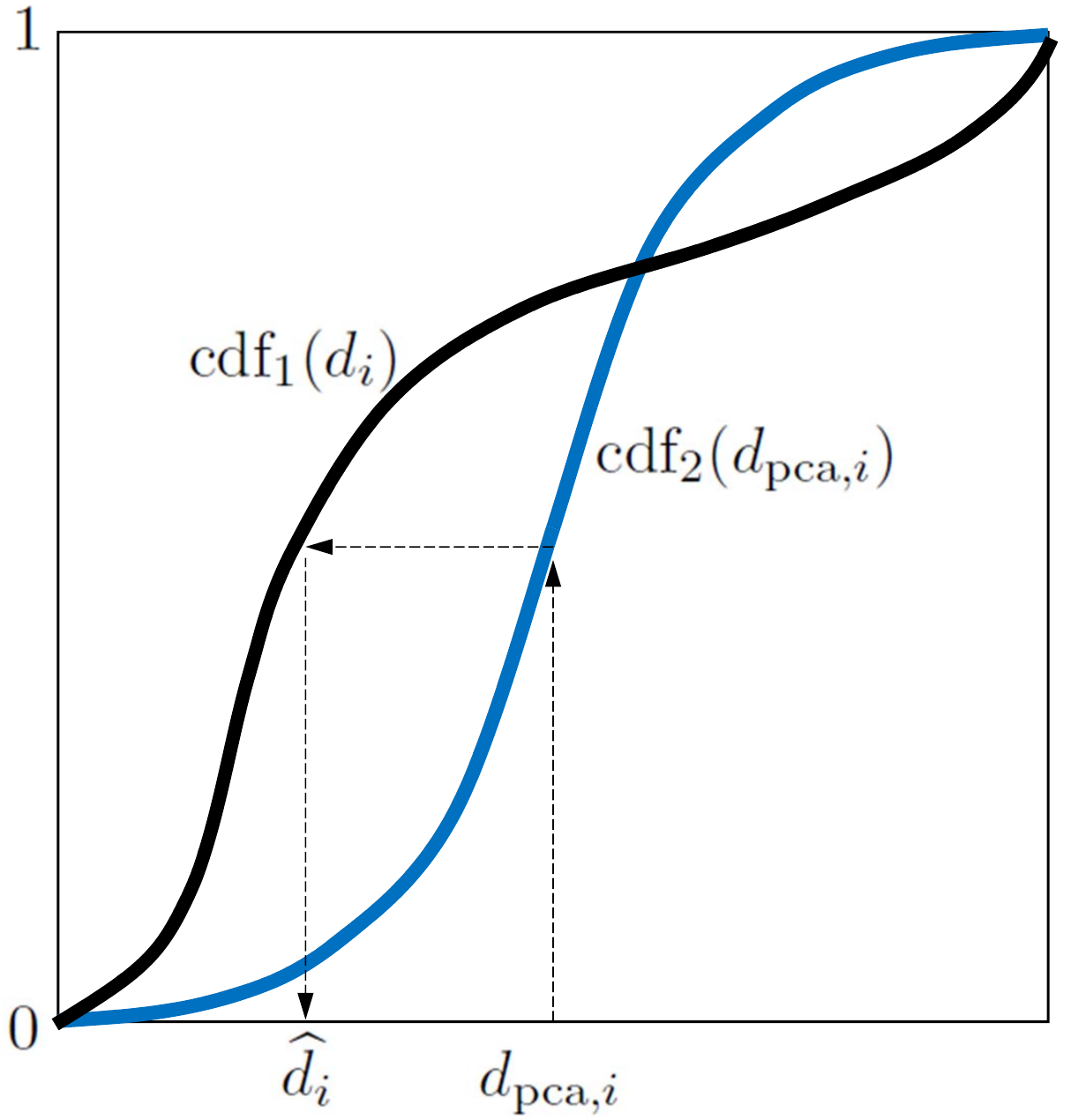}
\caption{Histogram transformation used in DSI.}
\label{Fig:CDFTrans}
\end{figure}

\citet{sun:17b} noted that the direct application of Eq.~\ref{Eq:PCA} may result in nonphysical values for the predicted data, for example negative production or pressure. According to the authors, this problem occurs mainly before water breakthrough time. Therefore, they proposed to apply a data transformation to the prior realizations of $\dsim$ before PCA. The transformation is based on shifting and compressing/stretching the time series. They claim and illustrate in a example that the transformed vectors, $\widehat{\dsim}$, have a more Gaussian prior distribution. However, this procedure is difficult to apply in cases with frequent changes in well controls. In \citep{sun:17c}, the authors proposed to use an inverse Gaussian anamorphosis procedure using the empirical cumulative density function (CDF) computed using the $N_e$ prior realizations of $\dsim$. Figure~\ref{Fig:CDFTrans} illustrates the process, where each component of the transformed vector $\widehat{\dsim}$ is computed as

\begin{equation}\label{Eq:CDFTrans}
  \widehat{d}_i = \textrm{cdf}_1\inv \left( \textrm{cdf}_2 \left(d_{\textrm{pca},i} \right) \right),
\end{equation}
where $\textrm{cdf}_1(\cdot)$ and $\textrm{cdf}_2(\cdot)$ are the CDF's of $\dsim$ and $\dsim_{\textrm{pca}}$, respectively.

The final step of DSI is to use RML to generate a posterior ensemble of predicted data. The RML objective function can be written in terms of the vector of PCA coefficients as

\begin{equation}\label{Eq:Oxrml}
  \mathcal{O}_{\textrm{rml}}(\x) = \frac{1}{2}\left( \dobs^\ast - \widehat{\dsim}_{\textrm{h}}(\x)\right)\trp \Ce\inv  \left( \dobs^\ast - \widehat{\dsim}_{\textrm{h}}(\x) \right) + \frac{1}{2}\left( \x - \x^\ast \right)\trp \left( \x - \x^\ast \right),
\end{equation}
where $\x^\ast \sim \mathcal{N}(\mathbf{0}, \mathbf{I})$. The minimization of Eq.~\ref{Eq:Oxrml} can be done with any optimization method. \citet{sun:17b} used the Broyden–Fletcher–Goldfarb–Shanno (BFGS) method \citep{nocedal:06}, in which case it is necessary to compute the gradient of $\mathcal{O}_{\textrm{rml}}(\x)$ with respect to $\x$. However, $\widehat{\dsim}_{\textrm{h}}(\x)$ is a nonlinear function of $\x$ with no general analytical form because of the transformation of Eq.~\ref{Eq:CDFTrans}. One alternative, is to compute numerical gradients, which may be computationally expensive if we have a large number of data points. \citet{jiang:18a} proposed to ignore the transformation of Eq.~\ref{Eq:CDFTrans} and use the analytical gradients computed with $\dsim_{\textrm{pca}}$ instead of $\widehat{\dsim}$. Our limited set of tests indicated that this procedure in fact improved the computational performance of our DSI implementation, which uses the limited-memory BFGS \citep{nocedal:06} implementation available in the C\# library Accord.NET \citep{souza:17a}.

\subsection{Data-Space Inversion with Ensemble Smoother}
\label{Sec:DSI-ESMDA}

The ensemble smoother (ES) was introduced by \citet{vanleeuwen:96} as an alternative to the sequential data assimilation scheme of EnKF. The first application of ES for history matching was presented by \citet{skjervheim:11a}, which concluded that the method is faster than EnKF with similar results. Despite the good results presented in \citep{skjervheim:11a}, some authors \citep{chen:12a,chen:13a,emerick:13b,emerick:13a} observed that ES tends to result in unreasonable data matches when applied to more complex history-matching problems. The reason for the poor performance is because ES is similar to applying a single Gauss-Newton iteration to minimize a RML-type of objective function with a sensitivity matrix estimated based on the prior ensemble \citep{reynolds:06}. This fact lead to the development of several iterative forms of ES; see, for example, \citep{chen:12a,chen:13a,emerick:13b,stordal:14b,stordal:15a,luo:15b}. Among the iterative forms of ES, the ensemble smoother with multiple data assimilation (ES-MDA) \citep{emerick:13b} is a popular choice. The popularity of ES-MDA can be attributed mainly to its good performance in history-matching problems \citep{emerick:13c,emerick:16a,maucec:16a,evensen:18a} and its simplicity of implementation. In fact, ES-MDA is essentially equivalent to repeat ES a few times with the data-error covariance matrix, $\Ce$, multiplied by coefficients $\alpha_k$'s to avoid overweighing the measurements. The choice of the coefficients and the number of repetitions of ES, $N_a$, must obey the condition $\sum_{k=1}^{N_a} \alpha_k\inv = 1$ to ensure that ES-MDA samples the correct posterior PDF in the linear-Gaussian case \citep{emerick:13b}. For nonlinear problems, the choice of the $\alpha_k$'s and $N_a$ has a major impact in the performance of the method. There are recent works proposing methods to select $\alpha_k$'s and $N_a$ \citep{duc:16a,emerick:16a,rafiee:17a,ma:17a,emerick:18b}. However, here we use simplest choice, which consists of selecting $N_a$ in advance and setting $\alpha_k = N_a$, for $k=1, \ldots, N_a$.

The application of ES-MDA to DSI, that is, to generate samples of the predicted data vector, $\dsim = [\dsim_{\textrm{h}}\trp, \dsim_{\textrm{f}}\trp]\trp$, given the vector of observations, $\dobs$, is straightforward. We can write the resulting DSI-ESMDA update equation as

\begin{equation}\label{Eq:DSI-ESMDA}
  \dsim_j^{k+1} = \dsim_j^{k} + \R \circ \K^k \left( \dobs + \sqrt{\alpha_k}\e_j^k - \dsim_{\textrm{h}}^k  \right),
\end{equation}
for $k = 1, \ldots, N_a$ and $j = 1, \ldots, N_e$, where $\K$ is a modified version of the Kalman gain given by

\begin{equation}\label{Eq:DSI-ESMDA-Kgain}
  \K^k = \DD^k \left(\DD^k_{\textrm{h}}\right)\trp \left( \DD_{\textrm{h}}^k \left(\DD^k_{\textrm{h}}\right)\trp +  \alpha_k \Ce \right)\inv.
\end{equation}
In the above equations, $\DD$ was defined before (Eq.~\ref{Eq:DD}) and $\DD_{\textrm{h}}$ includes only the predicted data corresponding to the historical period, that is,
\begin{equation}
  \DD_{\textrm{h}} = \frac{1}{\sqrt{N_e - 1}}\left[\dsim_{\textrm{h},1} -  \overline{\dsim_{\textrm{h}}}, \ldots, \dsim_{\textrm{h},N_e} -  \overline{\dsim_{\textrm{h}}} \right].
\end{equation}
The vector $\e$ is a sample from $\mathcal{N}(\mathbf{0},\Ce)$ and $\R$ is the localization matrix with ``$\circ$'' denoting the Schur (element-wise) product. The matrix inversion required in Eq.~\ref{Eq:DSI-ESMDA-Kgain} is done using the subspace inversion method \citep{evensen:04} as described in the Appendix section of \citep{emerick:12a}. This procedure involves the truncated SVD of a rescaled matrix $\Ce\invsqr \DD_{\textrm{h}}$ to avoid loss of relevant information when removing small singular values. The number of singular values retained is computed using the same energy criterium of Eq.~\ref{Eq:SVDEnergy}.

The reason for introducing the Schur product between the localization matrix and the Kalman gain in Eq.~\ref{Eq:DSI-ESMDA} is twofold: it regularizes the estimates of the Kalman gain removing long-distance spurious correlations and increases the degrees of freedom to assimilate data \citep{houtekamer:01}. In fact, localization proved to improve substantially the results of ensemble methods applied to history matching; see, for example, \citep{chen:10b,emerick:11e,emerick:11c}.

The regularization introduced by the Schur product is obtained by constructing the localization matrix using a correlation function with compact support. A common choice is the fifth-order compact correlation function from \citep{gaspari:99}, in which case each entry of the matrix $\R$ is computed using

\begin{equation}\label{Eq:GaspariCohn}
r\left( \frac{h}{L} \right) =\left\{ \begin{array}{ll}
-\frac{1}{4}\left(\frac{h}{L}\right)^{5}+\frac{1}{2}\left(\frac{h}{L}\right)^{4}+\frac{5}{8}\left(\frac{h}{L}\right)^{3}-\frac{5}{3}\left(\frac{h}{L}\right)^{2}+1, & \text{if}~0 \leq \frac{h}{L} \leq 1\\
\frac{1}{12}\left(\frac{h}{L}\right)^{5}-\frac{1}{2}\left(\frac{h}{L}\right)^{4}+\frac{5}{8}\left(\frac{h}{L}\right)^{3}+\frac{5}{3}\left(\frac{h}{L}\right)^{2}-5\left(\frac{h}{L}\right)+4-\frac{2}{3}\left(\frac{h}{L}\right)^{-1}, & \text{if}~1 \leq \frac{h}{L} \leq 2\\
0 & \text{if}~\frac{h}{L} > 2
\end{array}\right.,
\end{equation}
where $h$ is a ``distance'' and $L$ is a parameter known as ``critical length,'' which corresponds to the distance where the correlation function decays to approximately 0.21. Note that the matrix $\R$ has the same shape of the Kalman gain and that each row of the Kalman gain corresponds to the correction term applied to each variable updated with Eq.~\ref{Eq:DSI-ESMDA}, that is, entries of the vector $\dsim$. Each column of the Kalman gain corresponds to an observation used in the conditioning process. Therefore, $h$ in Eq.~\ref{Eq:GaspariCohn} corresponds to the distance between entries of $\dsim$ and $\dobs$. In the applications of interest of this paper, both $\dsim$ and $\dobs$ contain data from wells in the field. Therefore, the ratio $h/L$ is computed based on the spatial distance between wells. Besides the spatial distance, we also introduced the time difference between data points in the calculation of $h/L$ to account for the fact that early data have lower correlations with the predictions than the late data. The final expression for computing $h/L$ is given by

\begin{equation} \label{Eq:hL}
  \frac{h}{L} = \sqrt{\left(\frac{\Delta x^\prime}{L_{x^\prime}}\right)^{2}+\left(\frac{\Delta y^\prime}{L_{x^\prime}}\right)^{2} + \left(\frac{\Delta t}{T}\right)^{2}},
\end{equation}
where

\begin{equation}\label{Eq:Aniso}
  \left[\begin{array}{c}
\Delta x^\prime \\
\Delta y^\prime
\end{array}\right] = \left[\begin{array}{cc}
\cos \theta & -\sin \theta\\
\sin\theta & \cos \theta
\end{array}\right]\left[\begin{array}{c}
\Delta x \\
\Delta y
\end{array}\right].
\end{equation}
In the above equations, $\Delta x$ and $\Delta y$ correspond to the spatial distance between the wells in the orthogonal directions $x$ and $y$, while $\Delta x^\prime$ and $\Delta y^\prime$ are the corresponding distances rotated by an angle $\theta$. This allows to consider anisotropy in the localization function as discussed in \citep{emerick:11e}. $\Delta t$ is the time difference between data points. $L_{x^\prime}$, $L_{y^\prime}$ and $T$ are the critical lengths in the ``directions'' $x^\prime$, $y^\prime$ and $t$, respectively.

Besides removing long-distance spurious correlations, localization also increases the degrees of freedom to assimilate data. To show that, first consider the case without localization. Using Eq.~\ref{Eq:DSI-ESMDA-Kgain} in (\ref{Eq:DSI-ESMDA}), we can define the vector $\wj^k$ as

\begin{equation}
  \wj^k \equiv \left( \DD_{\textrm{h}}^k \left(\DD^k_{\textrm{h}}\right)\trp + \alpha_k \Ce \right)\inv \left( \dobs + \sqrt{\alpha_k}\e_j^k - \dsim_{\textrm{h}}^k  \right).
\end{equation}
Hence

\begin{eqnarray}\label{Eq:DSI-ESMDA-nolocal}
  \nonumber \dsim_j^{k+1} & = & \dsim_j^{k} + \DD^k \left(\DD^k_{\textrm{h}}\right)\trp \wj^k \\
  & = & \dsim_j^{k} + \frac{1}{N_e - 1} \sum_{i=1}^{N_e} \left(\dsim_i^k - \overline{\dsim^k} \right) \left(\dsim_{\textrm{h},i}^k - \overline{\dsim_{\textrm{h}}^k} \right)\trp \wj^k \\
  \nonumber & = & \dsim_j^{k} + \frac{1}{N_e - 1}\sum_{i=1}^{N_e} \beta_{i,j}^k \left(\dsim_i^k - \overline{\dsim^k} \right),
\end{eqnarray}
where $\beta_{i,j}^k \equiv \left(\dsim_{\textrm{h},i}^k - \overline{\dsim_{\textrm{h}}^k} \right)\trp \wj^k$ is a scalar. The vector $\overline{\dsim^k}$ is a linear combination of the vectors $\dsim_i^k$'s. Hence, it is possible to write $\dsim_j^{k+1}$ as

\begin{equation}\label{Eq:LinearComb}
  \dsim_j^{k+1} = \sum_{i=1}^{N_e} \gamma_{i,j}^k \dsim_i^{k},
\end{equation}
where $\gamma_{i,j}^k$'s are scalars. Eq.~\ref{Eq:LinearComb} means that $\dsim_j^{k+1}$ is a linear combination of the vectors $\dsim_i^k$'s. Applying this result recursively to all MDA iterations, we conclude that the final predictions from DSI-ESMDA without localization are simply linear combinations of the prior ones. This effectively means that we have at most $N_e$ coefficients (degrees of freedom) available to assimilate data. Note that the number of degrees of freedom of the standard DSI is also limited by the number of PCA coefficients which is $N_r \leq \min \{ N_d, N_e - 1\}$.

Considering now the case with localization, we can write

\begin{eqnarray}\label{Eq:DSI-ESMDA-local1}
  \nonumber \dsim_j^{k+1} & = & \dsim_j^{k} + \left( \R \circ \K^k  \right) \delta\dsim_j^k \\
  & = & \dsim_j^{k} + \sum_{i=1}^{N_{d,\textrm{h}}} \left(\mathbf{r}_i \circ \bm{\kappa}^k_i \right) \delta d_{i,j}^k,
\end{eqnarray}
where $\delta\dsim_j^k \equiv  \dobs + \sqrt{\alpha_k}\e_j^k - \dsim_{\textrm{h}}^k$ and $\delta d_{i,j}^k$ is its $i$th entry. $\mathbf{r}_i$ and $\bm{\kappa}^k_i$ correspond to the $i$th columns of the localization and Kalman gain matrices, respectively. Writing the update equation for the $n$th entry of the vector $\dsim_j^{k}$, we have

\begin{eqnarray}\label{Eq:DSI-ESMDA-local2}
  \nonumber d_{n,j}^{k+1} & = & d_{n,j}^k + \sum_{i=1}^{N_{d,\textrm{h}}} \kappa^k_{n,i}  r_{n,i} \delta d_{i,j}^k \\
  & = & d_{n,j}^k + \sum_{i=1}^{N_{d,\textrm{h}}} \kappa^k_{n,i} \eta_{n,i,j}.
\end{eqnarray}
For each $n$, we have a different coefficient $\eta_{n,i,j}$. Therefore, Eq.~\ref{Eq:DSI-ESMDA-local2} means that each component of $\dsim_j^{k+1}$ may be computed with a different linear combination of the $N_{d,\textrm{h}}$ columns of $\K^k$. Thus, localization expands the degrees of freedom to assimilate data.

The application of DSI-ESMDA is similar to the standard DSI, both methods require to run reservoir simulations only for the prior ensemble of models to generate the prior ensemble of predicted data. After that, DSI-ESMDA applies Eq.~\ref{Eq:DSI-ESMDA} $N_a$ times to generate the posterior ensemble. Note that the $N_a$ iterations are necessary because the prior ensemble is not Gaussian (if the prior is Gaussian, MDA is equivalent to a single ES update). Moreover, note that using Gaussian anamorphosis to transform variables only ensures that the marginal distributions of individual $d_i$'s are Gaussian, the joint distribution which is updated with DSI-ESMDA may not be Gaussian. In our tests, we noticed that only a few MDA iterations are required. In all cases presented in the next section, we use $N_a = 4$. DSI-ESMDA can also be applied to the same parameterization of DSI, that is, update the PCA coefficients with the data transformation of Eq.~\ref{Eq:CDFTrans}. However, our initial tests showed that this procedure did improve the results, actually the results were slightly worse. Moreover, using the PCA parameterization of DSI would prevent to use localization. The only correction we applied is to truncate in zero if the final predicted data is negative, which in our tests occurred with water production data before breakthrough. Note that we apply this truncation only to the final estimates, not during the DSI-ESMDA iterations.

\section{Test Cases}
\label{Sec:TestCases}

\subsection{Test Case 1}
\label{Sec:TestCase1}

\begin{figure}
\centering
	\includegraphics[width=0.6\linewidth]{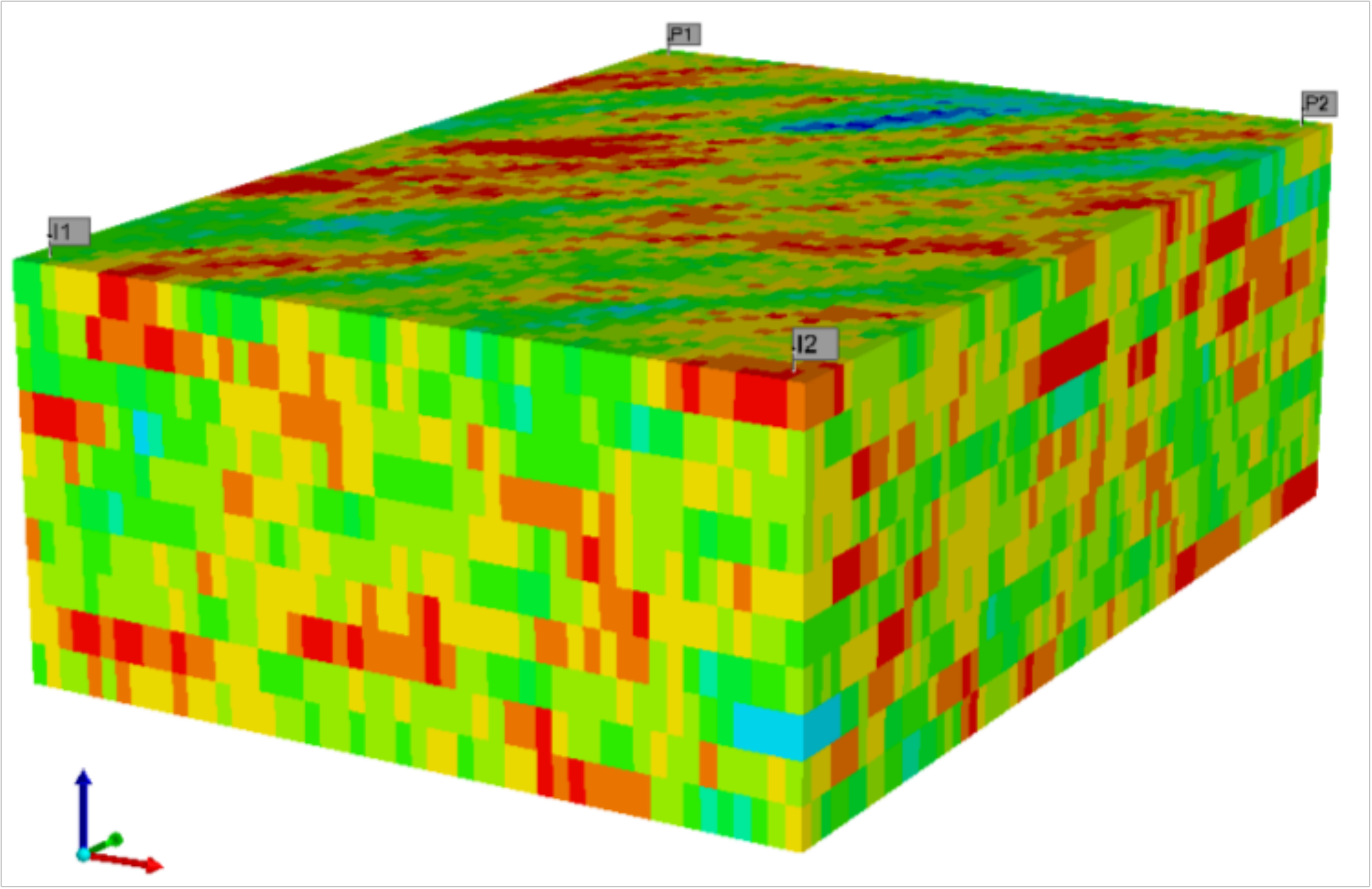}
\caption{Test case 1.}
\label{Fig:Case1Model}
\end{figure}

The first test problem is a synthetic reservoir case created with rock and fluid properties typically found in the Campos Basin, Brazil. The reservoir model contains 50 $\times$ 70 $\times$ 10 gridblocks with uniform size of 50~m $\times$ 50~m $\times$ 5~m. A reference (true) case was generated using sequential Gaussian simulation for modeling porosity and sequential Gaussian co-simulation for permeability using the porosity as secondary variable and a correlation coefficient of 0.95. Table~\ref{Tab:GeoCase1} summarizes the geostatistical parameters used to generate the reference model. The model has two vertical oil producing and two vertical water injection wells placed on the borders of the reservoir as illustrated in Fig.~\ref{Fig:Case1Model}. The producers operated under as specified bottom-hole pressure (BHP) of 25,000~kPa while the injectors operate with a BHP of 35,000~kPa.

\begin{table}
	\centering
	\caption{Geostatistical parameters used to create the reference model and the prior ensemble. Test case 1}
	\label{Tab:GeoCase1}
	\begin{tabular}{lcc}
		\toprule
		Parameter & Porosity & Log-permeability \\
		\midrule
		Mean & 0.22 & 7.2~ln-mD \\
        Standard deviation & 0.05 & 0.60~ln-mD \\
        Variogram type & Spherical & Spherical \\
        Variogram maximum range & 1124.0~m & 1124.0~m \\
        Variogram minimum range & 281.0~m & 281.0~m \\
        Azimuth & 45$^\circ$ & 45$^\circ$ \\
		\bottomrule
	\end{tabular}
\end{table}

Using the same geostatistical parameters of the reference model, we generated an ensemble of 500~realizations for testing DSI and DSI-ESMDA. The synthetic measurements correspond to oil and water rate at both production wells and water injection rate at both injectors. These observations were corrupted by adding random Gaussian noise with zero mean and standard deviation corresponding to 10\% of the data predicted by the reference model. For DSI, we kept 99\% of the singular value energy (Eq.~\ref{Eq:SVDEnergy}) which corresponded to 98~singular values. For DSI-ESMDA, we also kept 99\% of the singular value energy in the subspace inversion. We used $N_a = 4$ MDA iterations. Neither spatial nor temporal localizations were applied in this case. For comparisons, we also applied a standard history matching (that is, model-space inversion) using ES-MDA to update the same prior ensemble of 500 realizations with $N_a = 4$ without localization. It is worth noting that the results of the history matching with ES-MDA do not correspond to the reference solution, as the method is not guaranteed to sample the posterior PDF correctly for nonlinear problems. The objective is to compare the DSI results against a history matching procedure that is used in practice. There are rigorous sampling methods, such as Markov chain Monte Carlo and rejection sampling \citep{oliver:08bk} that could be used to generate the reference solution. Unfortunately, these methods are computationally prohibitive, even for a simple problem such as the one describe in this section. One alternative would be to simplify the problem by reducing the size of the model and the number of data points. However, this has already been done by \citep{sun:17b}, where it is shown that DSI obtained reasonable results compared to rejection sampling for a problem with few data points.

\begin{figure}
\centering
    \captionsetup{justification=centering}
    \subfloat[]{
      \includegraphics[width=0.5\textwidth]{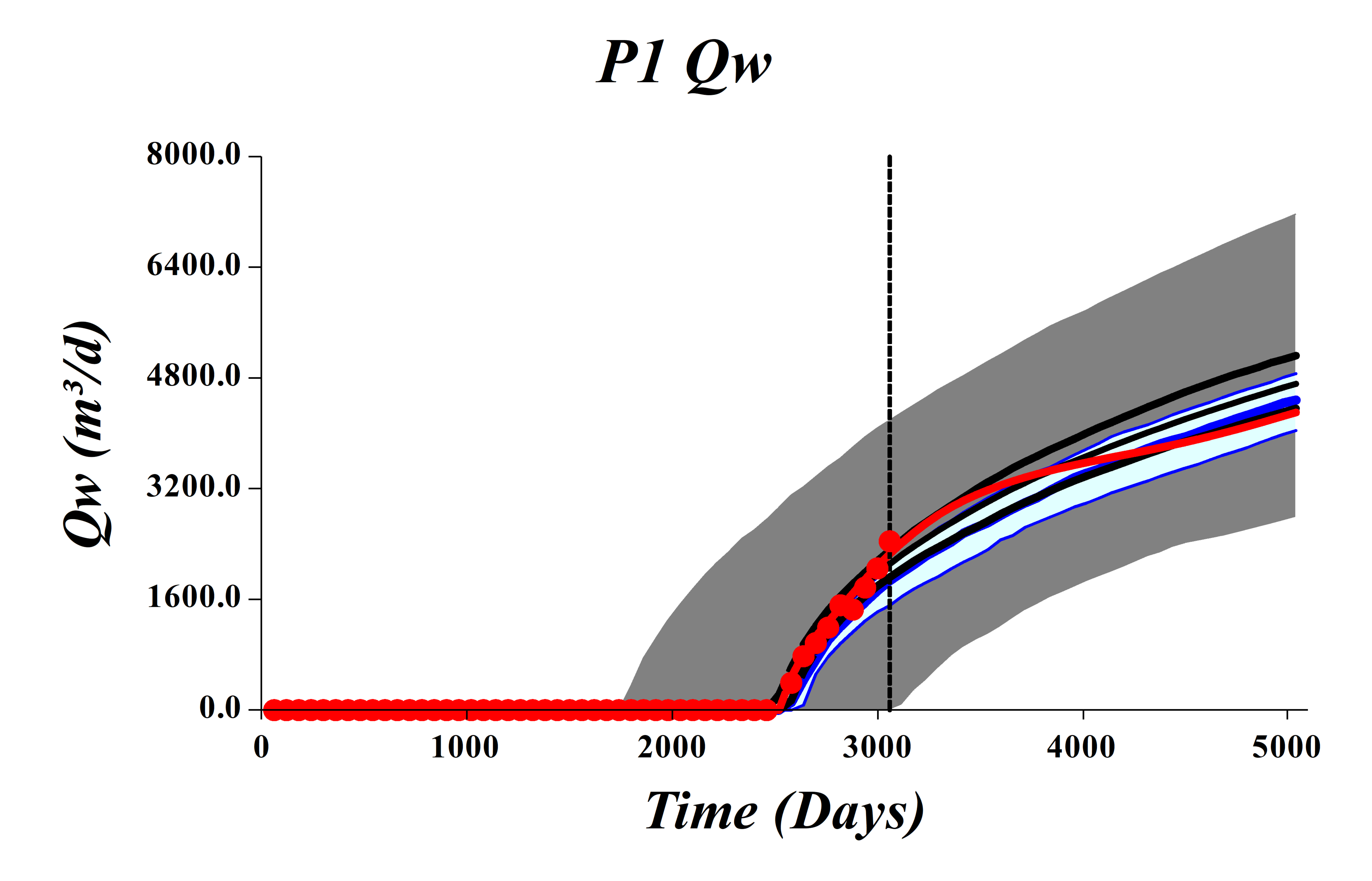}
    }
    \subfloat[]{
      \includegraphics[width=0.5\textwidth]{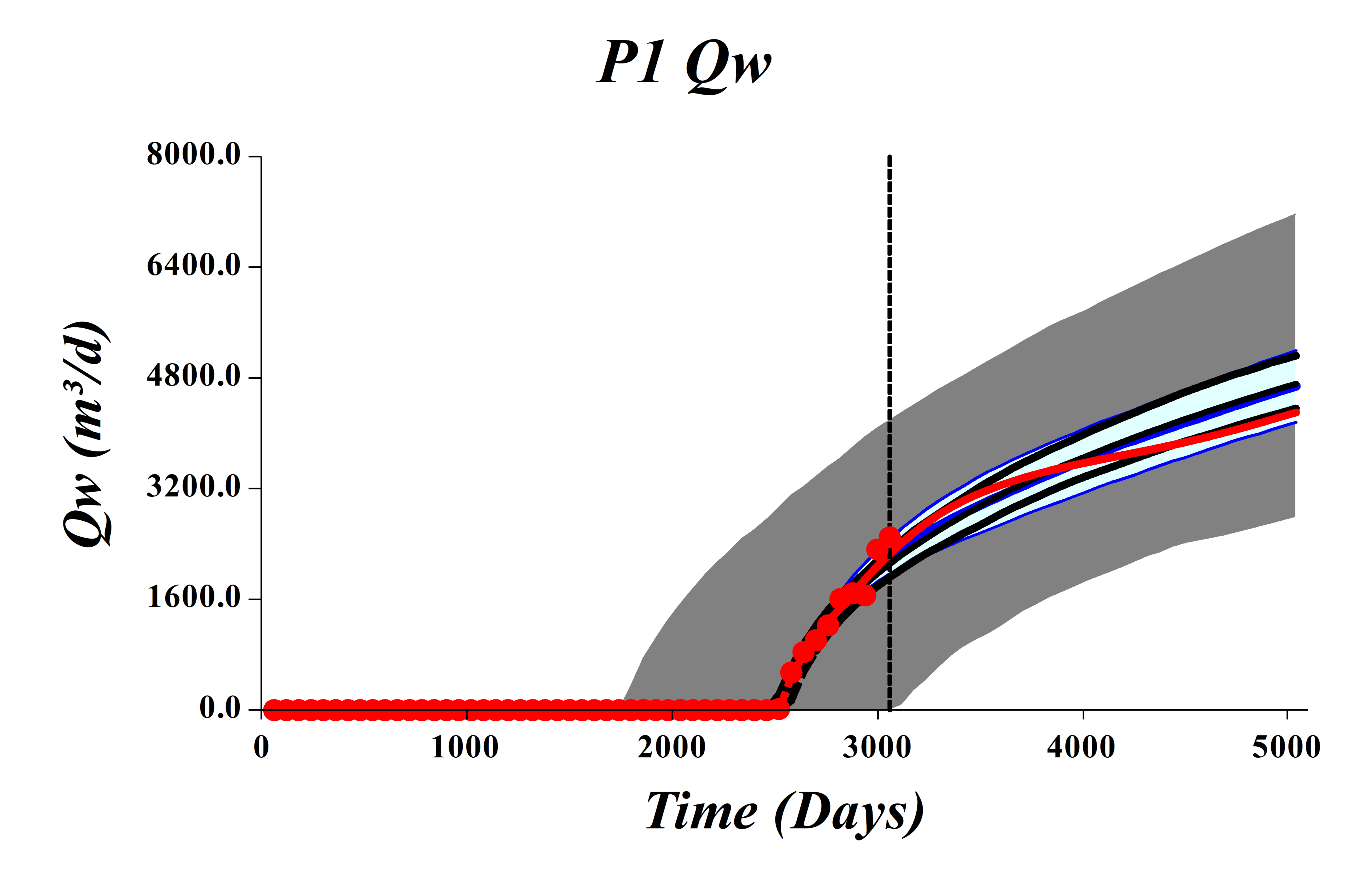}
    }
    \linebreak
    \subfloat[]{
      \includegraphics[width=0.5\textwidth]{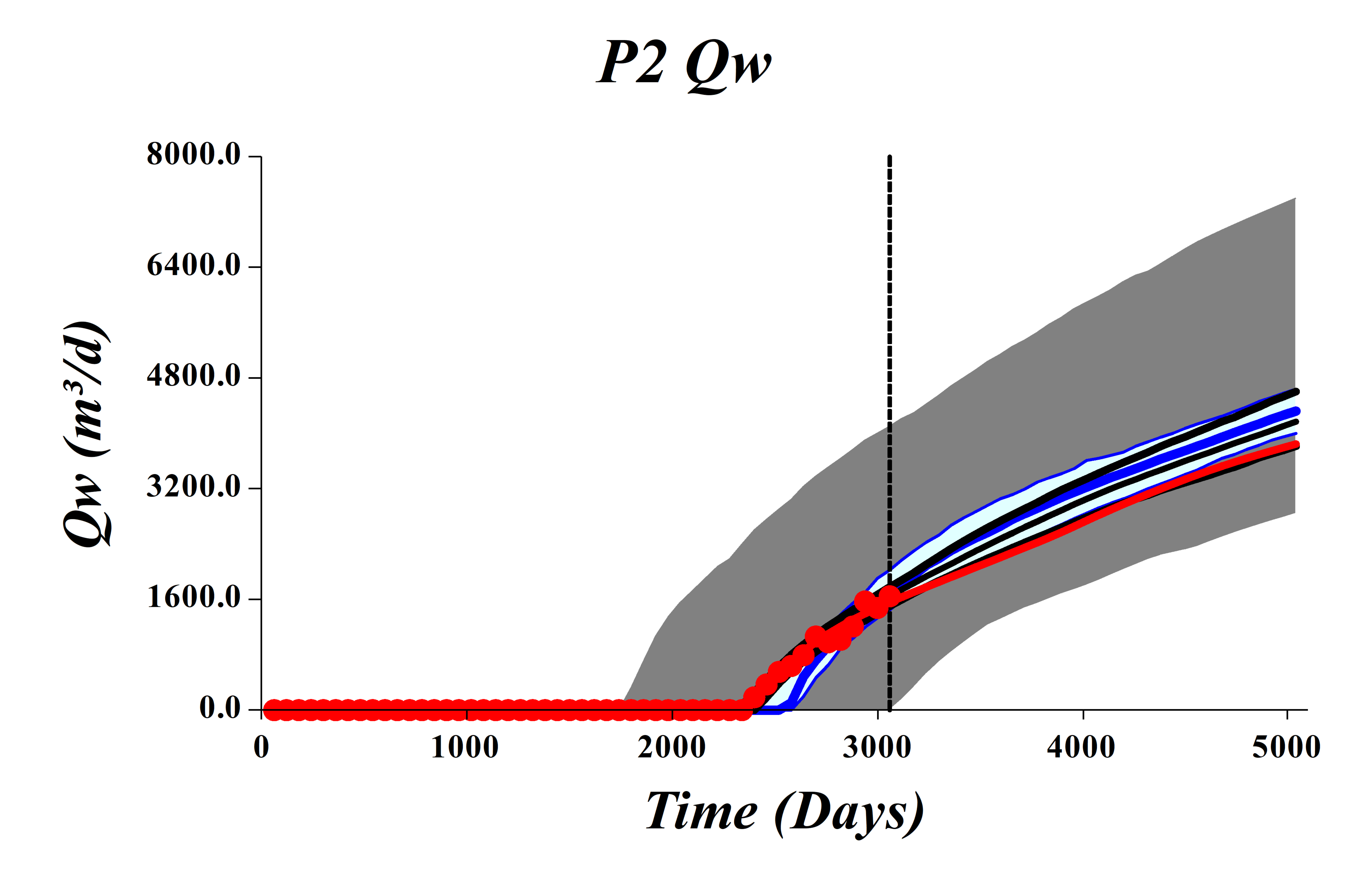}
    }
    \subfloat[]{
      \includegraphics[width=0.5\textwidth]{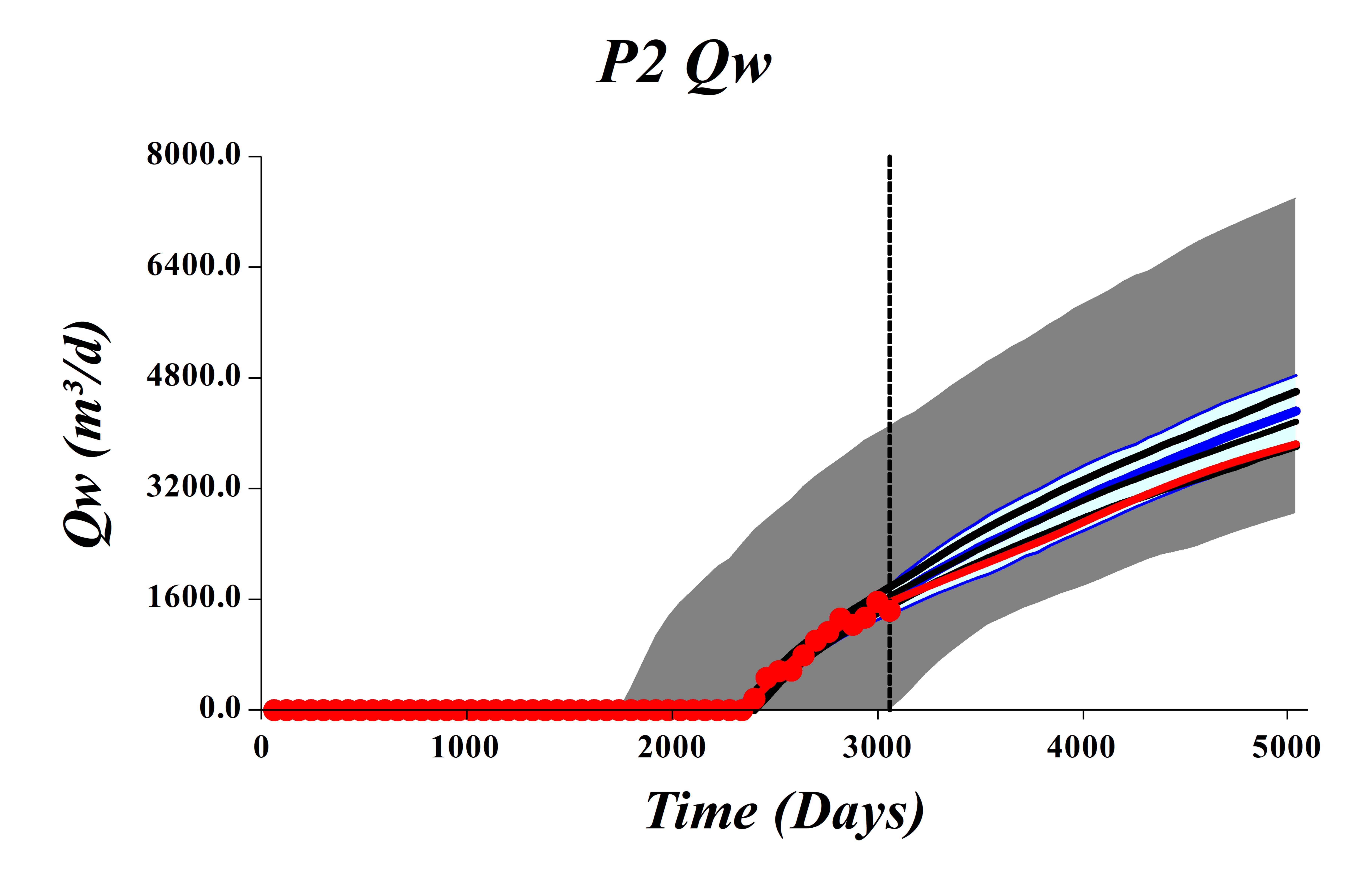}
    }

\captionsetup{justification=justified}
\caption{Water production rate in m$^\text{3}$/days for wells P1 (first row) and P2 (second row). Test case 1. (a) and (c) DSI, (b) and (d) DSI-ESMDA. The red dots are the observed data points and the red line is the prediction from the reference model. The grey region corresponds to the predictions within the percentiles P10--P90 obtained with the prior ensemble. The light blue region corresponds to the predictions within the percentiles P10--P90 obtained with DSI or DSI-ESMDA. The blue line corresponds to the percentile P50 obtained with DSI or DSI-ESMDA. The black lines correspond to the percentiles P10, P50 and P90 obtained by the history-matched models using ES-MDA. The vertical dashed line indicates the end of the history and beginning of the forecast period.}
\label{Fig:WaterRateTestCase1}
\end{figure}

\begin{figure}
\centering
    \captionsetup{justification=centering}
    \subfloat[]{
      \includegraphics[width=0.5\textwidth]{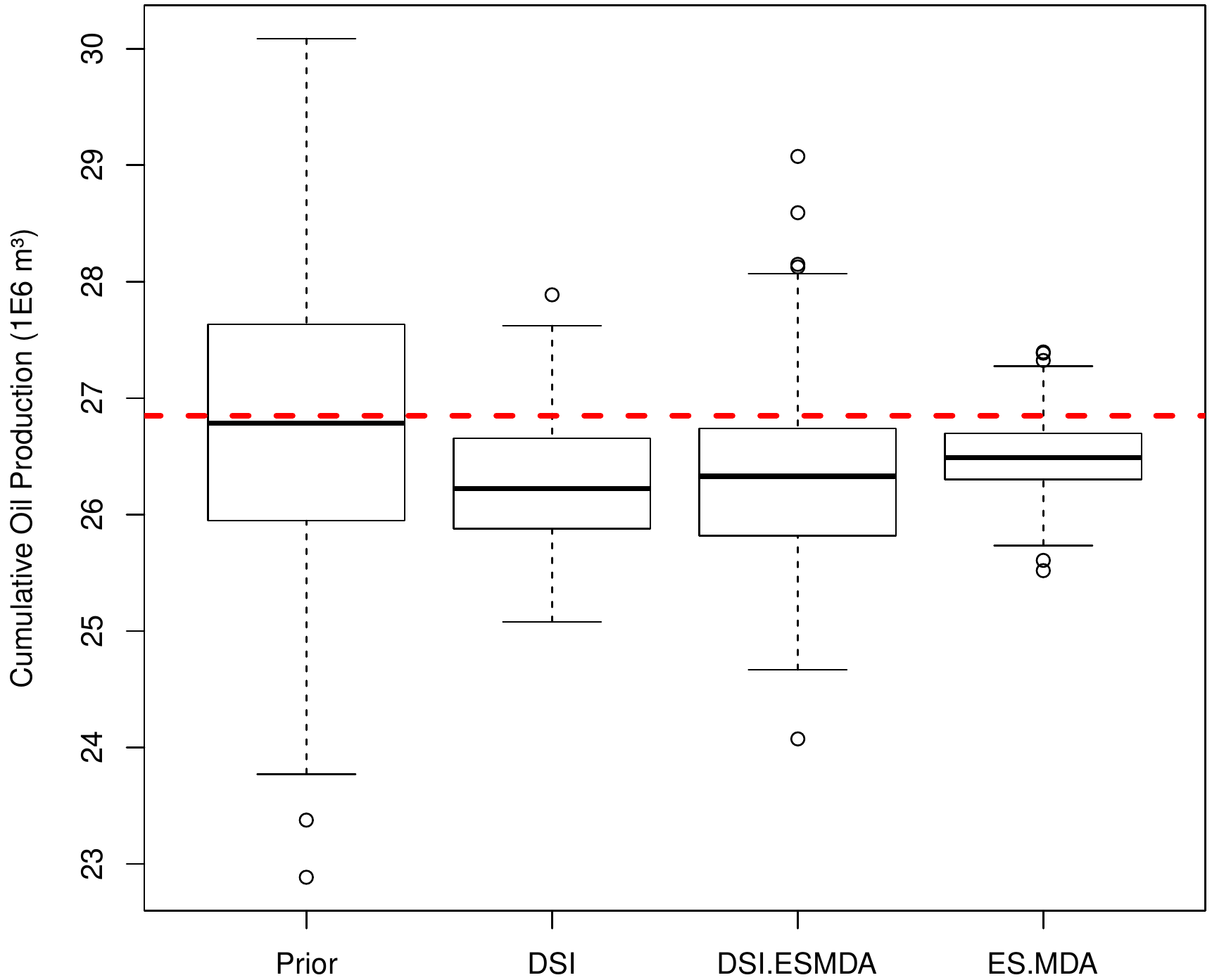}
    }
        \subfloat[]{
      \includegraphics[width=0.5\textwidth]{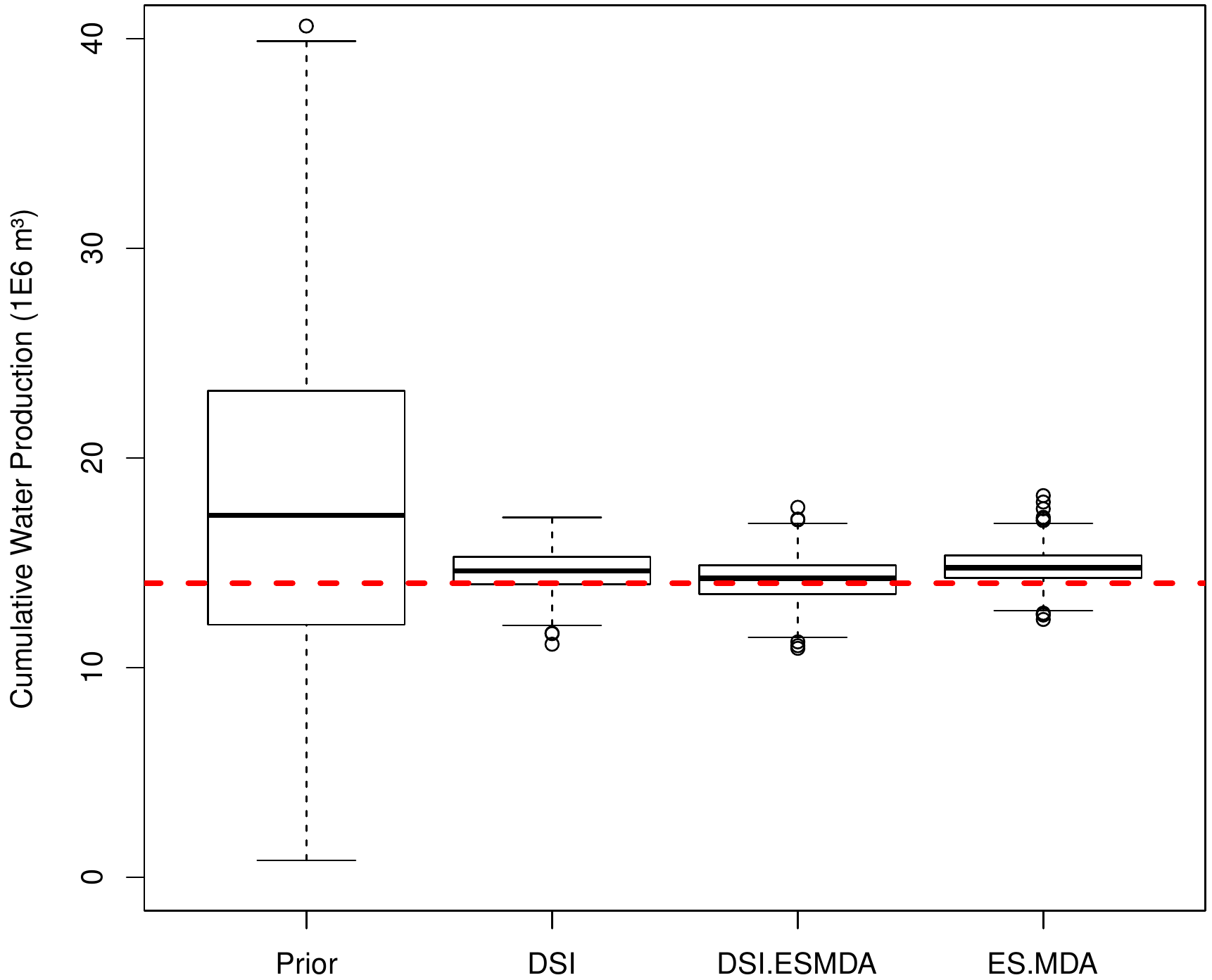}
    }
\captionsetup{justification=justified}
\caption{Field cumulative production in $\text{10}^\text{6}$~m$^\text{3}$. (a) Oil and (b) water. Test case 1. The dashed red line indicates the cumulative production of the reference case.}
\label{Fig:NpWpCase1}
\end{figure}

\begin{figure}
\centering
    \captionsetup{justification=centering}
    \subfloat[]{
      \includegraphics[width=0.5\textwidth]{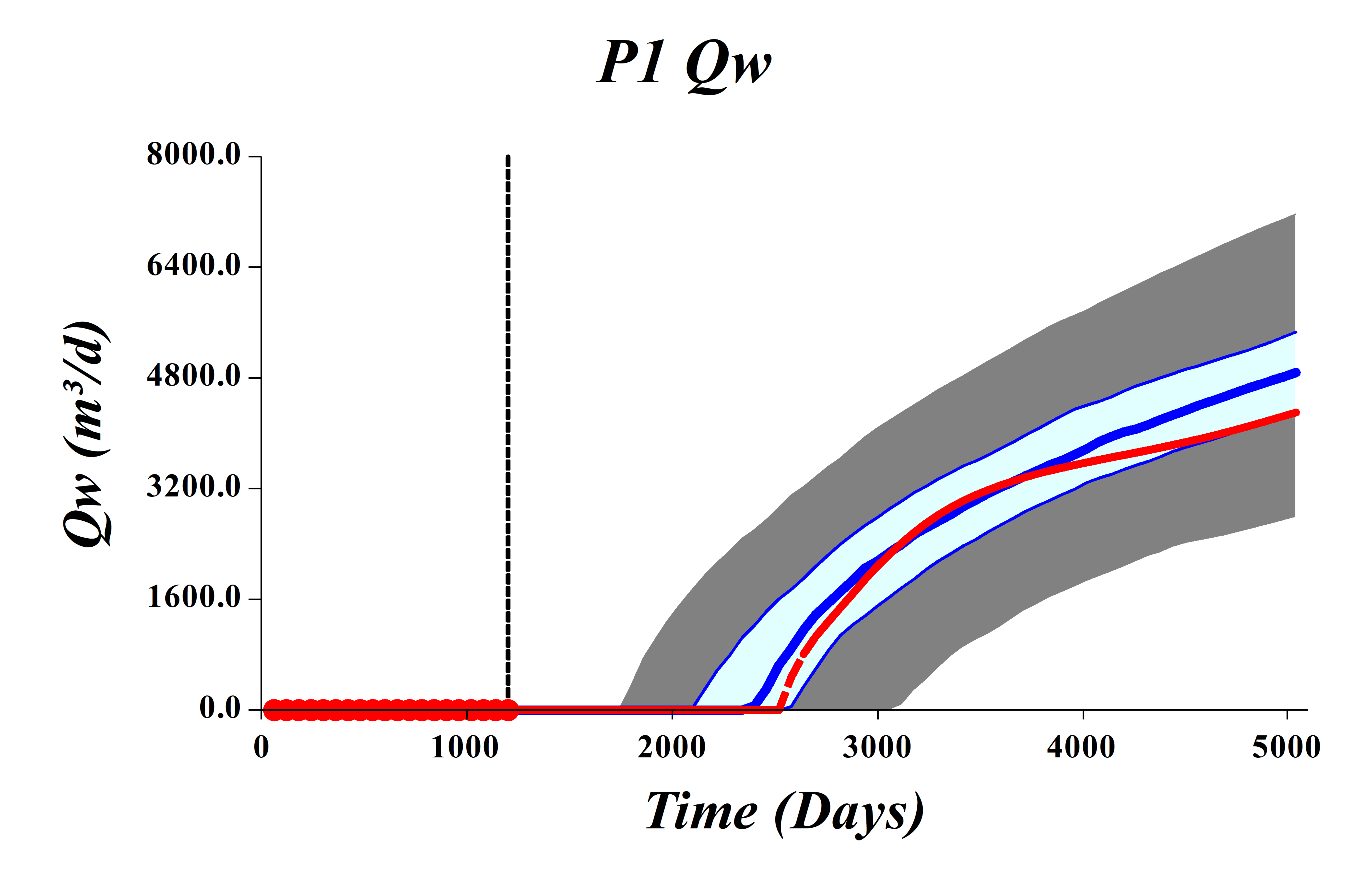}
    }
    \subfloat[]{
      \includegraphics[width=0.5\textwidth]{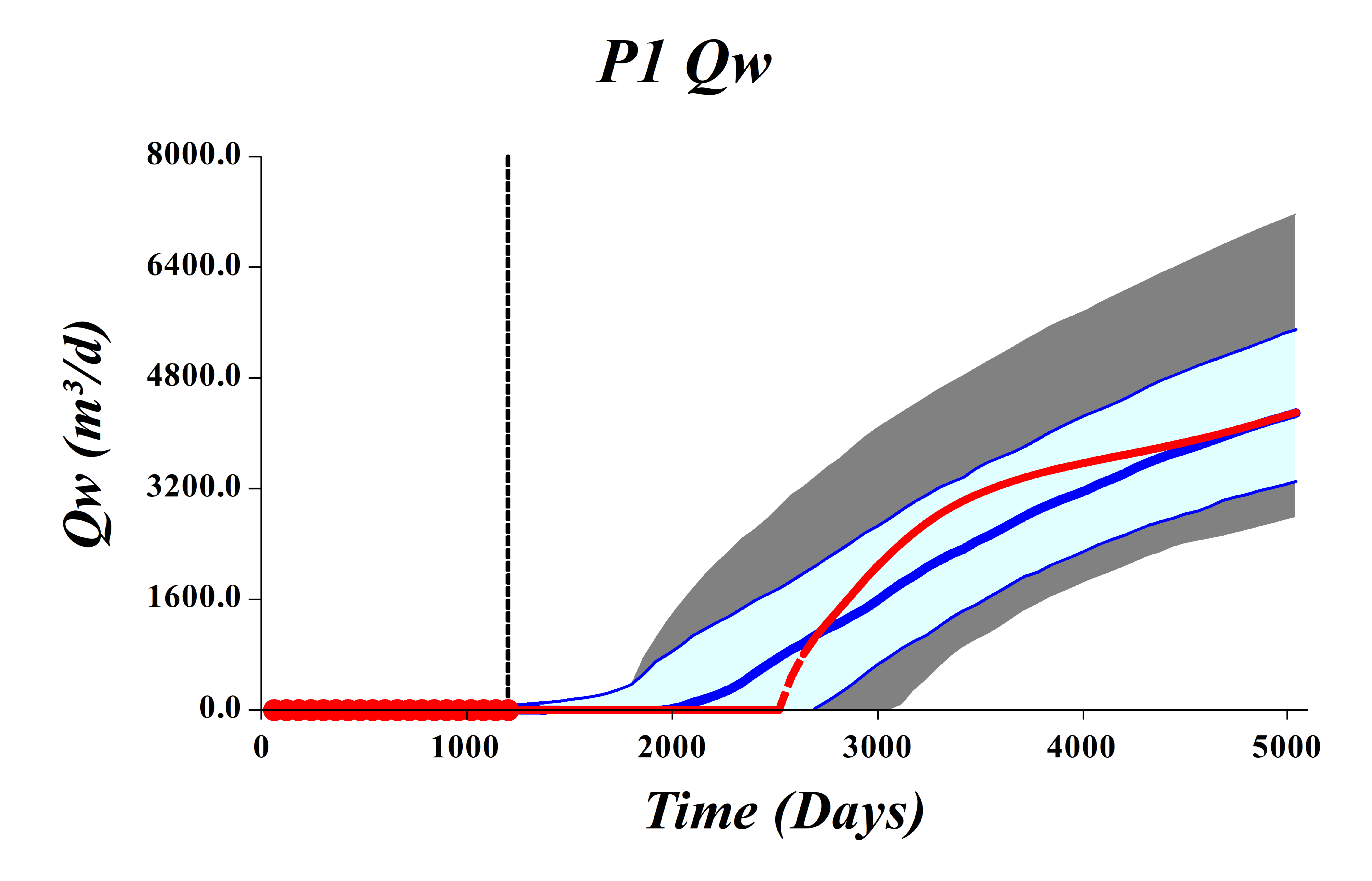}
    }
    \linebreak
    \subfloat[]{
      \includegraphics[width=0.5\textwidth]{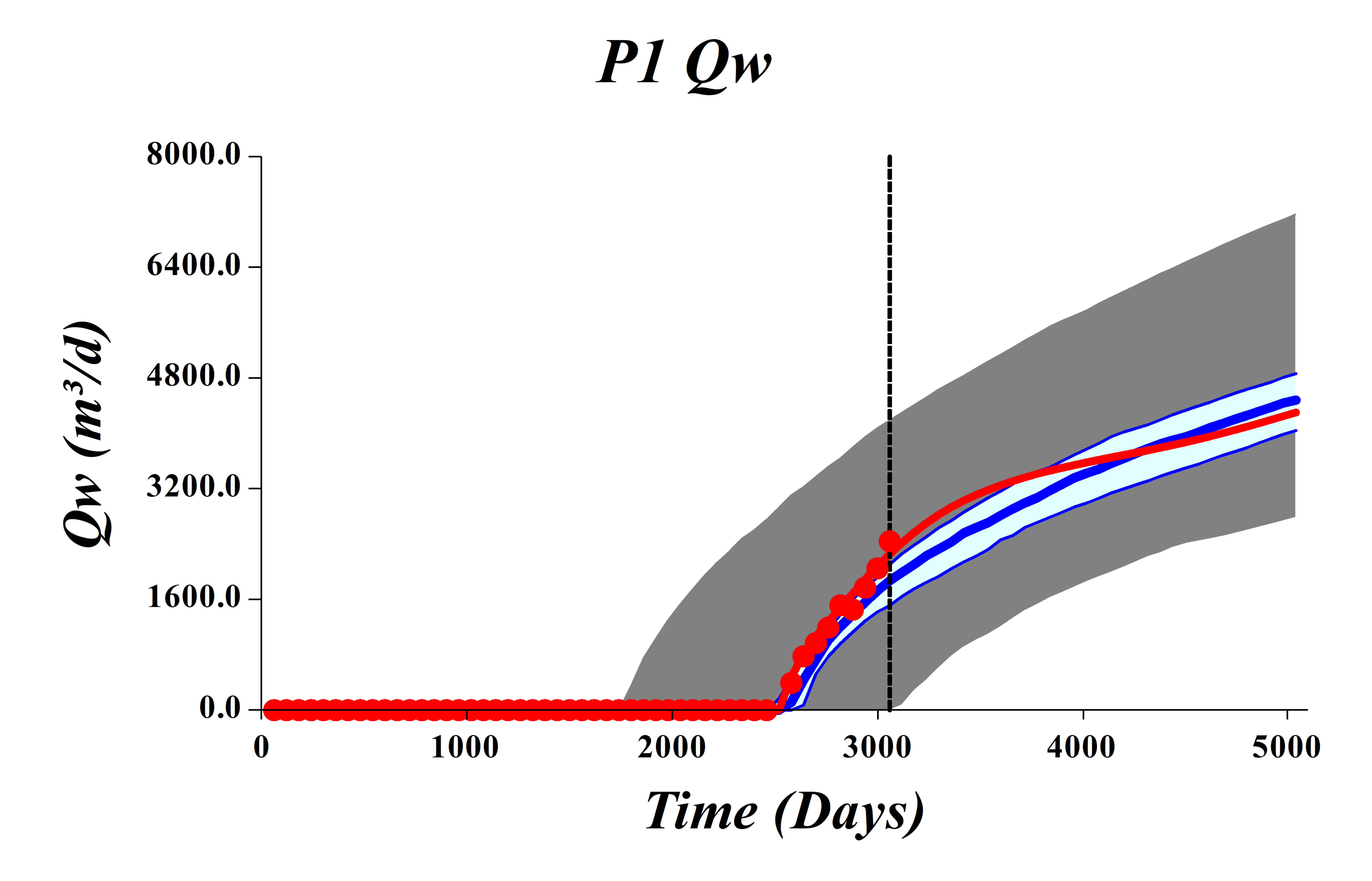}
    }
    \subfloat[]{
      \includegraphics[width=0.5\textwidth]{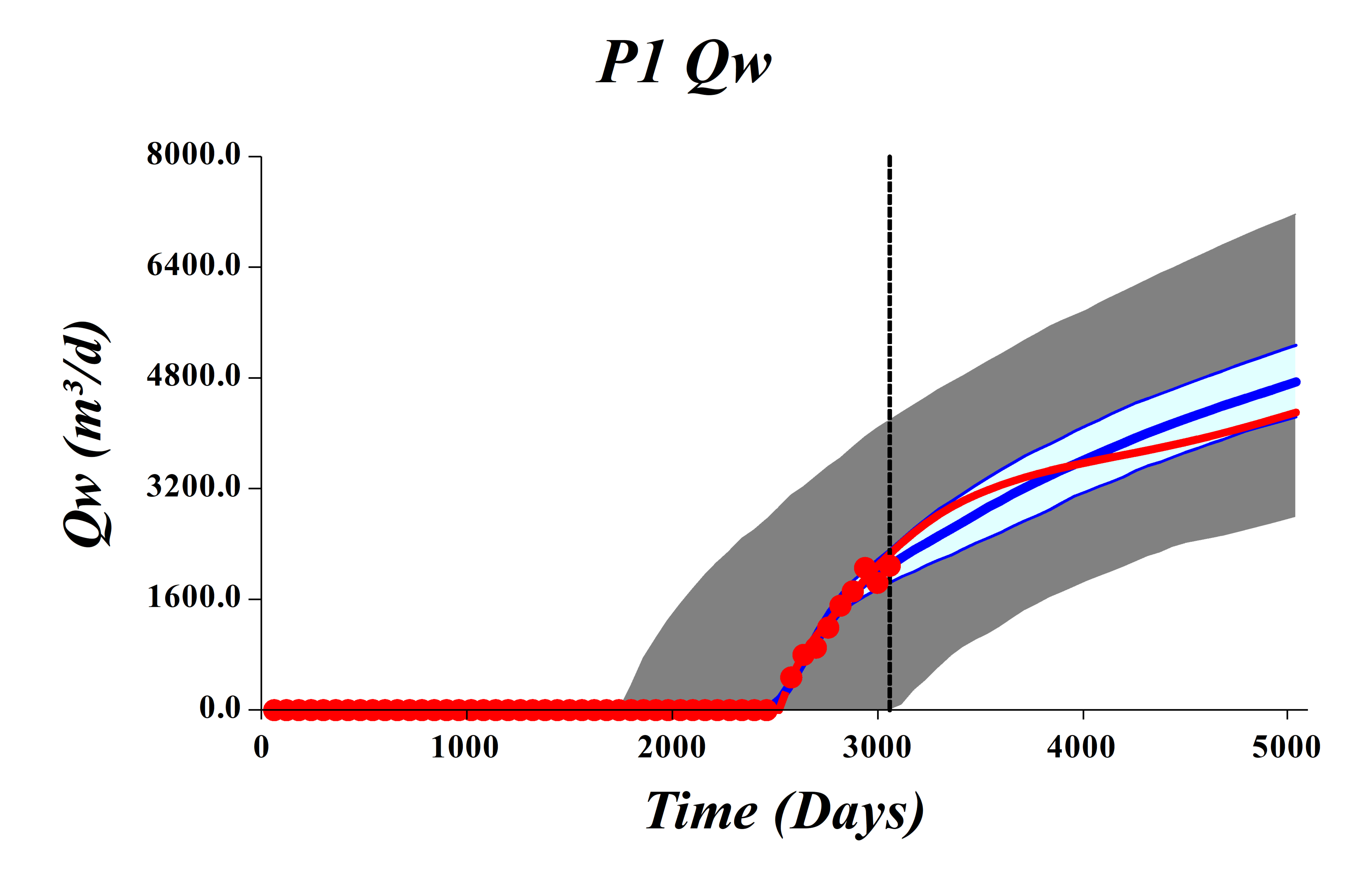}
    }
    \linebreak
    \subfloat[]{
      \includegraphics[width=0.5\textwidth]{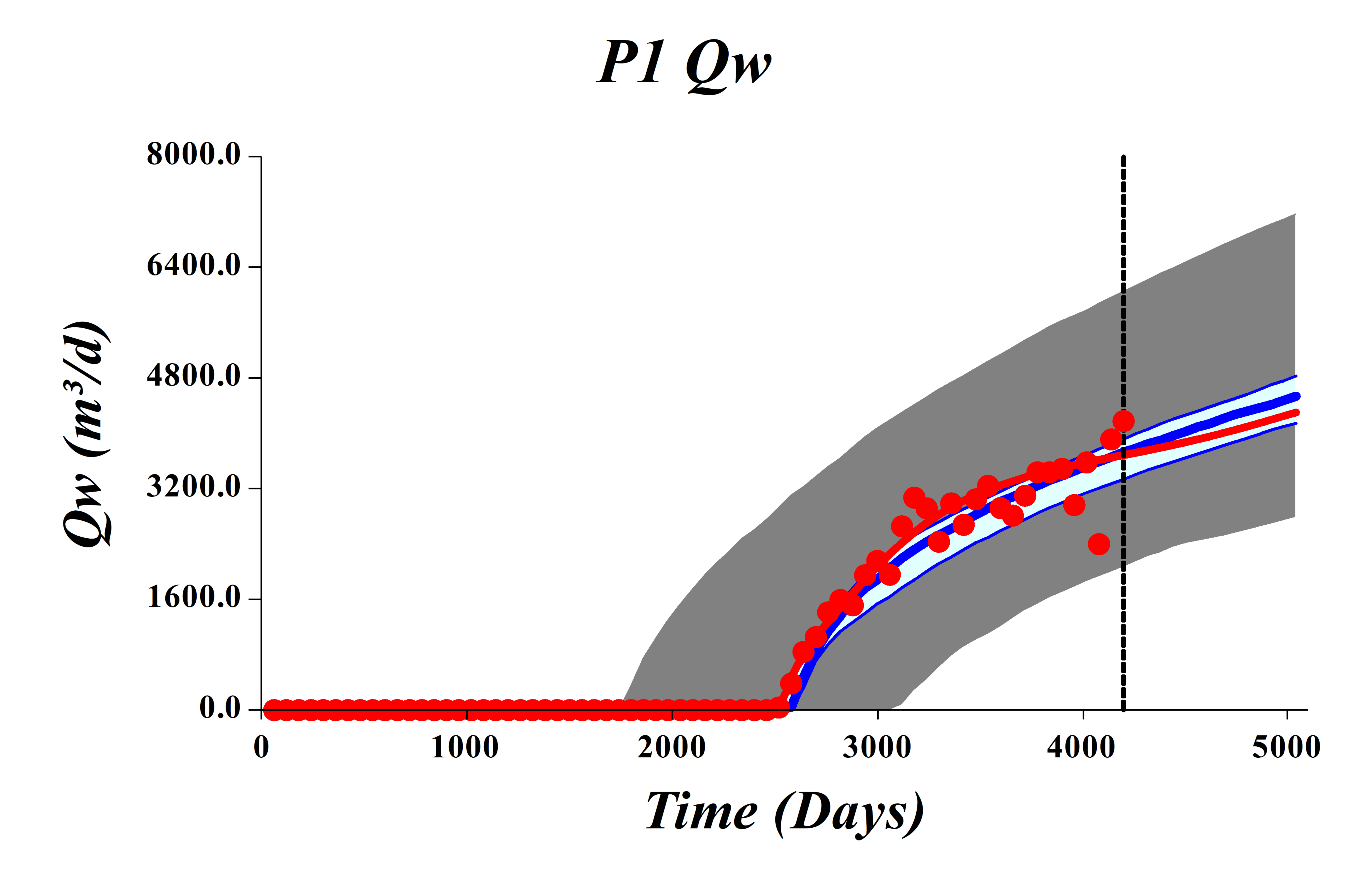}
    }
    \subfloat[]{
      \includegraphics[width=0.5\textwidth]{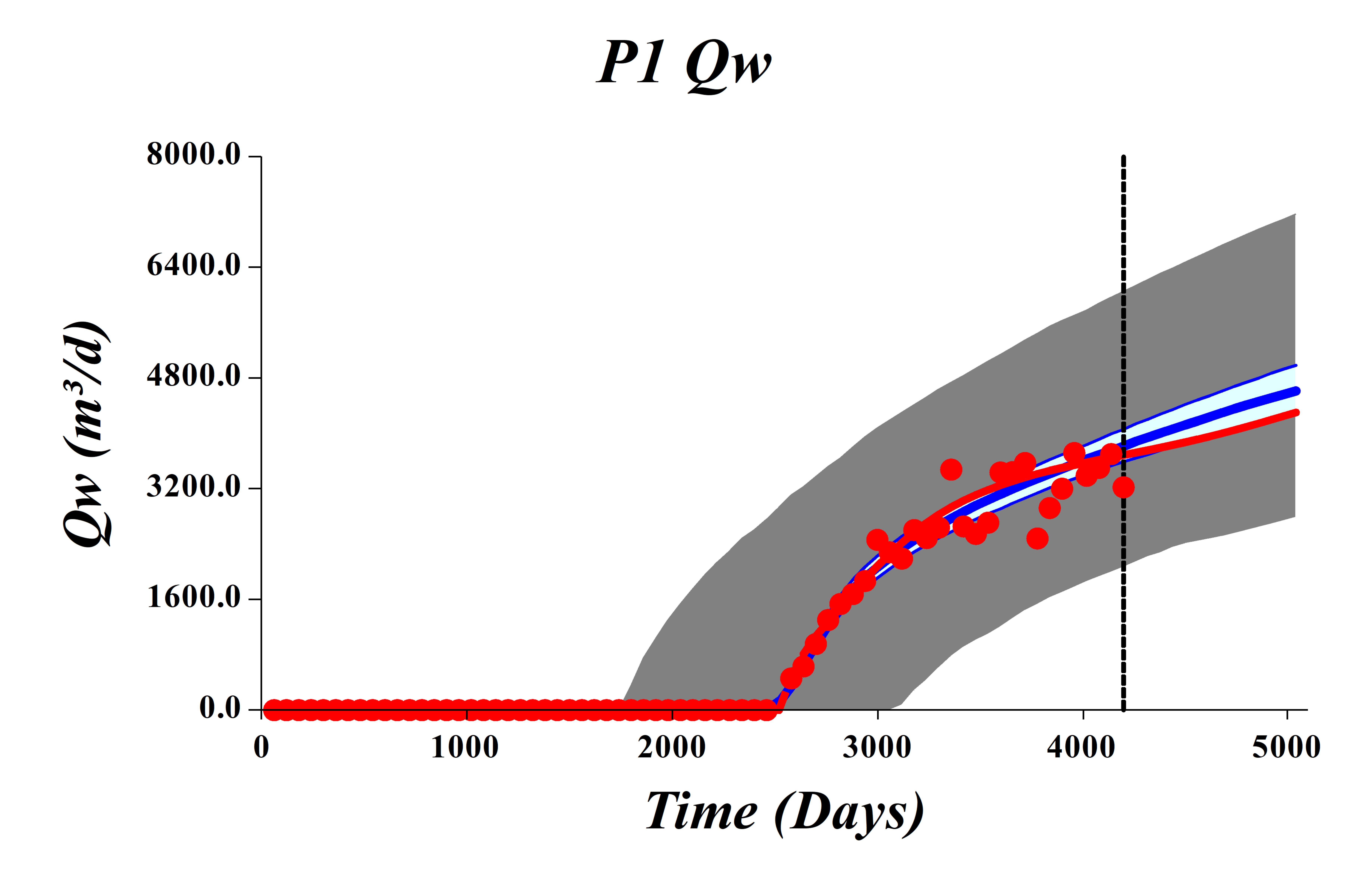}
    }
\captionsetup{justification=justified}
\caption{Water production rate in m$^\text{3}$/days for well P1 considering different number of observed data points 20 (first row), 51 (second row) and 70 data points (third row). Test case 1. (a), (c) and (e) DSI, (b), (d) and (f) DSI-ESMDA. See Fig.~\ref{Fig:WaterRateTestCase1} for description.}
\label{Fig:WaterRateTestHistCase1}
\end{figure}

\begin{figure}
\centering
    \captionsetup{justification=centering}
    \subfloat[]{
      \includegraphics[width=0.5\textwidth]{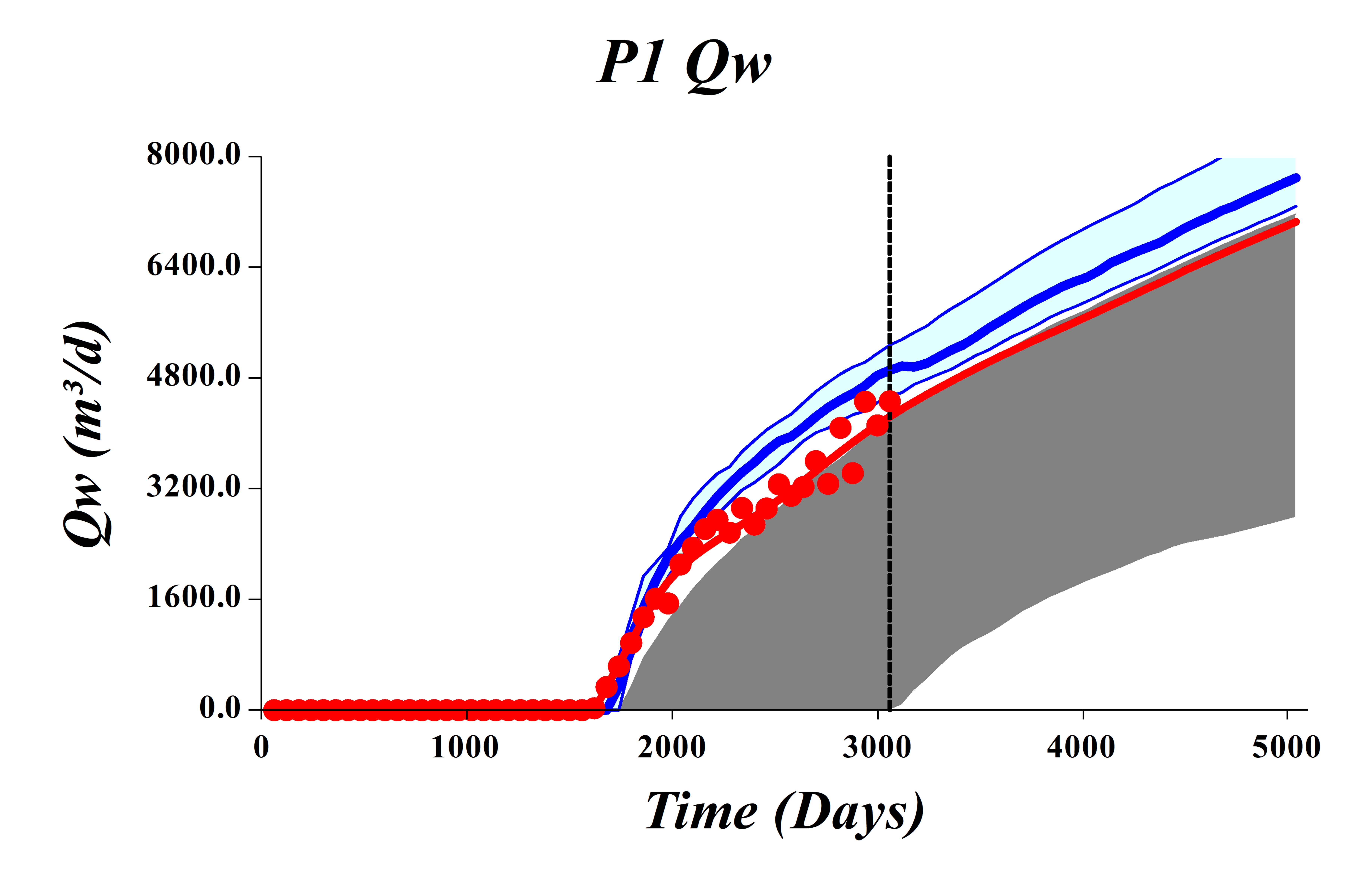}
    }
    \subfloat[]{
      \includegraphics[width=0.5\textwidth]{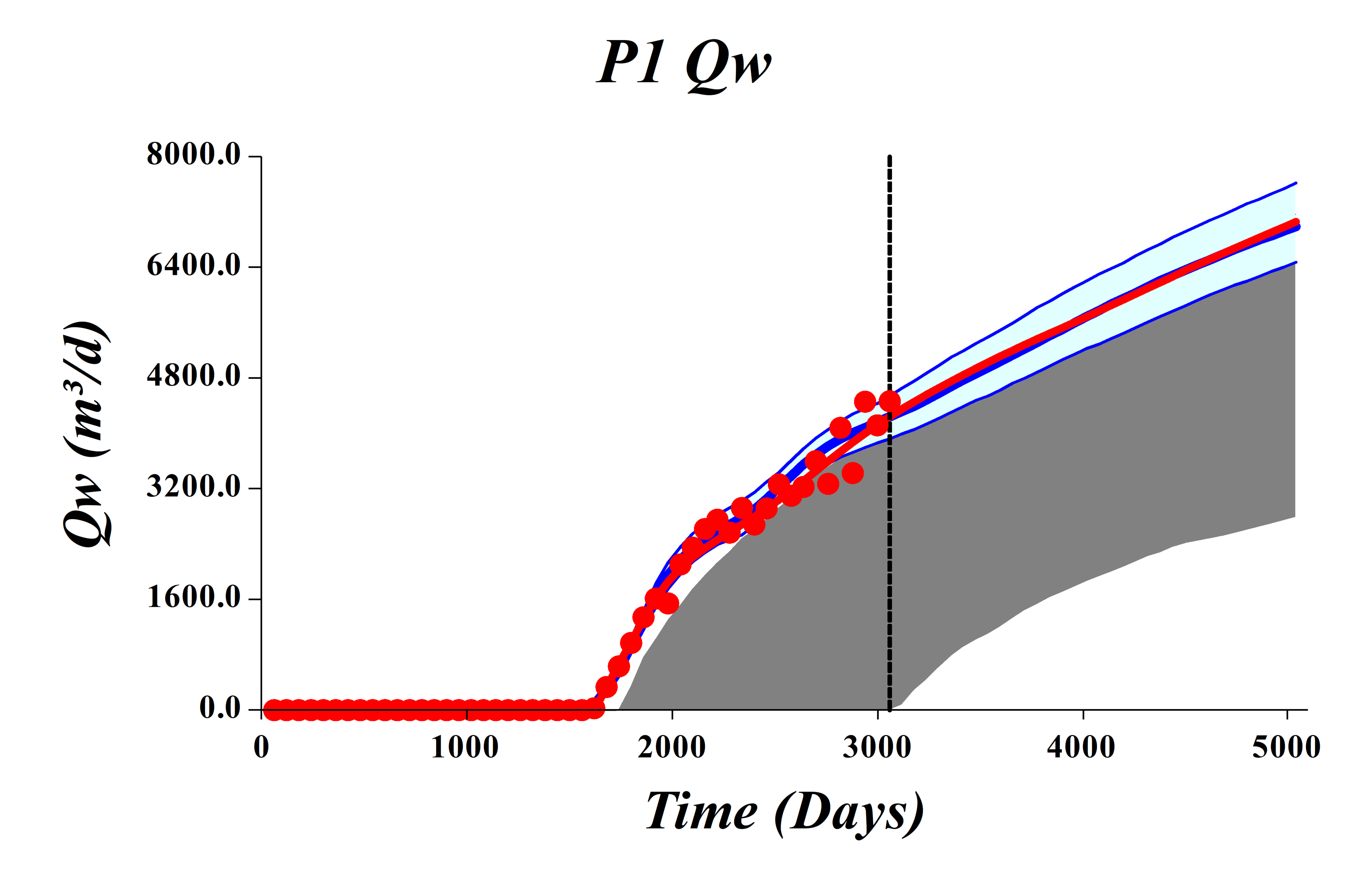}
    }
    \linebreak
    \subfloat[]{
      \includegraphics[width=0.5\textwidth]{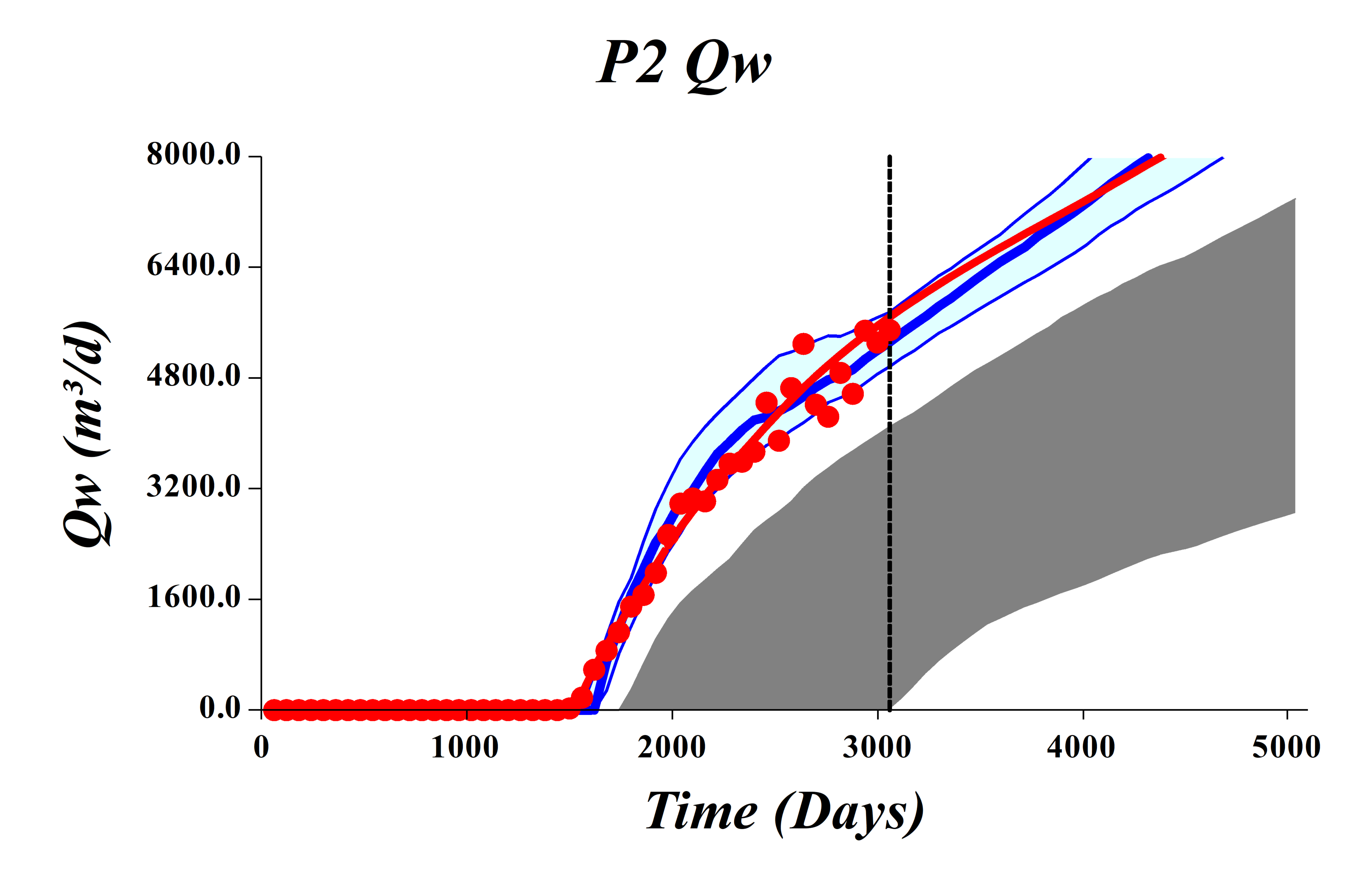}
    }
    \subfloat[]{
      \includegraphics[width=0.5\textwidth]{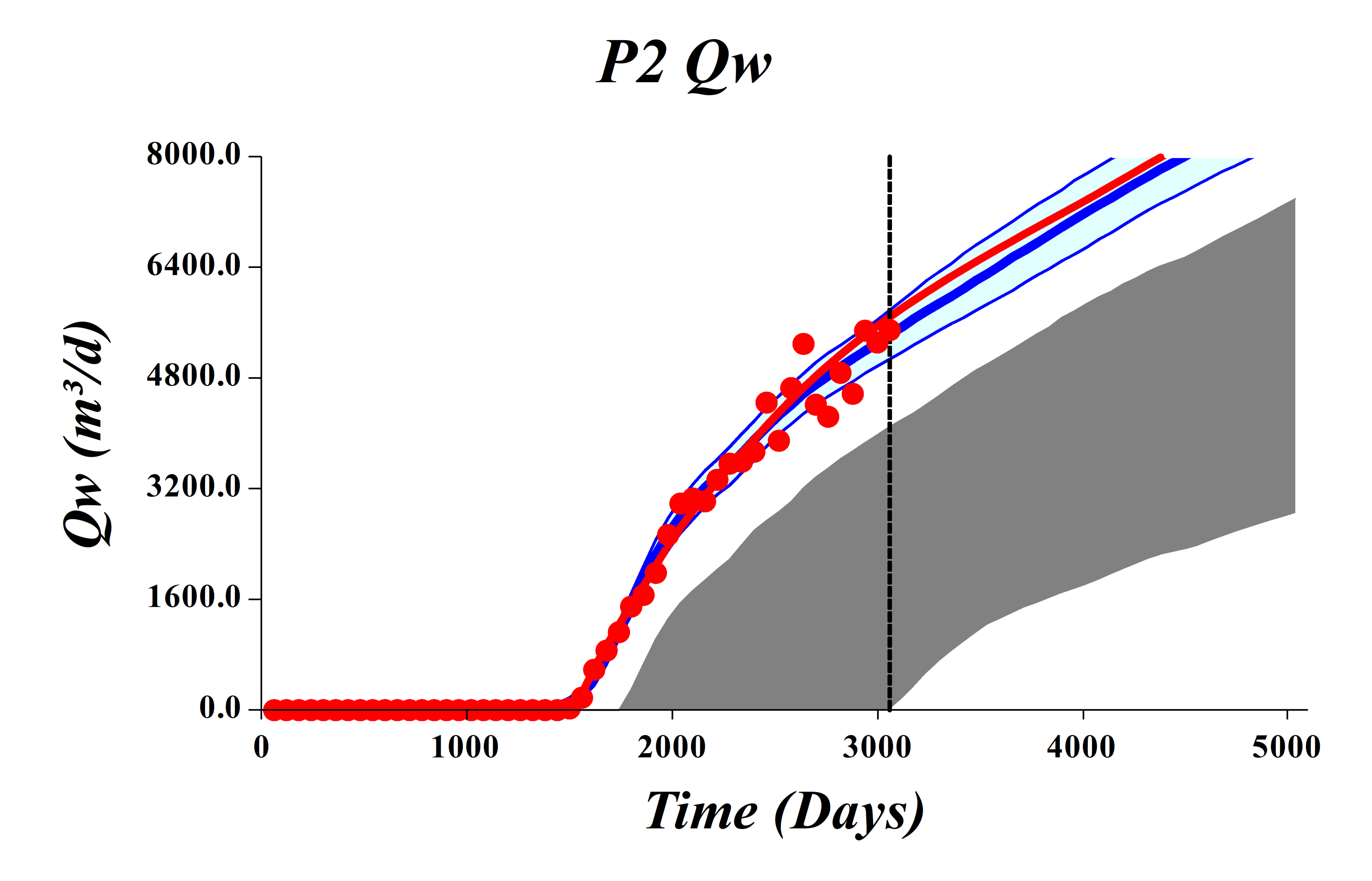}
    }

\captionsetup{justification=justified}
\caption{Water production rate in m$^\text{3}$/days for wells P1 (first row) and P2 (second row) for a case with poor coverage of predictions from the prior ensemble. Test case 1. (a) and (c) DSI, (b) and (d) DSI-ESMDA. See Fig.~\ref{Fig:WaterRateTestCase1} for description.}
\label{Fig:WaterRateBiasedReferenceTestCase1}
\end{figure}

Figure~\ref{Fig:WaterRateTestCase1} shows the predicted water production rate for both wells of the model obtained with DSI, DSI-ESMDA. For comparisons, we included in each plot the water rate predicted by the history-matched models with ES-MDA. The results in this figure indicate that both DSI methods obtained similar predictions, which are in reasonable agreement with the standard ES-MDA. Figure~\ref{Fig:NpWpCase1} shows the field cumulative production of oil and water predicted by the prior ensemble and by each method. The cumulative production from the reference case is also included in this figure for comparisons. This figure shows that the uncertainty ranges predicted by the methods are similar. Figure~\ref{Fig:WaterRateTestHistCase1} shows the predicted water production rate for the first well considering three different sizes of the historical period. The main difference is observed for the case with no water breakthrough in the production history (cases with 20 data points in Fig.~\ref{Fig:WaterRateTestCase1}). For this case, DSI-ESMDA predicts a larger uncertainty range in the water production rate. In order to test the robustness of the methods, we selected another reference model with predictions outside the P10--P90 range from the prior ensemble and the results are presented in Fig.~\ref{Fig:WaterRateBiasedReferenceTestCase1}. Despite of the poor coverage of the prior ensemble, both methods obtained forecasts in reasonable agreement with the reference. However, we note that DSI-ESMDA obtained a sightly better data match for the first well (Figs.~\ref{Fig:WaterRateBiasedReferenceTestCase1}a and b) and a forecast range centered around the reference.

Both methods require to run reservoir simulations of the prior ensemble, which is clearly the dominating computational time. DSI also requires to solve the RML optimizations with the BFGS method, while DSI-ESMDA requires to apply the analysis (Eq.~\ref{Eq:DSI-ESMDA}) $N_a = 4$ times for each ensemble member. We measured the CPU running time of the inversion part of each method and the results are reported in Table~\ref{Tab:CPUTestCase1}. All computations were performed in the same computer (Intel Core i5-4690 CPU 3.5~GHz and 16~GB RAM). The results in this table show that DSI-ESMDA is notably faster than DSI, requiring only two seconds against approximately 25 minutes for DSI.

\begin{table}
	\centering
	\caption{CPU time for inversion (70 data points). Test case 1}
	\label{Tab:CPUTestCase1}
	\begin{tabular}{lc}
		\toprule
		Method & CPU time (seconds) \\
		\midrule
		DSI & 1506 \\
        DSI-ESMDA & 2 \\
		\bottomrule
	\end{tabular}
\end{table}

\clearpage

\subsection{Test Case 2}
\label{Sec:TestCase2}

\begin{figure}
	\centering
	\includegraphics[width=0.7\linewidth]{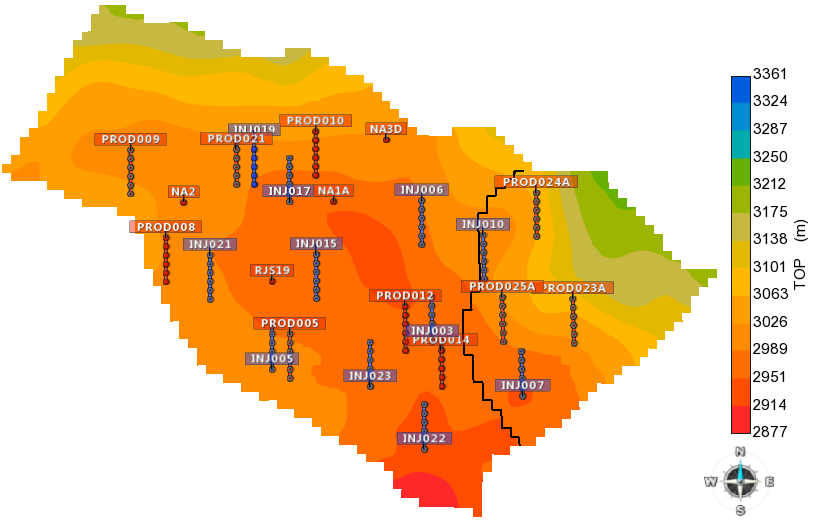}
	\caption{Position of the wells. UNISIM-I-H}
	\label{Fig:UNISIM-Wells}
\end{figure}

The second test case is a more realistic history-matching problem known as UNISIN-I-H \citep{avansi:15a}. This problem is based on actual data from Namorado Field (Campos Basis, Brazil). The UNISIM-I-H model has 81 $\times$ 58 $\times$ 20 gridblocks, but only 37,000 are active. All gridblocks have a uniform size of 100~m $\times$ 100~m $\times$ 8~m. The original dataset is available for download at \citep{unisim-i-h:13}. The dataset consists of 500 realizations of petrophysical properties (porosity, permeability in the three orthogonal directions and net-to-gross ratio). Besides the 500 petrophysical realizations, the UNISIM-I-H case also includes six global parameters, whose prior uncertainties were modeled as independent triangle distributions with values shown in Table~\ref{Tab:UNISIM-ScalarPar}. In this field there are 25 long horizontal wells (14 producers and 11 water injectors). Figure~\ref{Fig:UNISIM-Wells} shows the position of the wells projected in the first layer of the model. The oil producing wells are perforated close to the top and the water injection wells are perforated near to the bottom of the reservoir. The observed data were generated by adding random Gaussian data-error to the data predicted by a reference fine-scale model (UNISIM-I-R). The UNISIM-I-R was constructed with a higher level of geological details in a grid with 3.5 million active gridblocks \citep{avansi:15a}. The observed data correspond to monthly ``measurements'' of oil and water rate and the noise level was assumed as 10\% of the data predicted by the model UNISIM-I-R. All wells are controlled by specified BHP during the historical and forecast periods.

\begin{table}
\caption{Prior distribution of global parameters. UNISIM-I-H}
\label{Tab:UNISIM-ScalarPar}
\begin{scriptsize}
\begin{center}
\begin{tabular}{lccc}
\toprule
Parameter & Mode & Min. & Max. \\
\midrule
Rock compressibility in $(\text{cm}^2/\text{kgf})$ & $5.3\times 10^{-5}$ & $1\times 10^{-5}$ & $9.6\times 10^{-5}$ \\
Depth of the oil-water contact at the East block (Fig.~\ref{Fig:UNISIM-Wells}) & 3174 & 3169 & 3179 \\
Maximum water relative permeability & 0.35 & 0.15 & 0.55 \\
Corey exponent of water relative permeability & 2.4 & 2.0 & 3.0 \\
Multiplier for vertical permeability & 1.5 & 0.0 & 3.0 \\
\bottomrule
\end{tabular}
\end{center}
\end{scriptsize}
\end{table}

For this test problem, we applied DSI, DSI-ESMDA and DSI-ESMDA with localization. We also applied a model-based inversion (history matching) with ES-MDA for comparisons. The dimension of the matrix $\DD$ for the UNISIM-I-H case is $12825 \times 500$ and the SVD retaining 99\% of the singular values energy resulted in a vector of PCA coefficients with 395 elements for DSI. For DSI-ESMDA, we used $N_a = 4$ data assimilations and kept 99\% of the singular values in the subspace inversion. DSI-ESMDA with localization considered a critical length computed with $L_{x^\prime} = L_{y^\prime} = 2000$~m and $T = 6000$~days. These values were selected based on our previous experience testing ES-MDA for history matching this benchmark problem. We did not try different configurations for the critical length. The history matching with ES-MDA also used $N_a = 4$ and a spatial localization with critical length of 2000~m. Note that in the history matching case, the Kalman gain is applied to update model parameters (porosity, permeability, etc.). Therefore, it does not seem appropriate to introduce the ``time distance'' to compute the localization matrix in this case.

In ours tests, the DSI implementation failed to achieve reasonable data matches for the UNISIM-I-H case. This result came as a surprise because the same problem was not observed with DSI-ESMDA using the same prior ensemble. We conducted a series of experiments trying to solve this problem, but the reason for the poor performance of the optimizations is still unclear. Among our experiments, we tried to use both numerical and the approximate analytical gradients during the minimizations, however, the large majority of the optimization did not converge properly. We also introduced a rescaling procedure before the truncated SVD. Basically, instead of computing the SVD of the matrix $\DD$ (as in Eq.~\ref{Eq:SVD}), we applied SVD to a scaled version

\begin{equation}\label{Eq:SVD-rescaled}
  \widehat{\DD} \equiv \Ce\invsqr \DD = \widehat{\mathbf{U}} \widehat{\mathbf{\Sigma}} \widehat{\mathbf{V}}\trp,
\end{equation}
and computed the square root of $\C_{\dsim}$ for PCA as

\begin{equation}
  \C_{\dsim}\sqr = \Ce\sqr \widehat{\mathbf{U}} \widehat{\mathbf{\Sigma}}.
\end{equation}
The rationality for this procedure is to avoid losing relevant information during the truncation of small singular values because $\DD$ may be poorly scaled. Unfortunately, this procedure did not result in significant improvements. We also tried the optimizations keeping all 499 singular values to preserve the maximum number of degrees freedom provided by the prior ensemble to match data. In this case, we observed an even worse performance. Conversely, we also tried to reduce the number of singular values by keeping only 95\% and 90\% of the energy, which resulted in 237 and 154 singular values, respectively. The idea was to check if using a small number of PCA coefficients would help to regularize the optimizations. Again, no noticeable improvements were observed.

In order to evaluate the quality of the final data matches, we computed the data-mismatch objective function normalized by the number of observed data points ($N_{d,\textrm{h}} = 9900$) for 500 posterior realizations obtained with DSI, DSI-ESMDA, DSI-ESMDA with localization and the history matching with ES-MDA. The normalized data mismatch objective function was computed as

\begin{eqnarray}\label{Eq:ON}
  \nonumber \mathcal{O}_{N,d} & = & \frac{1}{2N_{d,\textrm{h}}}\left( \dobs - \dsim_{\textrm{h}}\right)\trp \Ce\inv \left( \dobs - \dsim_{\textrm{h}}\right)\\
  & = & \frac{1}{2N_{d,\textrm{h}}} \sum_{i=1}^{N_d} \left( \frac{d_{\textrm{obs},i} - d_{\textrm{h},i}}{\sigma_{e,i}}\right)^2,
\end{eqnarray}
where $\sigma_{e,i}$ is the data-error standard deviation of the $i$th datum. The last equality in (\ref{Eq:ON}) holds only for independent data errors. In the Appendix section of \citep{oliver:18b} it is shown that the expectation of $\mathcal{O}_{N,d}$ for a set of RML samples in a linear problem should be one half. Even though this value is valid only for linear problems, it serves as a reference. For example, if we have a predicted curve where the difference with all observed data points is exactly one standard deviation of the data error, then the normalized objective function of this curve is exactly 0.5. Analogously, we have $\mathcal{O}_{N,d} = 2$ and $\mathcal{O}_{N,d} = 4.5$ for two and three standard deviations, respectively. This effectively means that it would be difficult to justify in practice the acceptance of a data match with $\mathcal{O}_{N,d}$ significantly larger than five. Table~\ref{Tab:UNISIM-OF} shows the mean and standard deviation of $\mathcal{O}_{N,d}$ for each case. Clearly, DSI resulted too large values for $\mathcal{O}_{N,d}$, indicating the lack of convergence of the minimizations. The objective functions from DSI-ESMDA and DSI-ESMDA with localization have mean roughly ten times smaller. The history matching case with ES-MDA also resulted in a mean objective function significantly smaller than DSI. In terms of computational cost, the DSI method required approximately 28 hours to complete, while the DSI-ESMDA (with and without localization) require around 45 minutes each.

\begin{table}
\caption{Normalized data-mismatch objective function. UNISIM-I-H}
\label{Tab:UNISIM-OF}
\begin{center}
\begin{tabular}{lcc}
\toprule
Case & Mean & Standard deviation\\
\midrule
Prior & 456.7 & 333.462 \\
DSI & 26.6 & 7.379 \\
DSI-ESMDA & 2.7 & 0.003 \\
DSI-ESMDA with localization & 1.9 & 0.008 \\
ES-MDA & 3.0 & 0.109 \\
\bottomrule
\end{tabular}
\end{center}
\end{table}

\begin{figure}
\centering
    \captionsetup{justification=centering}
    \subfloat[]{
      \includegraphics[width=0.33\textwidth]{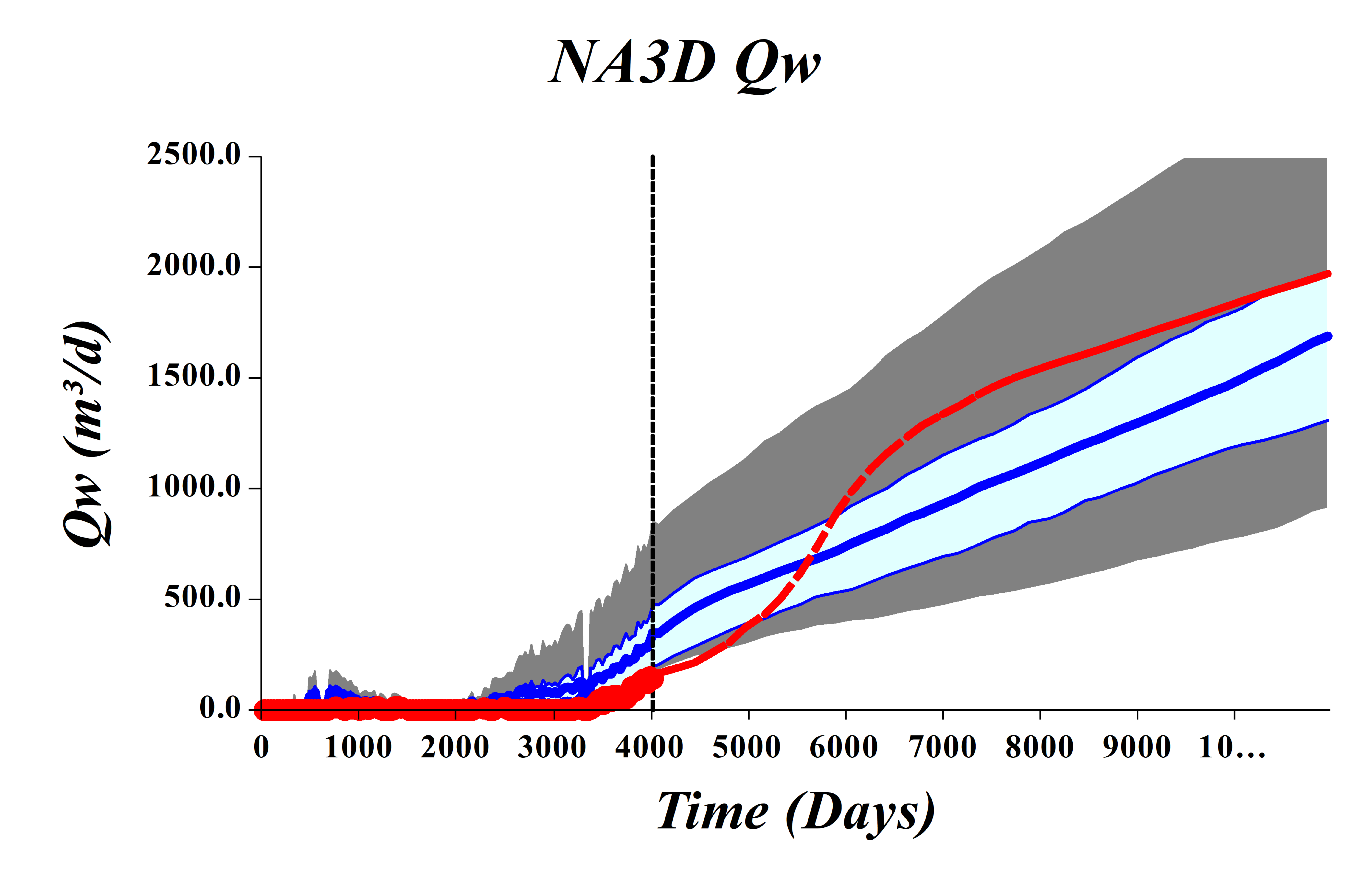}
    }
    \subfloat[]{
      \includegraphics[width=0.33\textwidth]{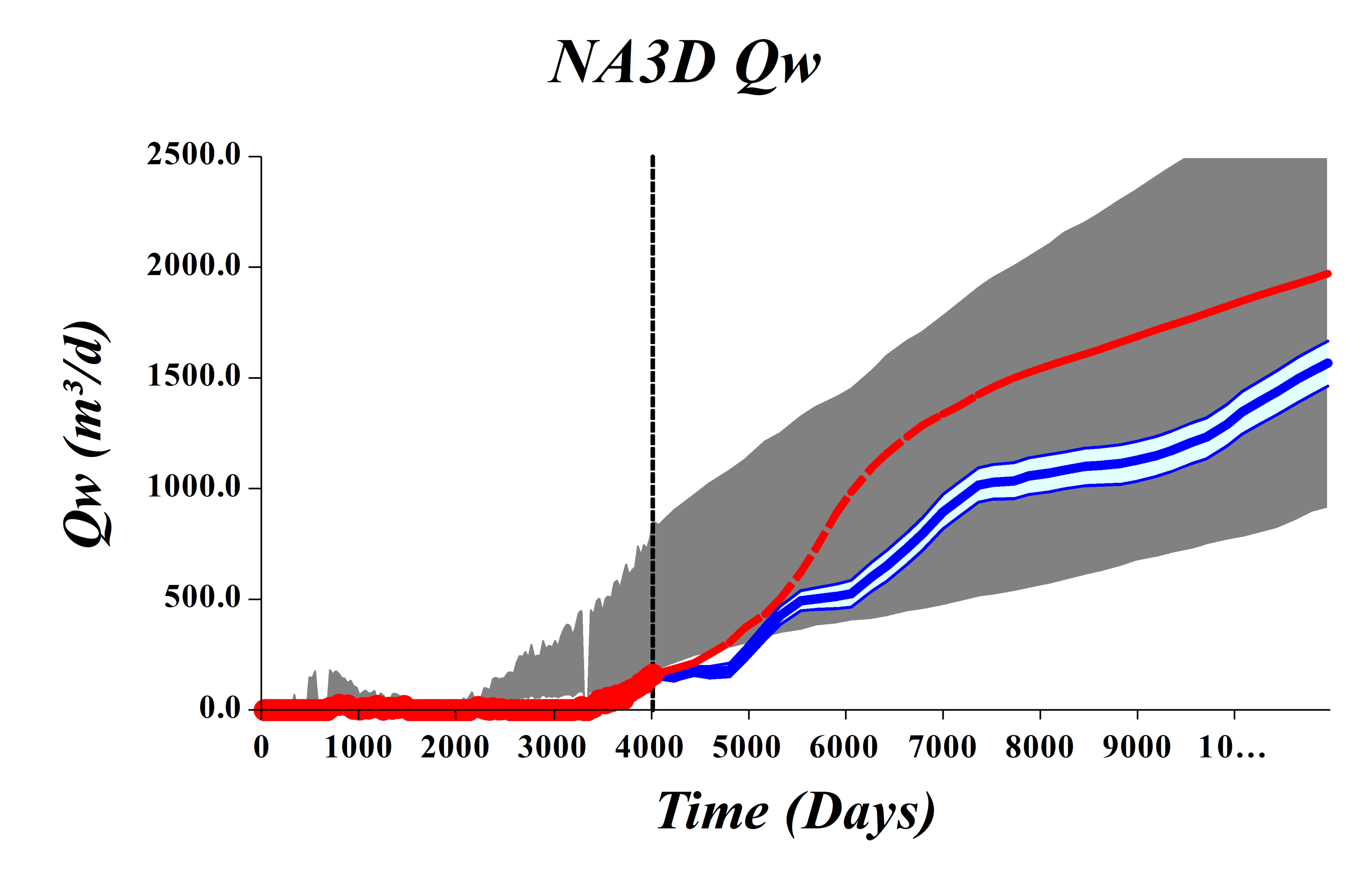}
    }
    \subfloat[]{
      \includegraphics[width=0.33\textwidth]{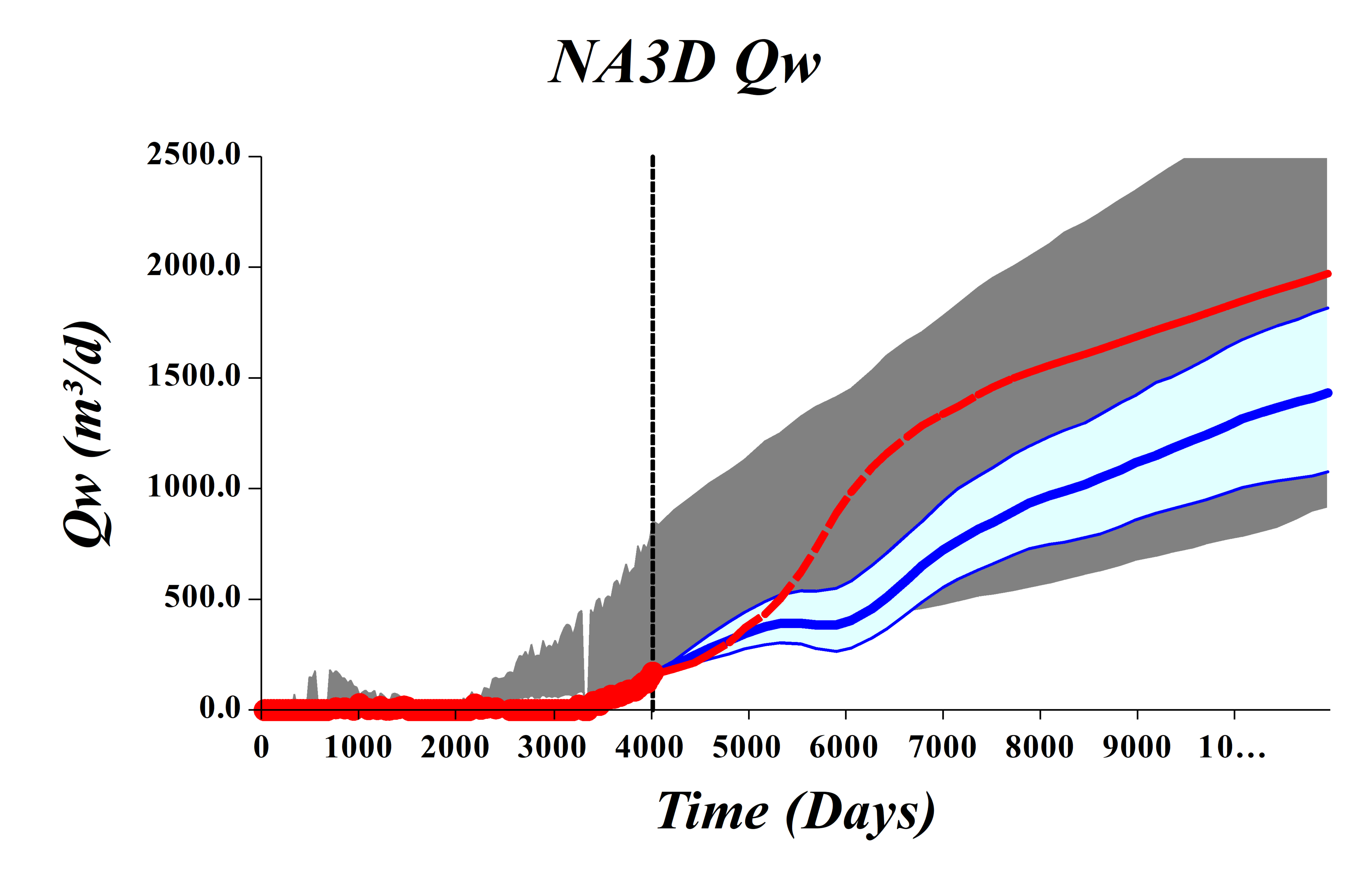}
    }
    \linebreak
        \subfloat[]{
      \includegraphics[width=0.33\textwidth]{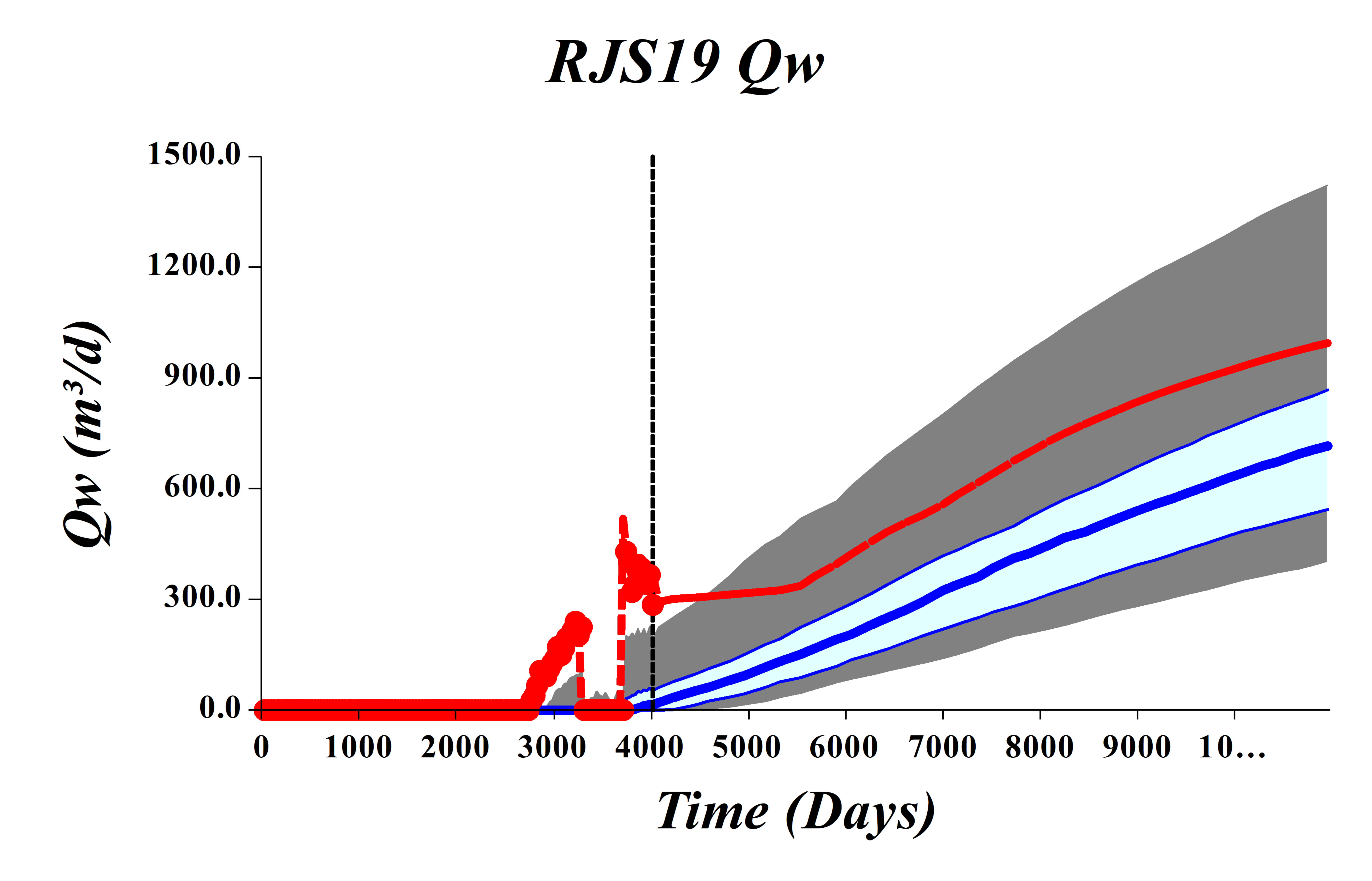}
    }
    \subfloat[]{
      \includegraphics[width=0.33\textwidth]{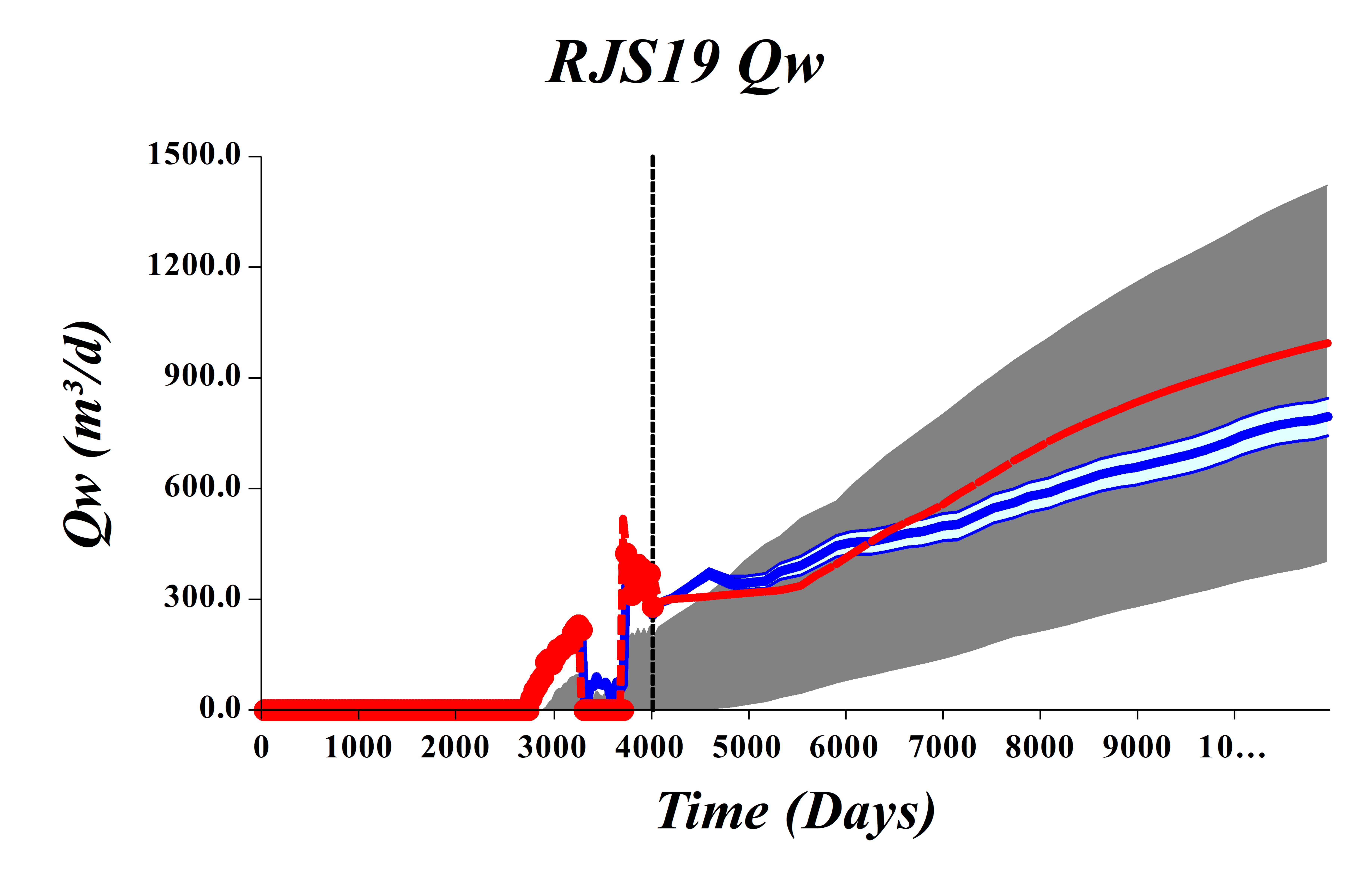}
    }
    \subfloat[]{
      \includegraphics[width=0.33\textwidth]{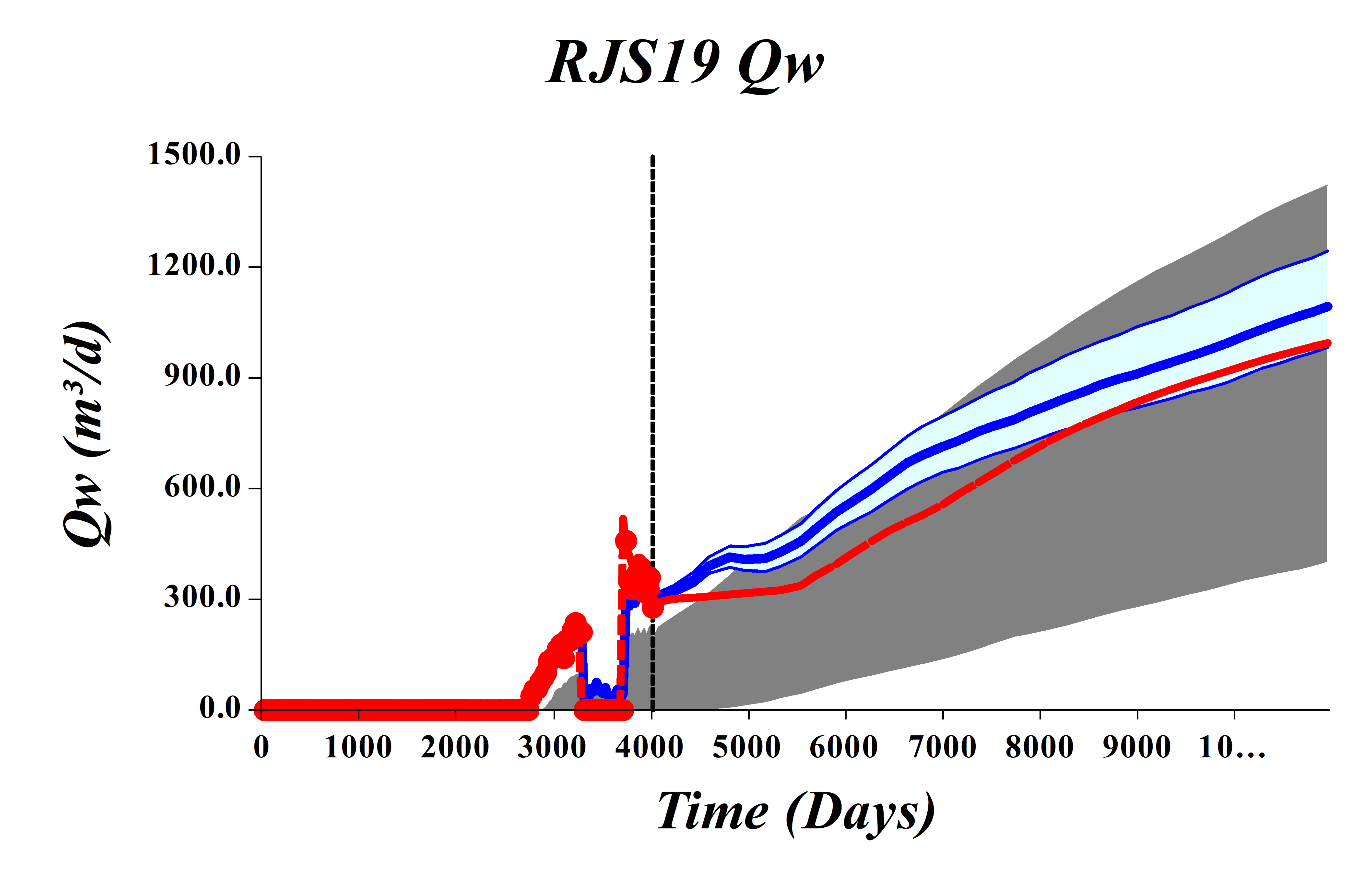}
    }
    \linebreak
        \subfloat[]{
      \includegraphics[width=0.33\textwidth]{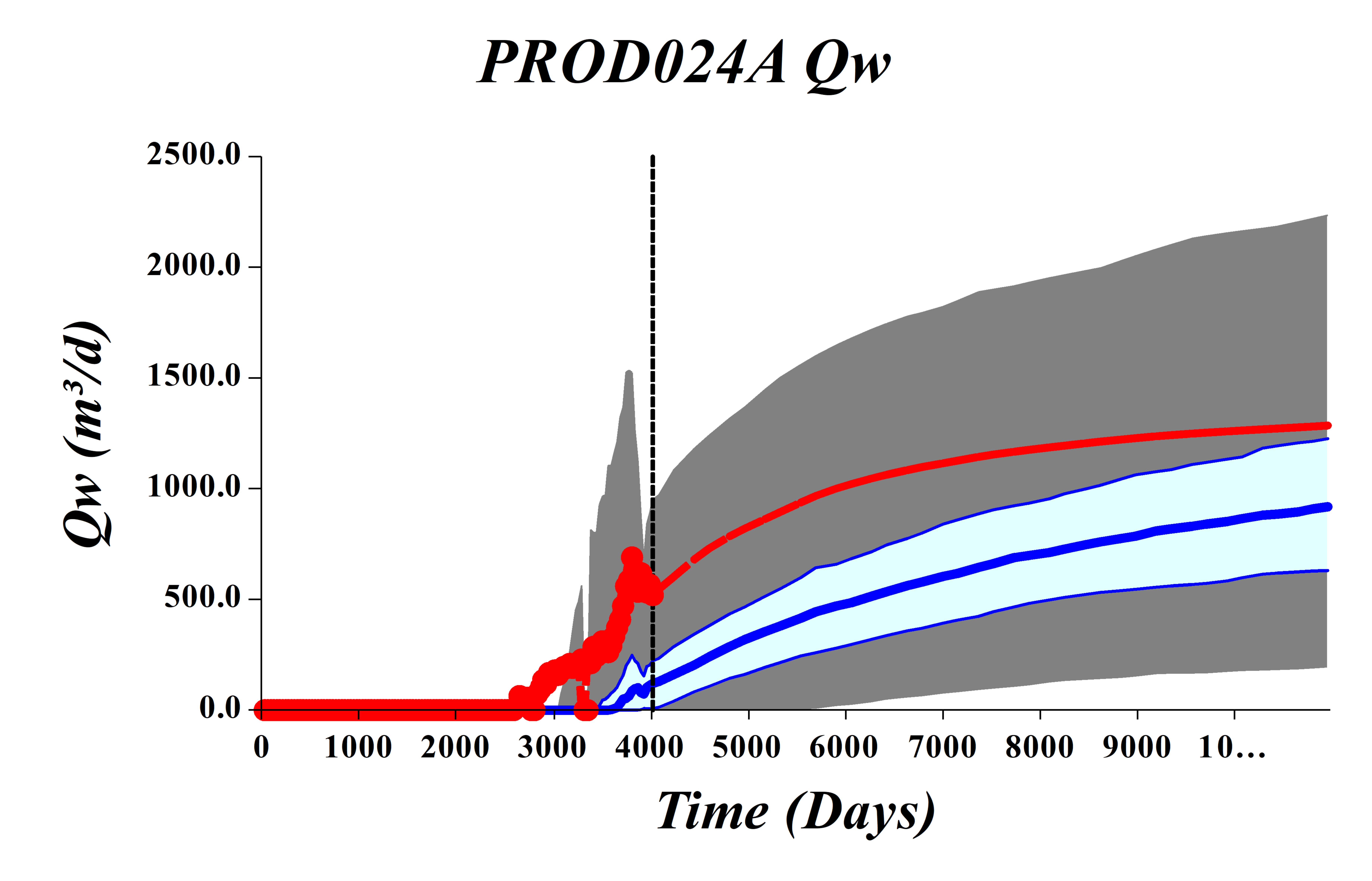}
    }
    \subfloat[]{
      \includegraphics[width=0.33\textwidth]{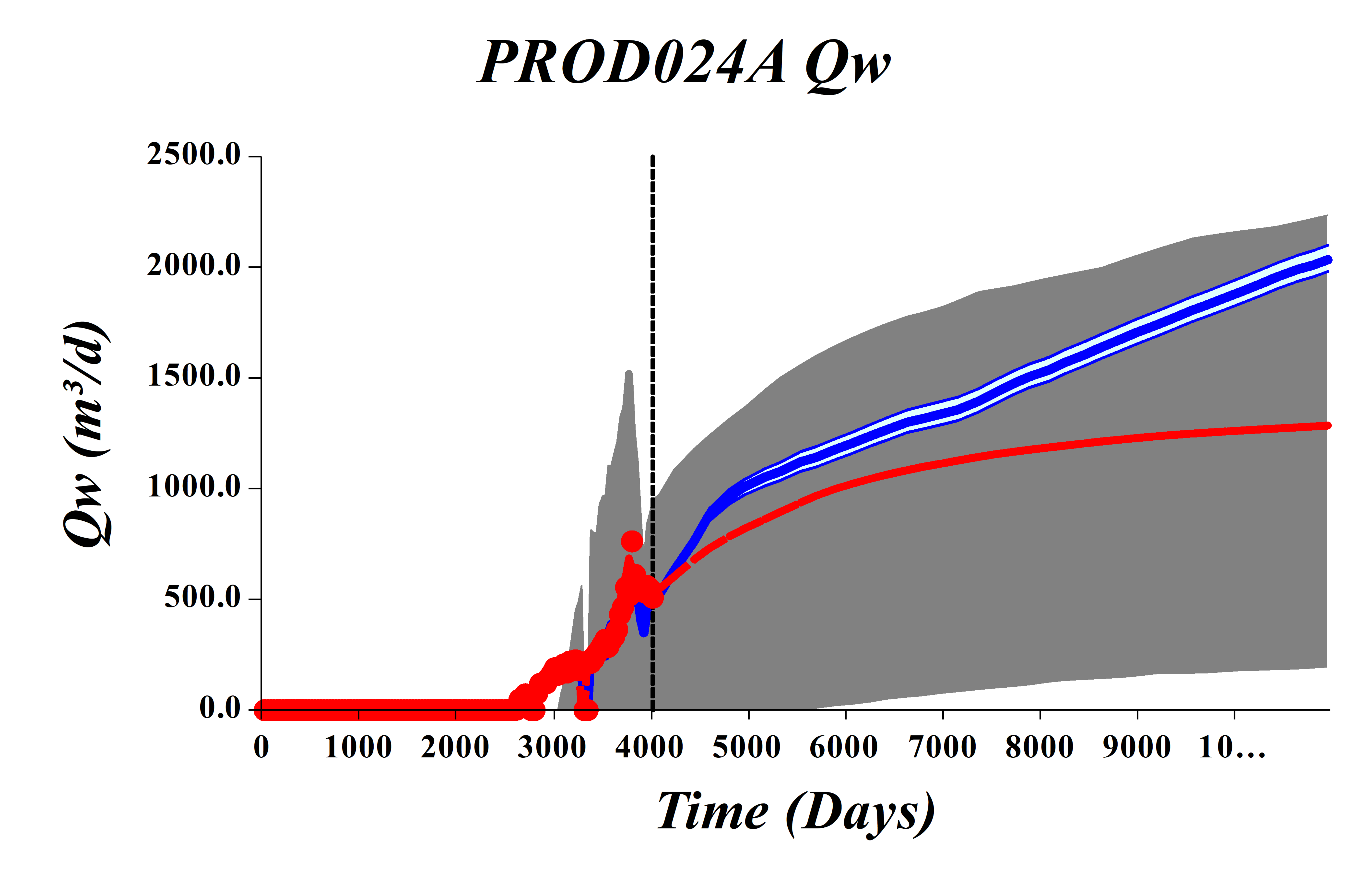}
    }
    \subfloat[]{
      \includegraphics[width=0.33\textwidth]{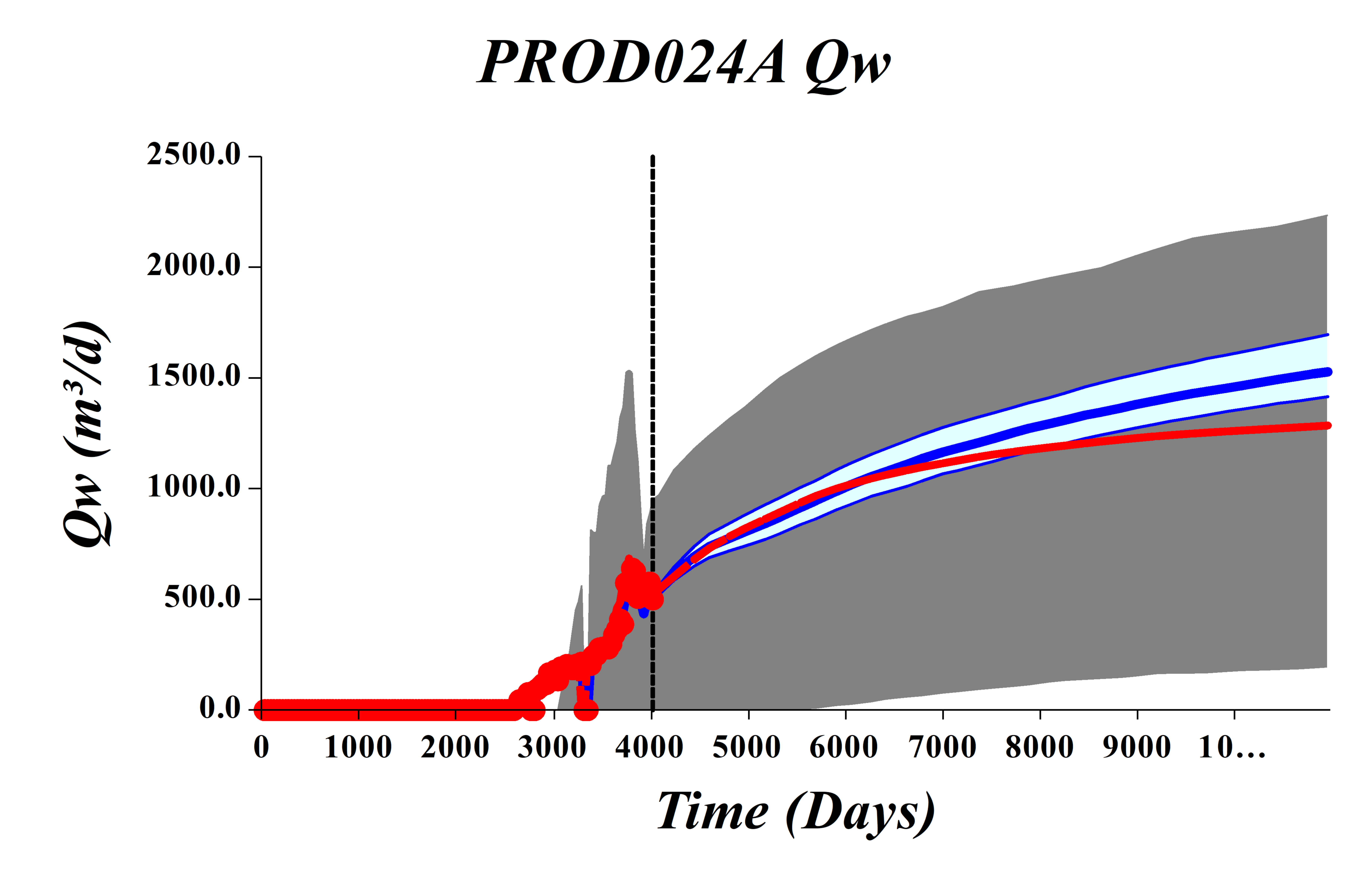}
    }
\captionsetup{justification=justified}
\caption{Water production rate in m$^\text{3}$/days for three wells with good coverage of predictions from the prior ensemble. Well NA3D (first row), RJS19 (second row) and PROD24A (third row). (a), (d) and (g) DSI. (b), (e) and (h) DSI-ESMDA, (c), (f) and (i) DSI-ESMDA with localization. UNISIM-I-H. See Fig.~\ref{Fig:WaterRateTestCase1} for description.}
\label{Fig:WaterRateGoodPriorUNISIM}
\end{figure}

\begin{figure}
\centering
    \captionsetup{justification=centering}
    \subfloat[]{
      \includegraphics[width=0.33\textwidth]{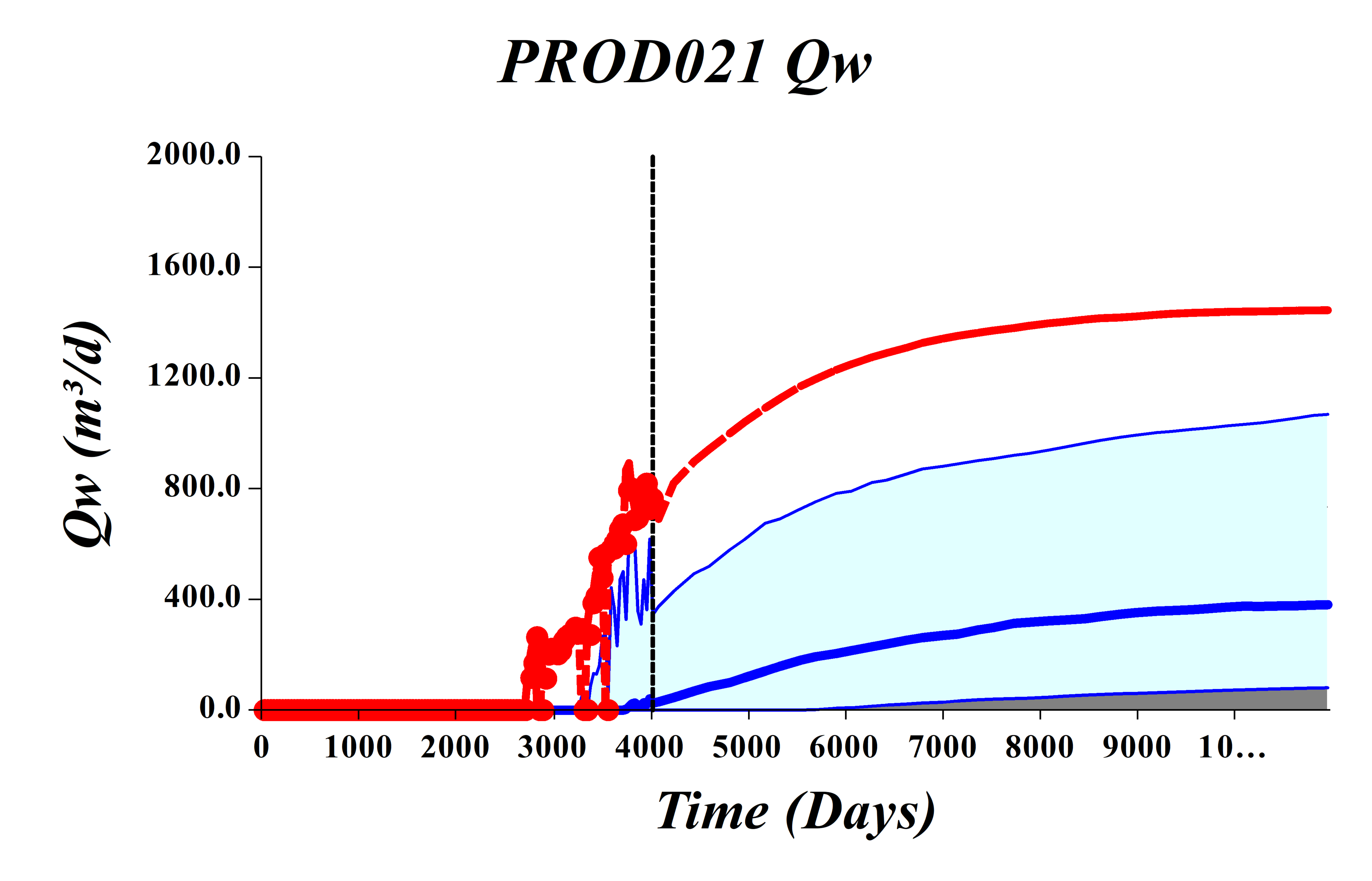}
    }
    \subfloat[]{
      \includegraphics[width=0.33\textwidth]{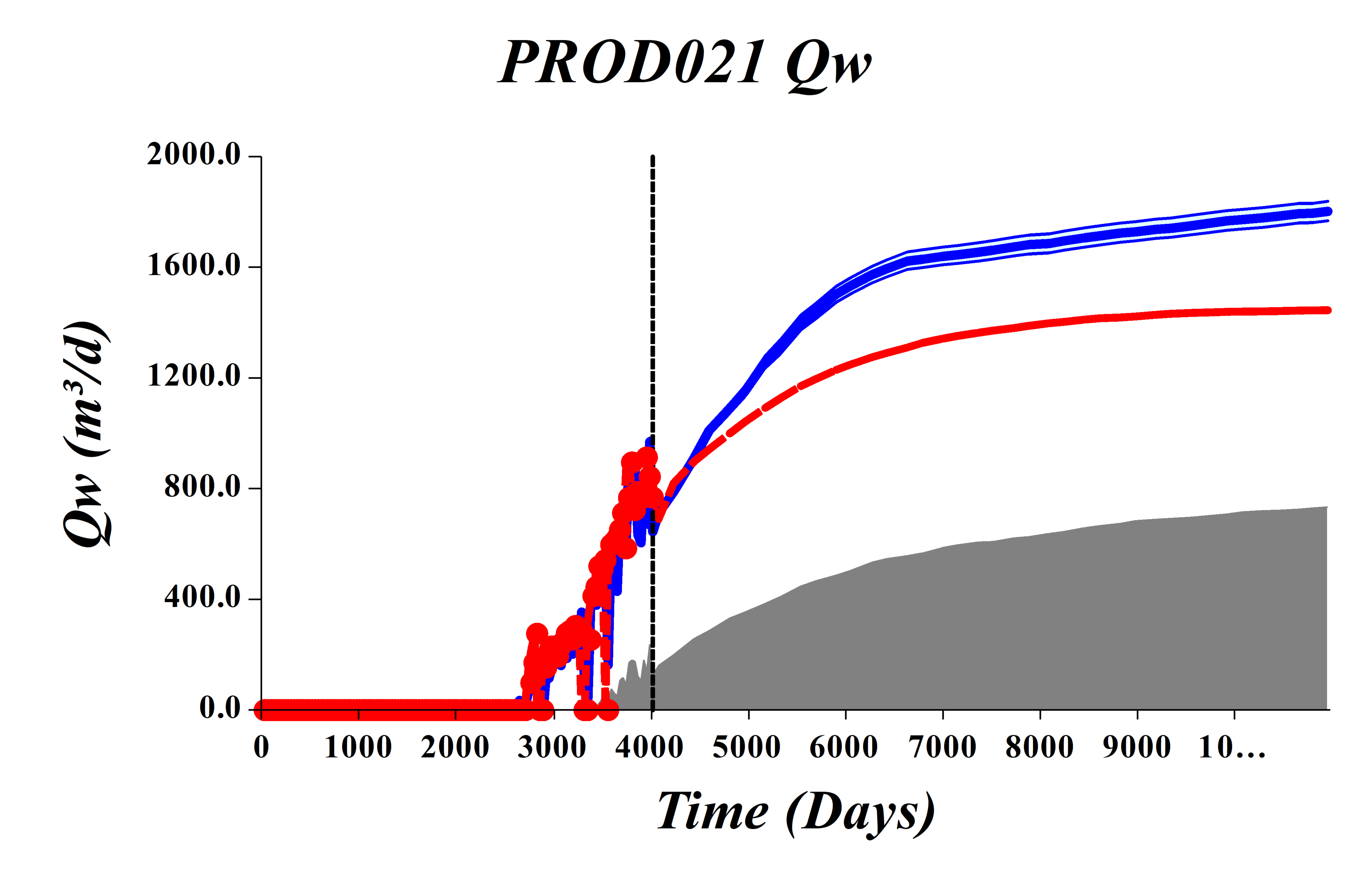}
    }
    \subfloat[]{
      \includegraphics[width=0.33\textwidth]{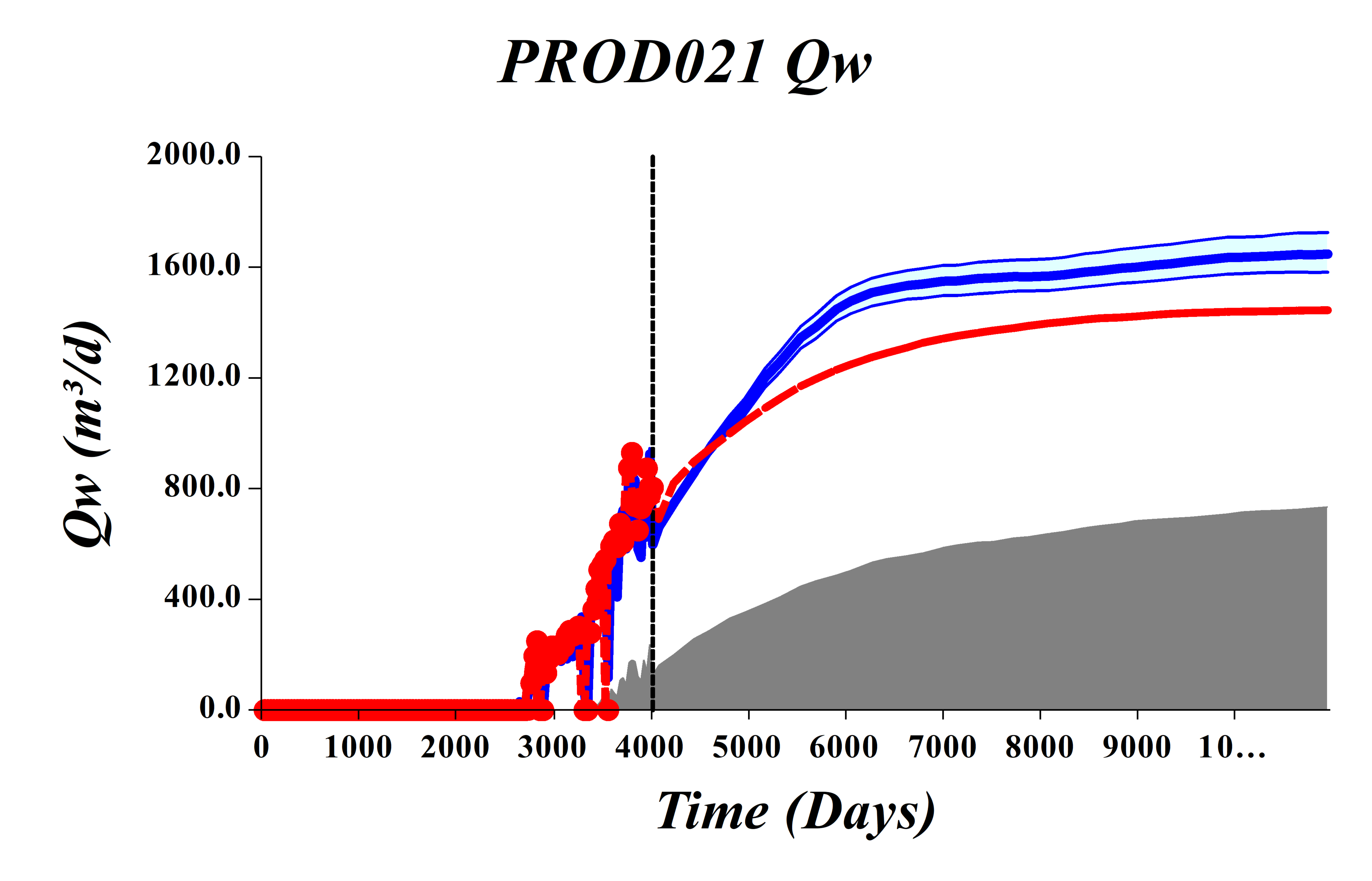}
    }
    \linebreak
    \subfloat[]{
      \includegraphics[width=0.33\textwidth]{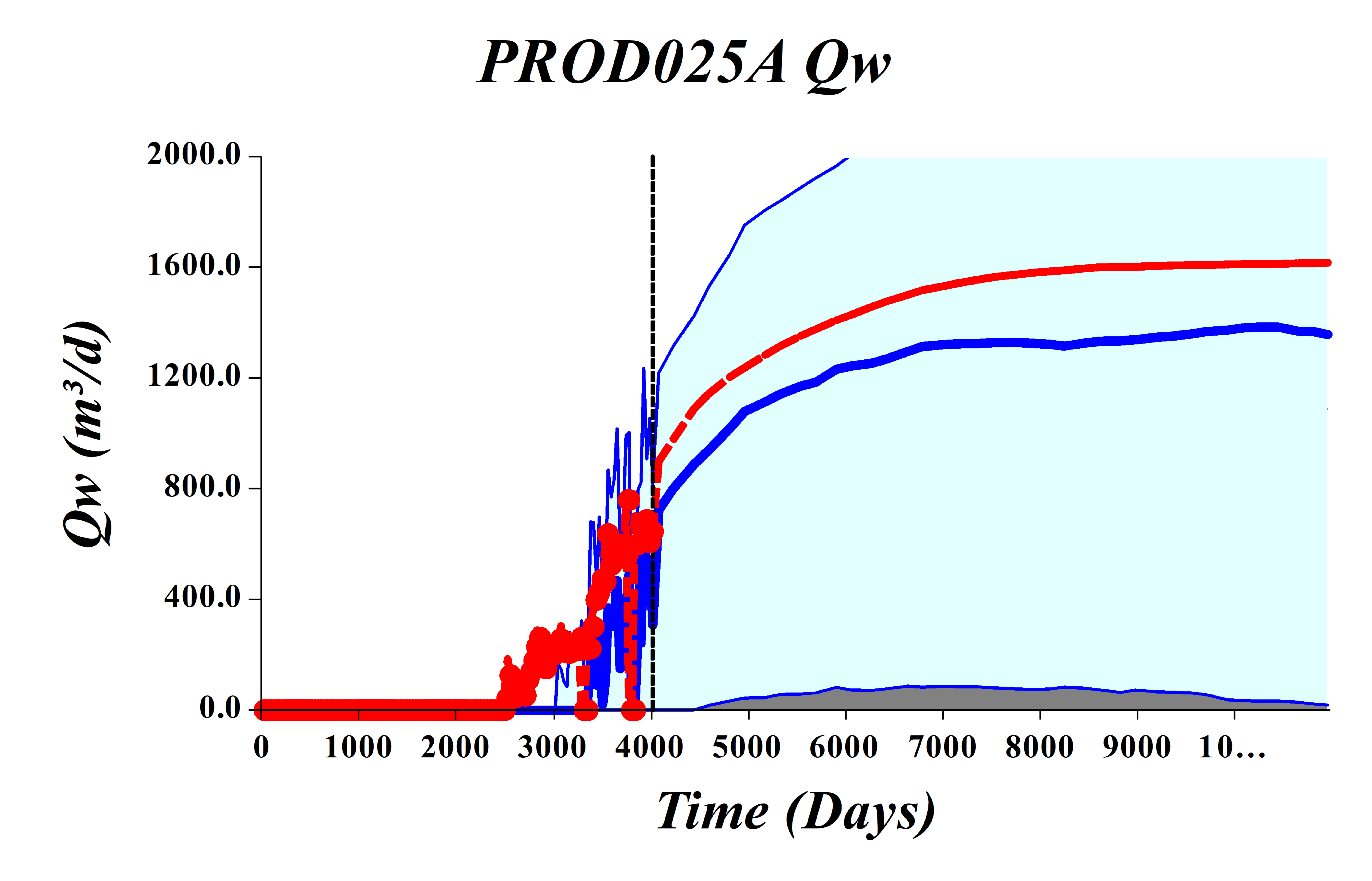}
    }
    \subfloat[]{
      \includegraphics[width=0.33\textwidth]{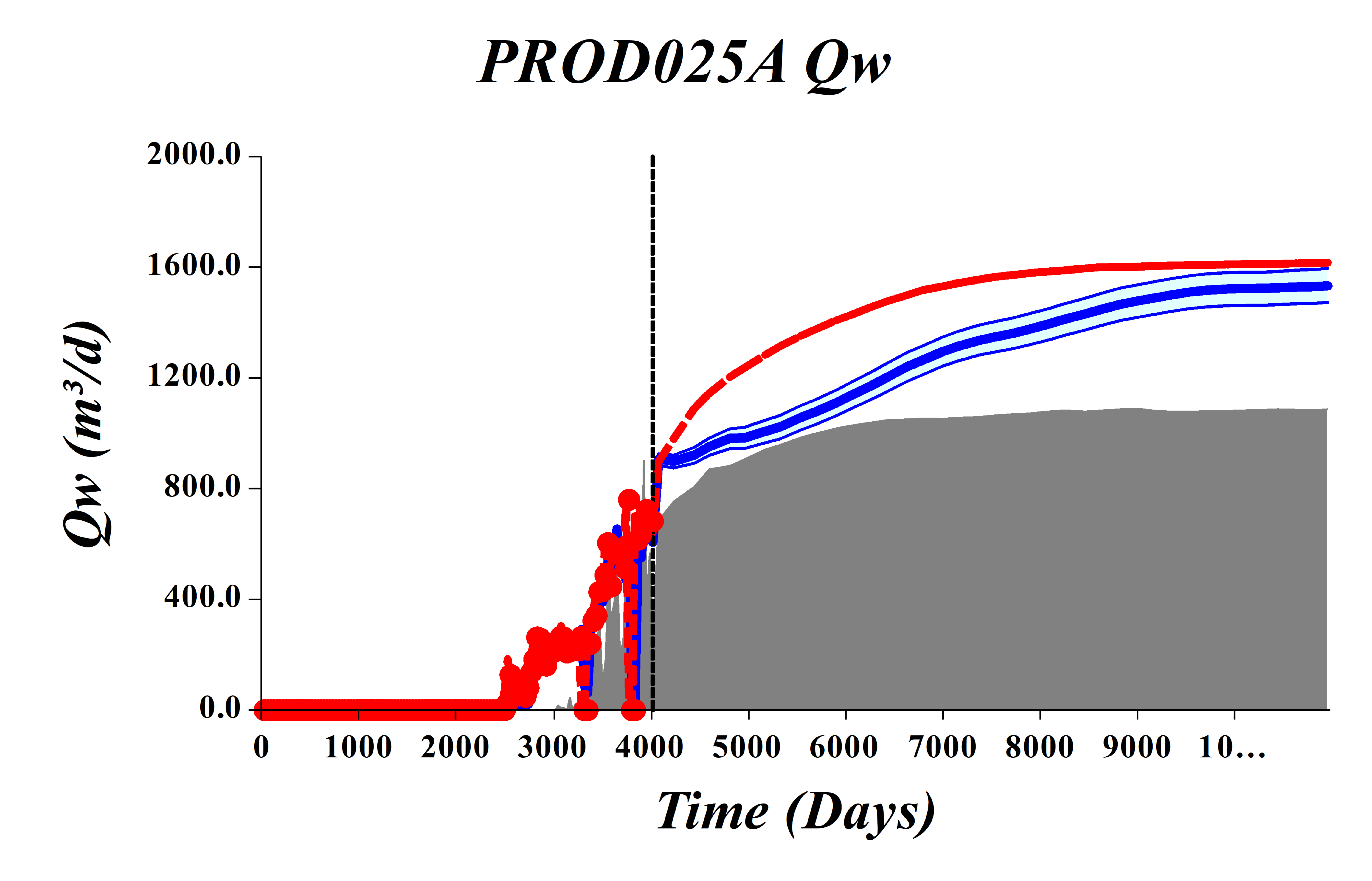}
    }
    \subfloat[]{
      \includegraphics[width=0.33\textwidth]{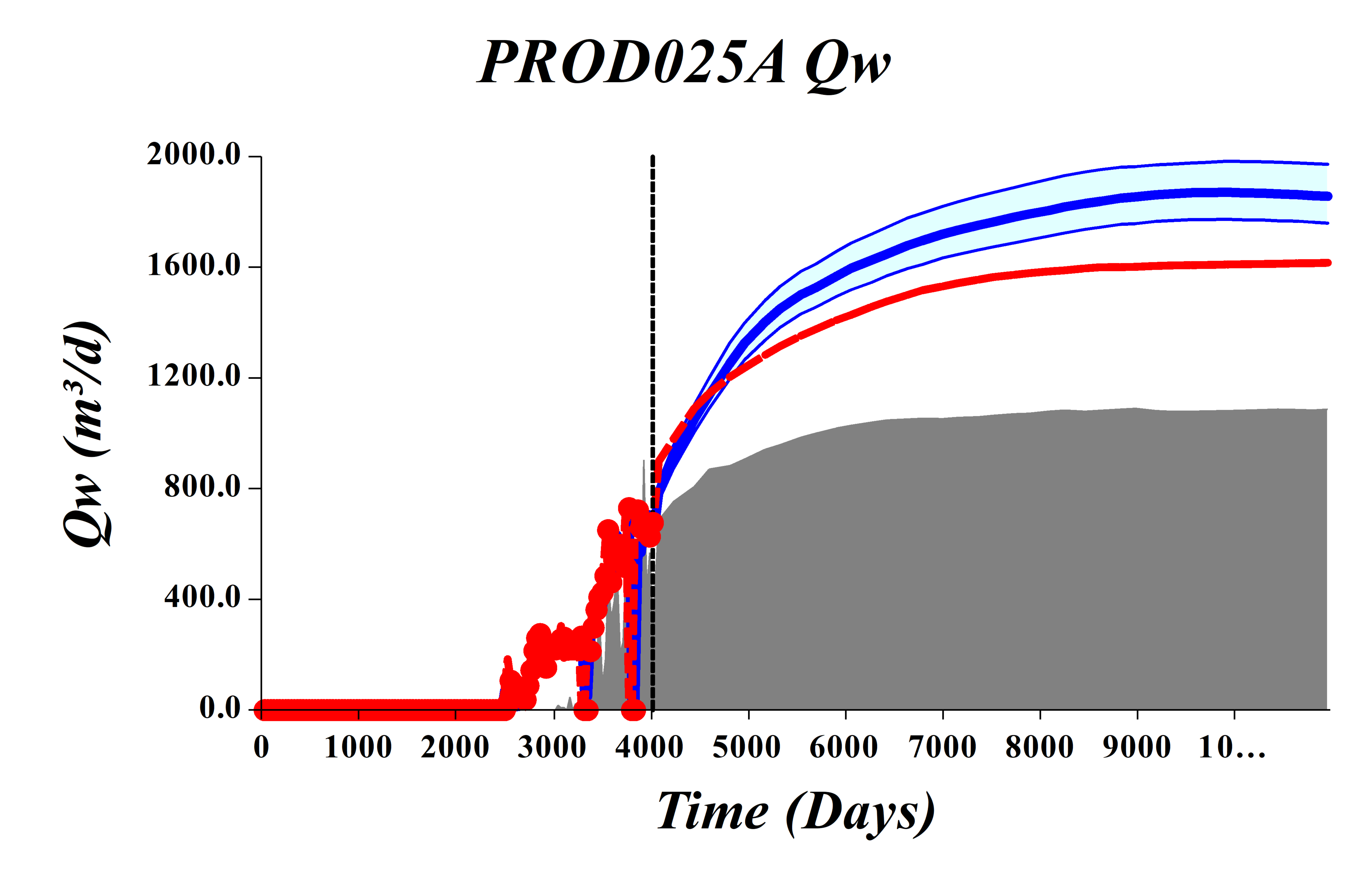}
    }
    \linebreak
    \subfloat[]{
      \includegraphics[width=0.33\textwidth]{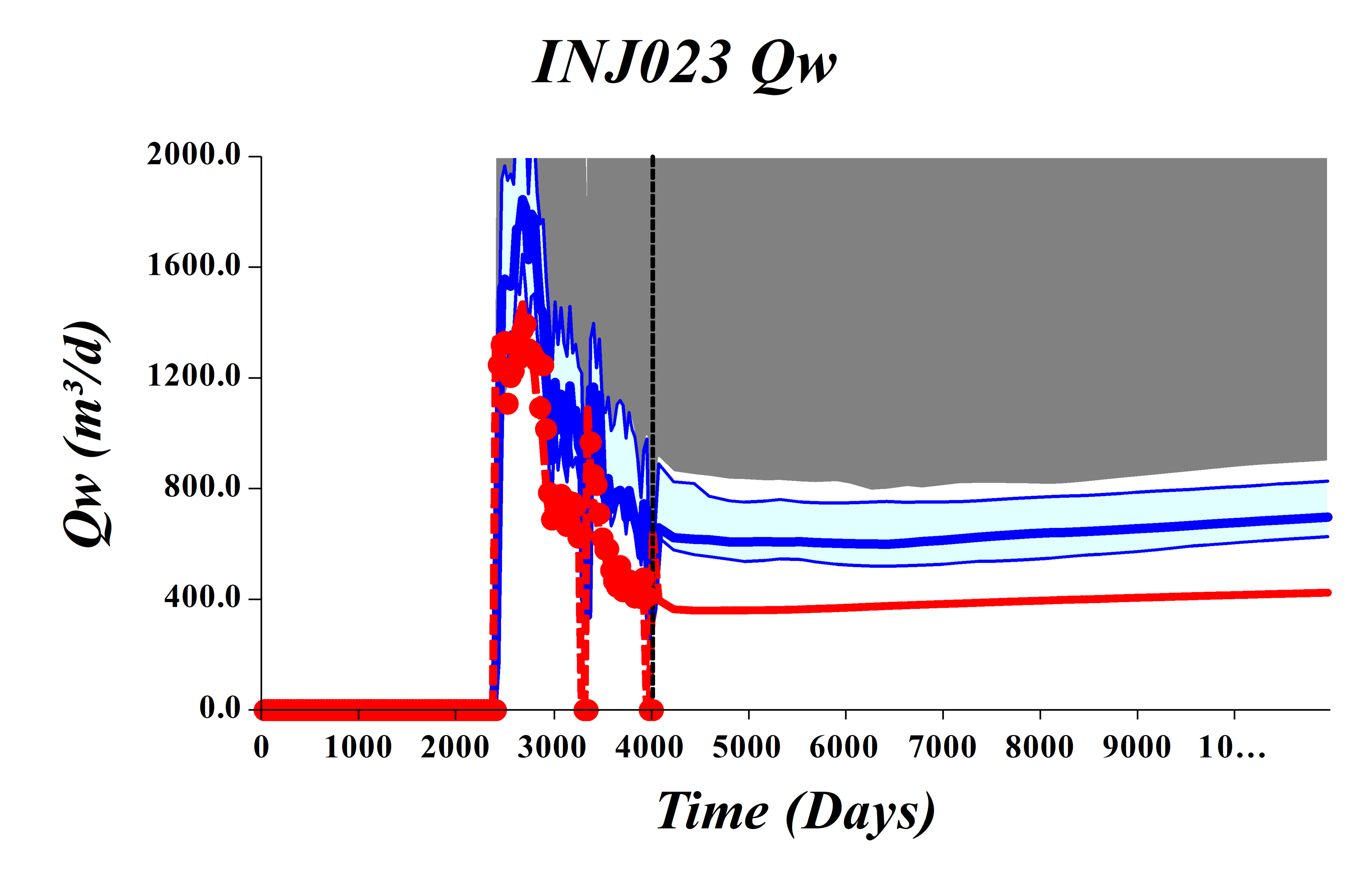}
    }
    \subfloat[]{
      \includegraphics[width=0.33\textwidth]{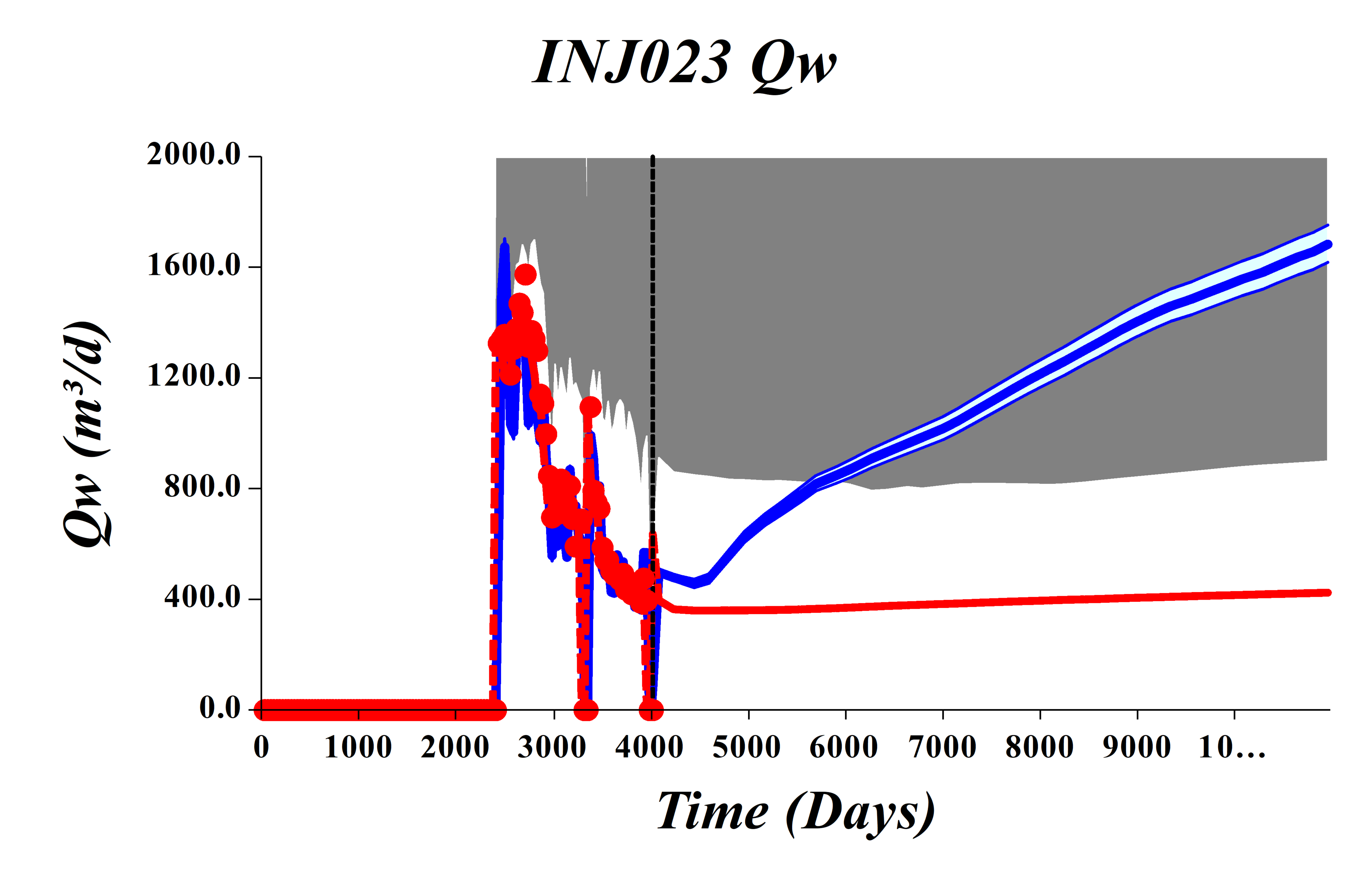}
    }
    \subfloat[]{
      \includegraphics[width=0.33\textwidth]{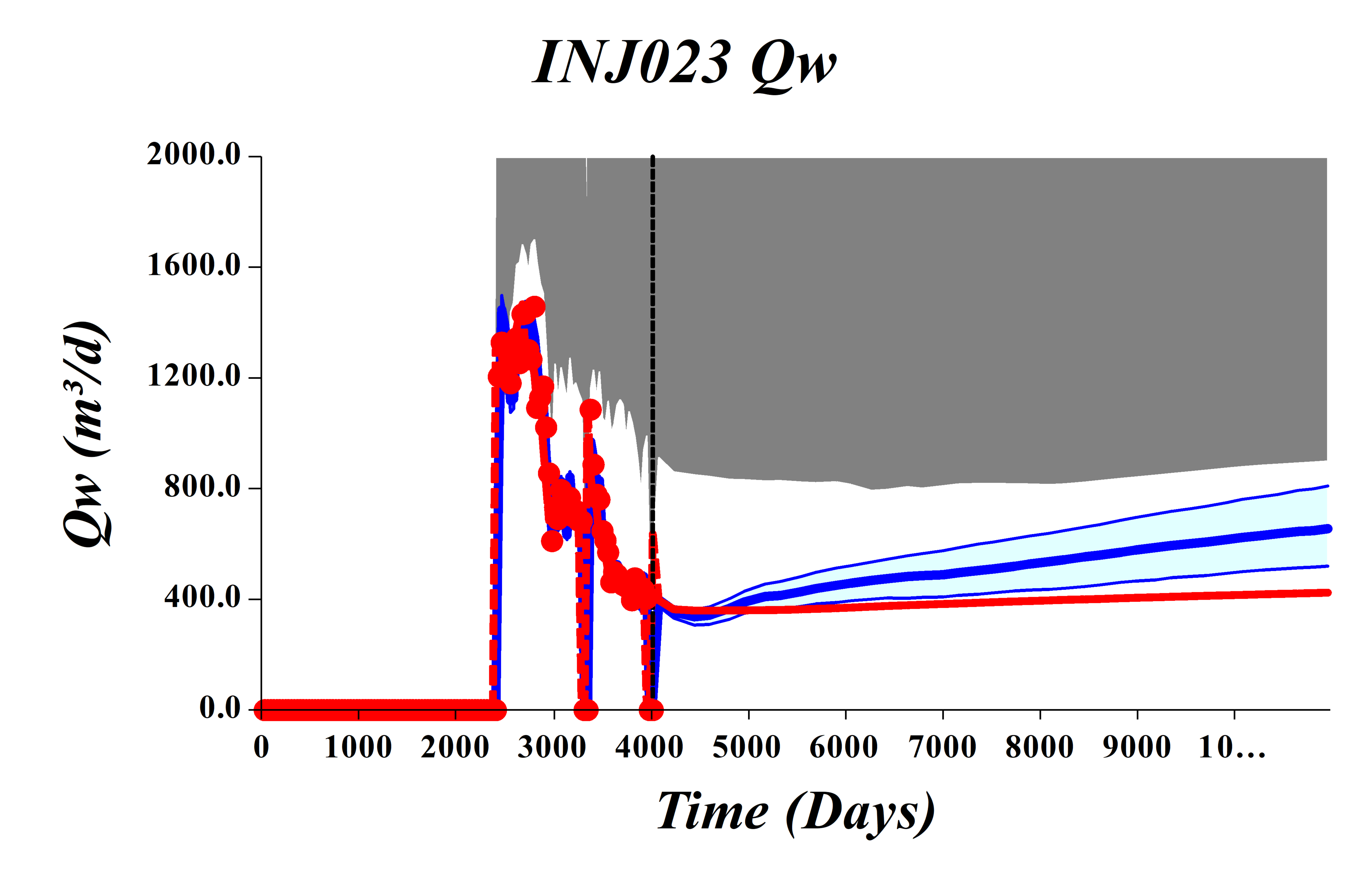}
    }
\captionsetup{justification=justified}
\caption{Water rate in m$^\text{3}$/days for three wells with poor coverage of predictions from the prior ensemble. Well PROD021 (first row), PROD025A (second row) and INJ023 (third row). (a), (d) and (g) DSI. (b), (e) and (h) DSI-ESMDA, (c), (f) and (i) DSI-ESMDA with localization. UNISIM-I-H. See Fig.~\ref{Fig:WaterRateTestCase1} for description.}
\label{Fig:WaterRateBadPriorUNISIM}
\end{figure}

\begin{figure}
\centering
    \captionsetup{justification=centering}
    \subfloat[]{
      \includegraphics[width=0.5\textwidth]{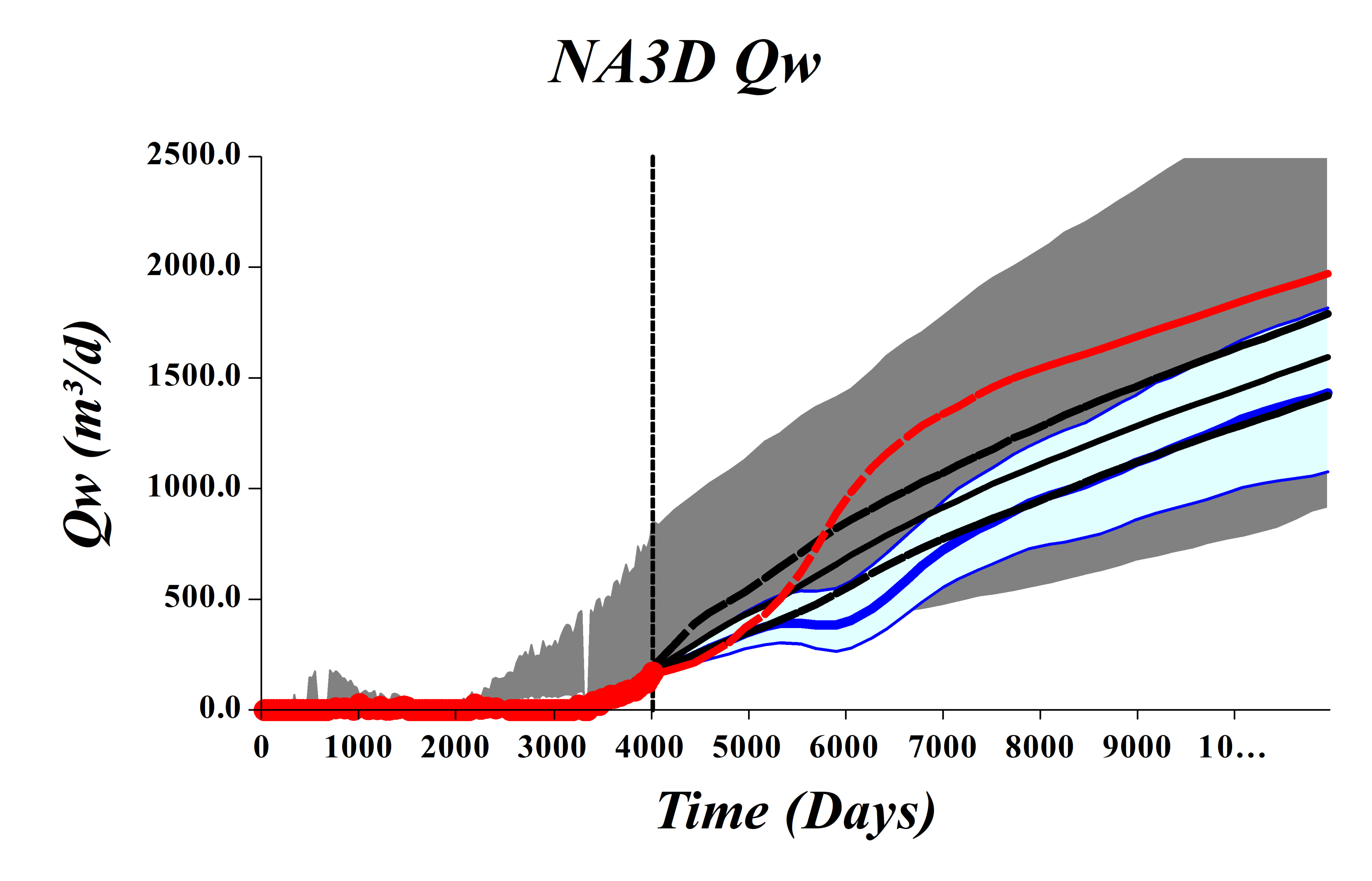}
    }
    \subfloat[]{
      \includegraphics[width=0.5\textwidth]{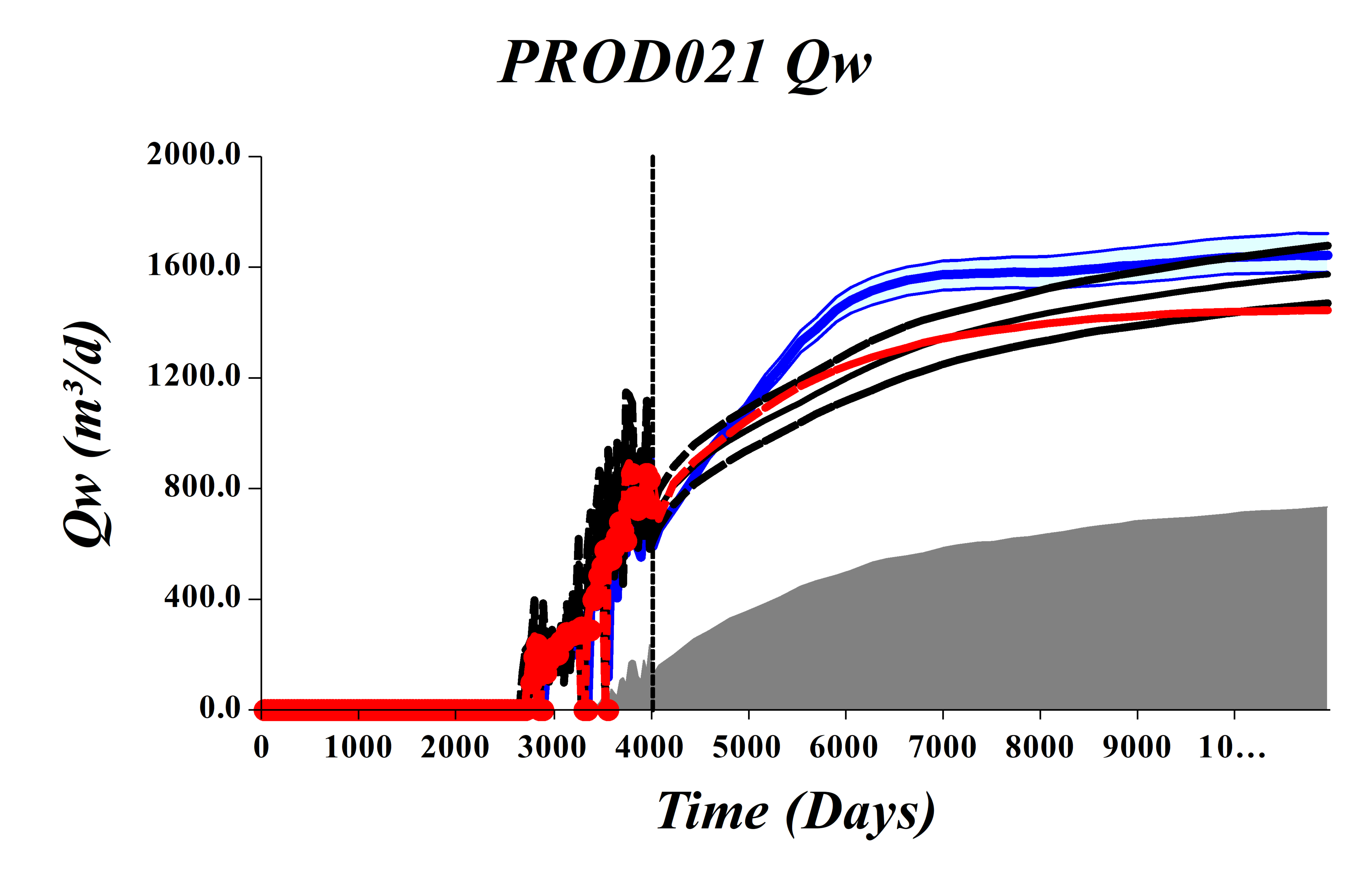}
    }
    \linebreak
    \subfloat[]{
      \includegraphics[width=0.5\textwidth]{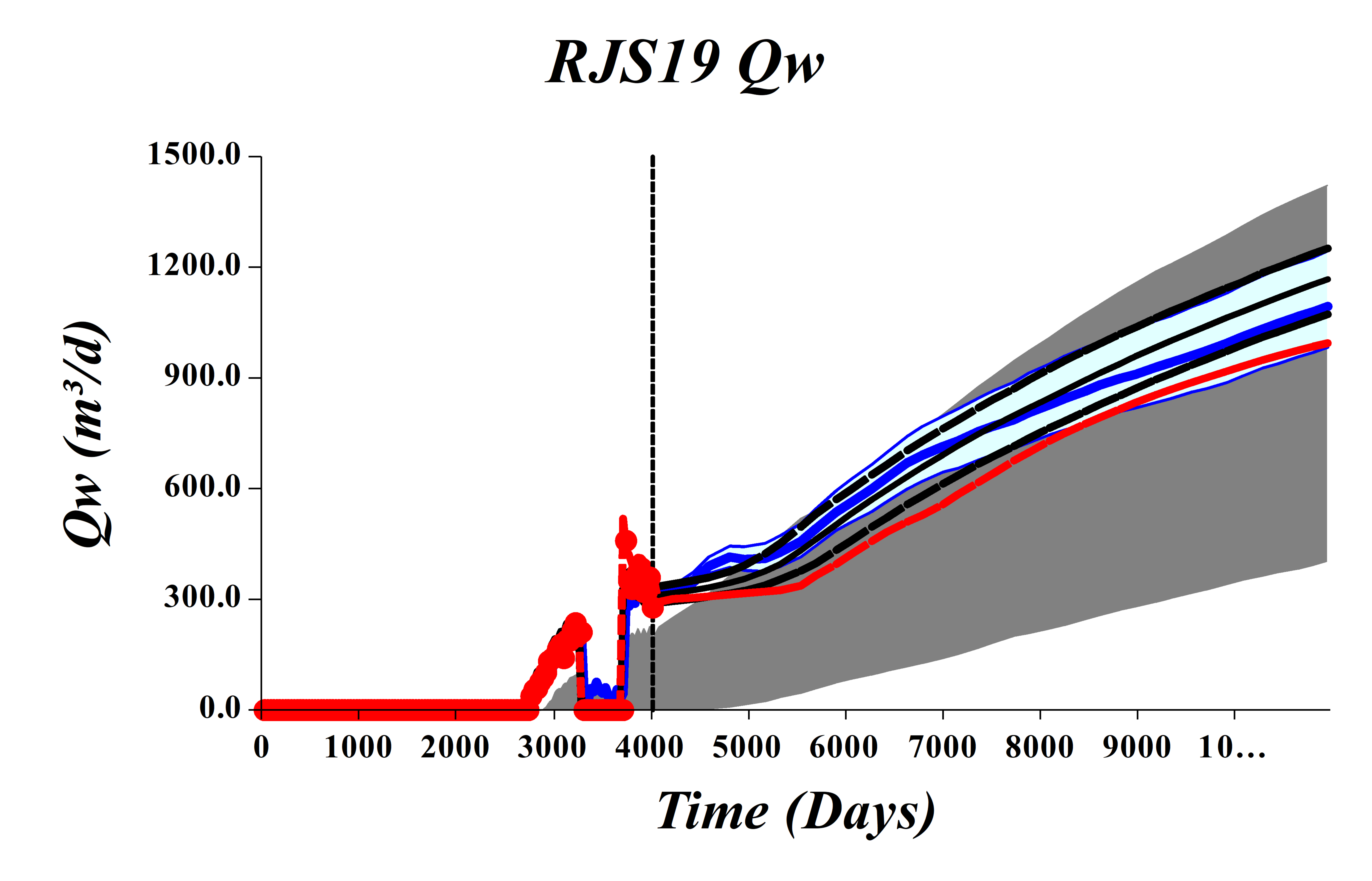}
    }
    \subfloat[]{
      \includegraphics[width=0.5\textwidth]{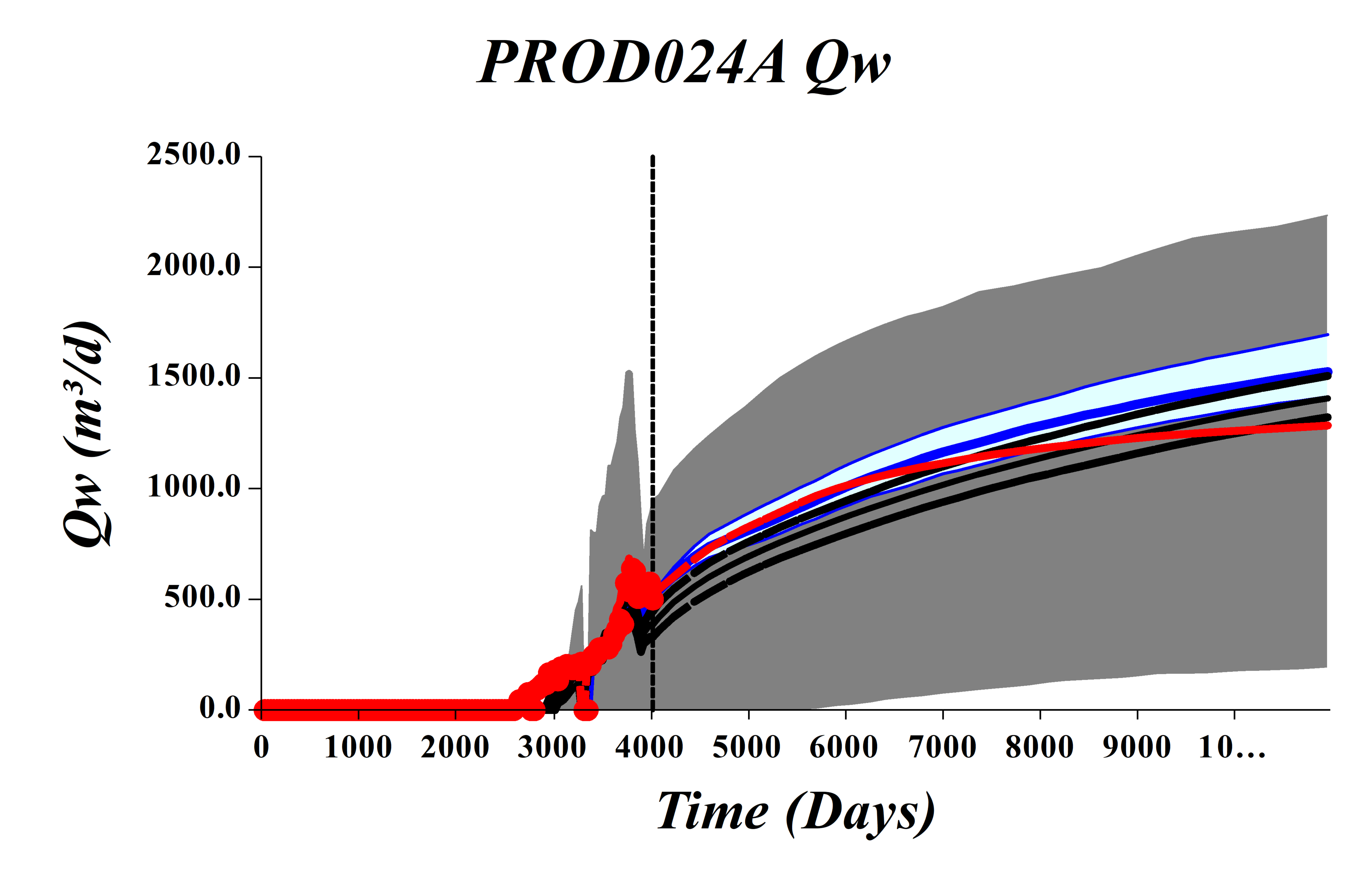}
    }
    \linebreak
    \subfloat[]{
      \includegraphics[width=0.5\textwidth]{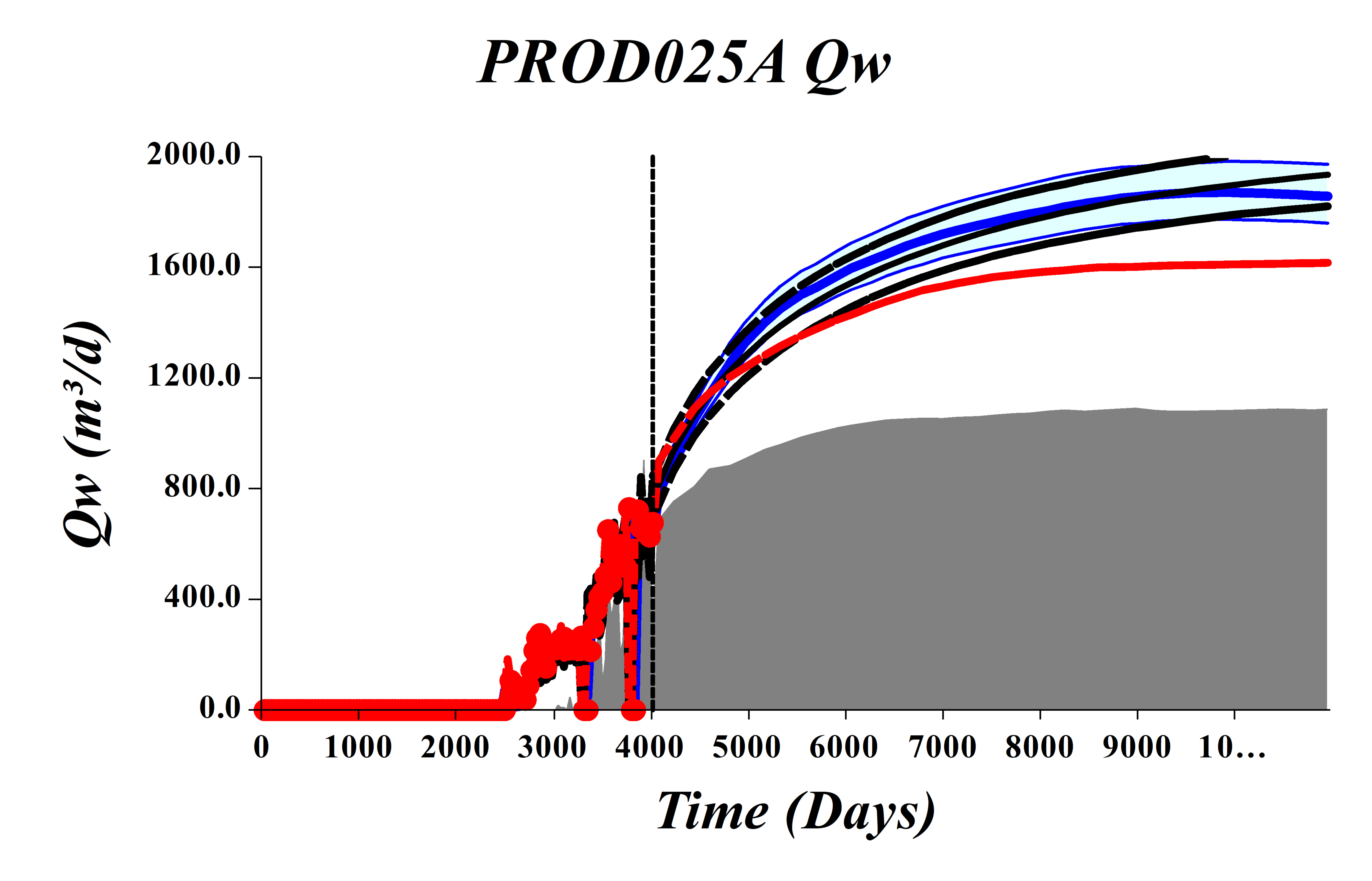}
    }
    \subfloat[]{
      \includegraphics[width=0.5\textwidth]{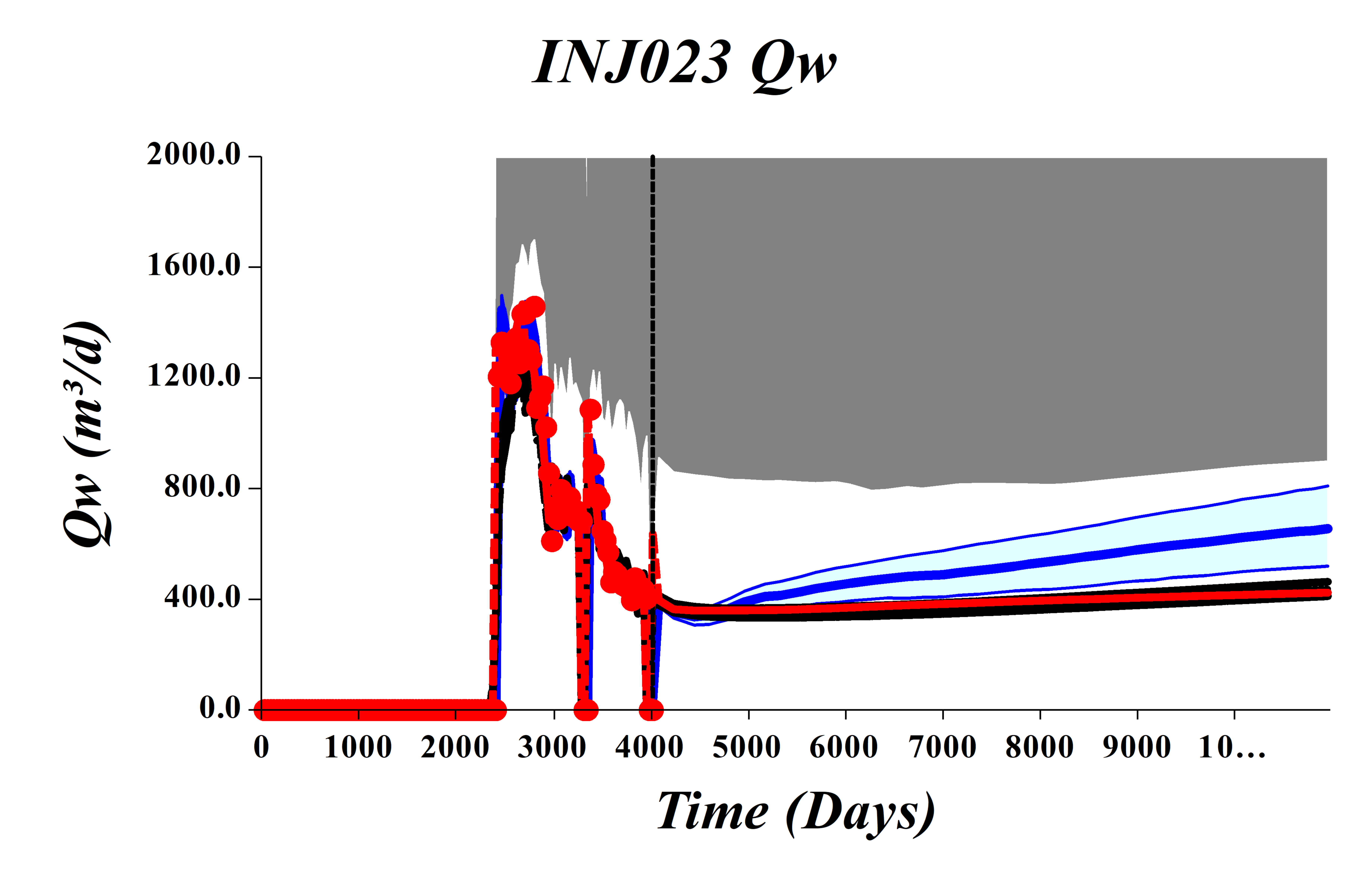}
    }
\captionsetup{justification=justified}
\caption{Water rate in m$^\text{3}$/days for six wells (NA3D, PROD021, RJS19, PROD024A, PROD025A, INJ023) obtained with DSI-ESMDA with localization and history matching with ES-MDA. UNISIM-I-H. See Fig.~\ref{Fig:WaterRateTestCase1} for description.}
\label{Fig:WaterRateESMDAUNISIM}
\end{figure}

\begin{figure}
\centering
    \captionsetup{justification=centering}
    \subfloat[]{
      \includegraphics[width=0.9\textwidth]{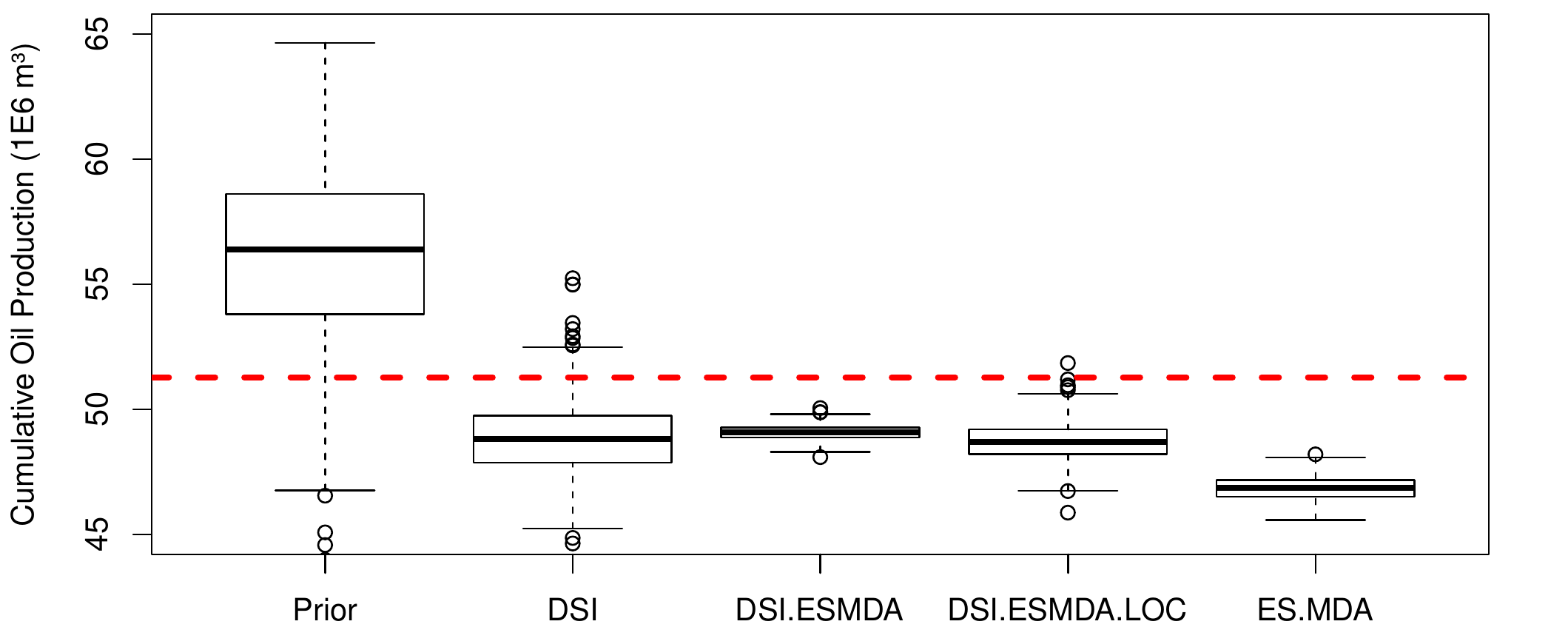}
    }
    \linebreak
        \subfloat[]{
      \includegraphics[width=0.9\textwidth]{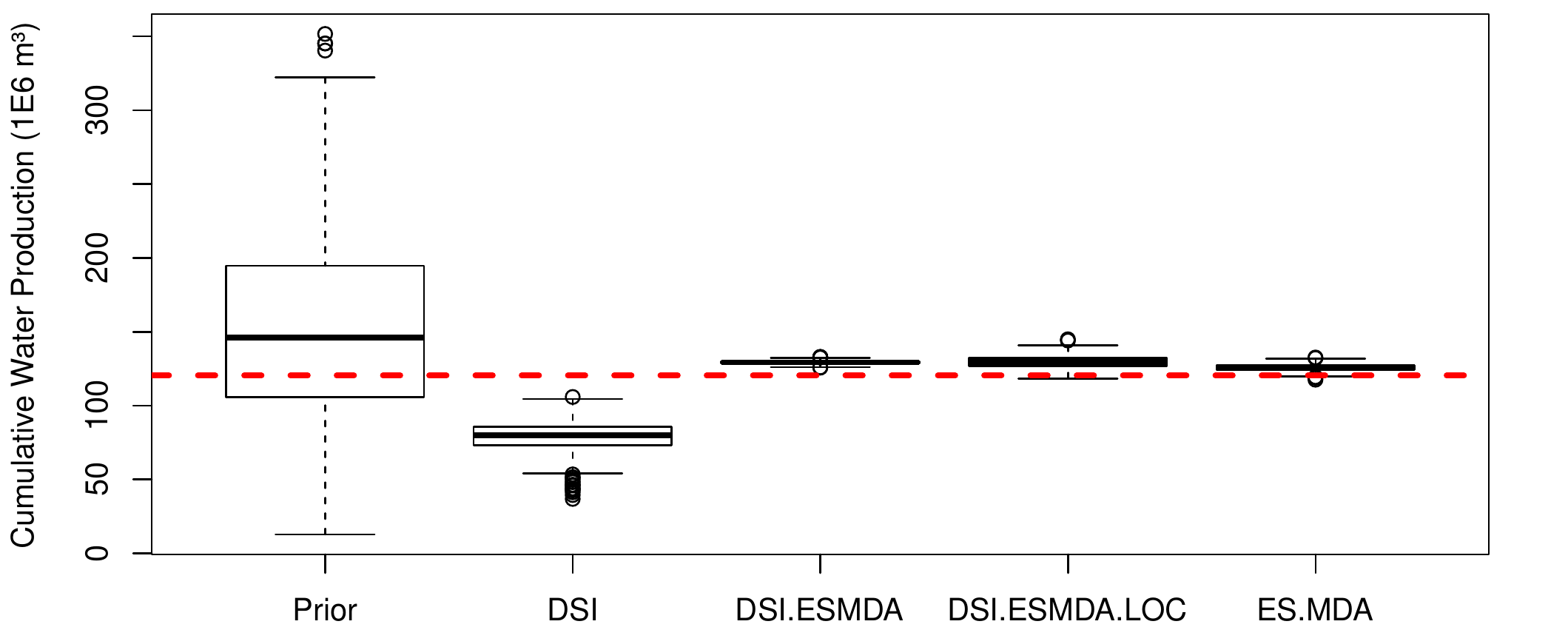}
    }
\captionsetup{justification=justified}
\caption{Field cumulative production in $\text{10}^\text{6}$~m$^\text{3}$. (a) Oil and (b) water. UNISIM-I-H. The dashed red line indicates the cumulative production of the reference case.}
\label{Fig:NpWpUNISIM}
\end{figure}

The prior ensemble for the UNISIM-I-H case was provided by the authors of the benchmark and it is noteworthy that even though this is a synthetic problem, the predictions from the prior ensemble do not span the data from the reference case for several wells. This is clearly not an ideal situation for applying DSI or even any history-matching method. Ideally, the prior ensemble should provide a reasonable estimate of the prior uncertainty, at the very least, it should to be able to span the observations. In practice, we should revise the prior ensemble before using for data assimilation. However, here for the purposes of the investigation, we decided to test the methods with this deficient prior ensemble. Figures~\ref{Fig:WaterRateGoodPriorUNISIM} and \ref{Fig:WaterRateBadPriorUNISIM} show the results obtained by DSI, DSI-ESMDA and DSI-ESMDA with localization for three wells with good and poor coverage of the predictions from the prior ensemble, respectively. We selected these wells because they are representative of the results observed in this problem. Figure~\ref{Fig:WaterRateGoodPriorUNISIM} indicates that DSI failed to match data even for the wells with good prior coverage. As a result, the prediction from the reference model lies outside the predicted P10--P90 range. For the wells with poor coverage (Fig.~\ref{Fig:WaterRateBadPriorUNISIM}), DSI seems to fail to reduce the uncertainty range properly, for example, the posterior predictions for wells PROD021 and PROD025A have almost the same uncertainty range of the prior ones. DSI-ESMDA obtained reasonable data matches for all wells, however, the predictions are clearly too narrow as they do not span the reference. DSI-ESMDA with localization improved the results significantly, although it was not able to span the reference predictions, especially for the wells with poor prior coverage. Figure~\ref{Fig:WaterRateESMDAUNISIM} compares the results of DSI-ESMDA with localization with the history matching with ES-MDA. It is interesting to note that the history matched models also suffer from the lack of representativeness of the prior ensemble; see, for example, the plots (a) and (d) in Fig.~\ref{Fig:WaterRateESMDAUNISIM}. Figure~\ref{Fig:NpWpUNISIM} shows boxplots of field cumulative oil and water production predicted by the prior ensemble and by each method. The cumulative production from the reference case is also included in this figure for comparisons. Overall, the predictions are biased compared to the reference production, which is probably explained by the problems with the prior ensemble.

\clearpage

\section{Field Case}
\label{Sec:FieldCases}

\begin{figure}
\centering
	\includegraphics[width=0.8\linewidth]{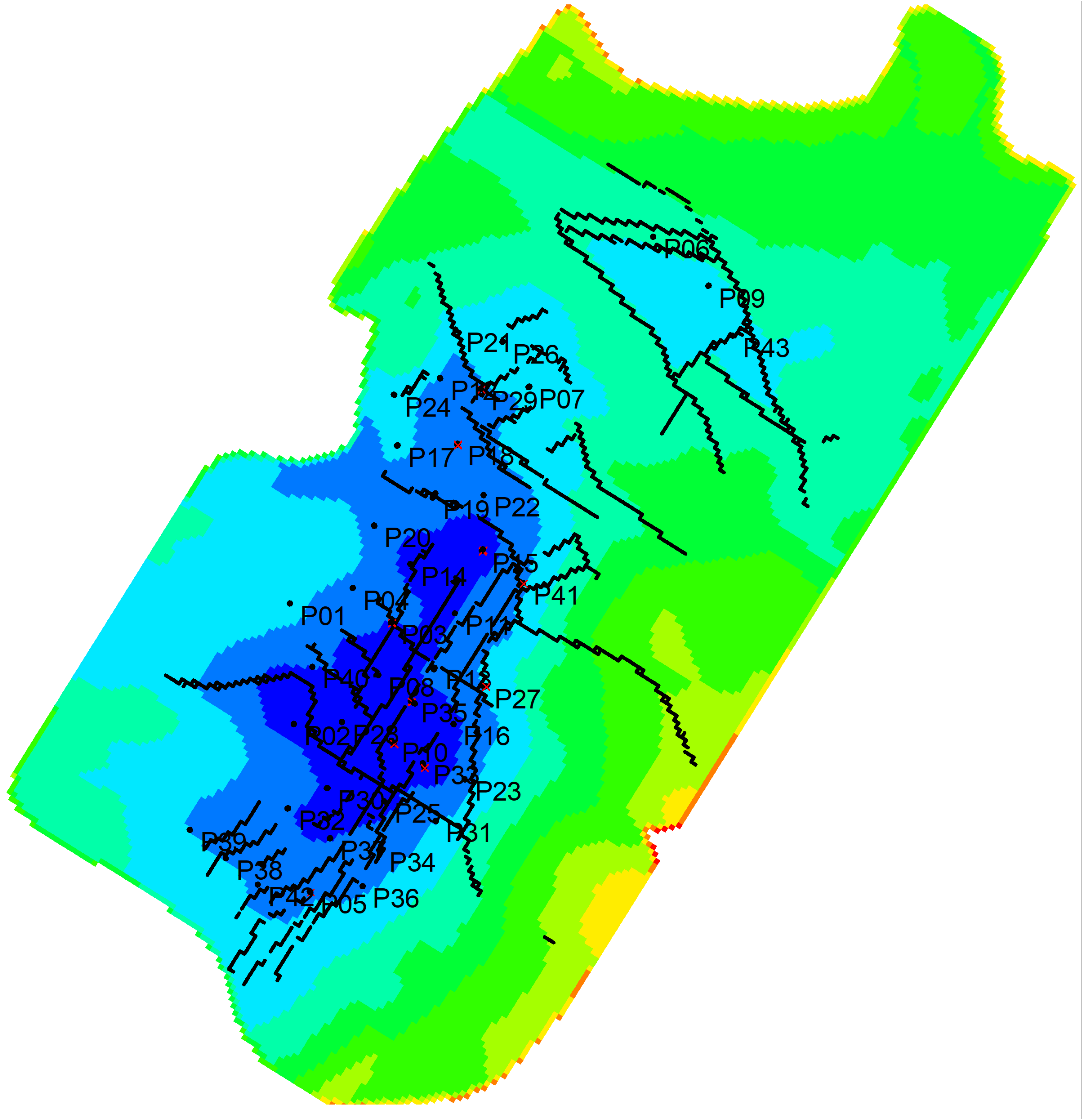}
\caption{Top view of the field case model showing the large number of wells and faults.}
\label{Fig:FieldCaseModel}
\end{figure}

The field case corresponds to a large offshore turbidite reservoir in Campos Basis. Currently there are 24~wells producing oil in the field from a total of 43 vertical wells that have been drilled during 18~years of operation. The main recovery mechanism is pressure maintenance by a large aquifer. A total of 500 prior realizations of the reservoir were built by the asset team of the field using the best practices available to date in the Company. These realization include facies, porosity, permeability and net-to-gross ratio for each one of 11 producing zones. The engineers of the field also selected the oil-water relative permeability curves as parameters with relevant uncertainty. The relative permeability curves in six different regions of the model were parameterized using Corey functions with the exponents and the maximum relative permeability of water selected as uncertainty parameters with Gaussian prior distributions. The model has approximately 500,000 active gridblocks (total of 114 $\times$ 238 $\times$ 60  with an average size of 75~m $\times$ 75~m $\times$ 2~m). This field has a large number of sealing faults as illustrated in Fig.~\ref{Fig:FieldCaseModel}. All wells are controlled by specified total liquid production rate during the historical period and the forecast. The observed data corresponds to monthly measurements of oil and water rate at the well with a total of 26703 data points.

\begin{figure}
\centering
    \captionsetup{justification=centering}
    \subfloat[]{
      \includegraphics[width=0.47\textwidth]{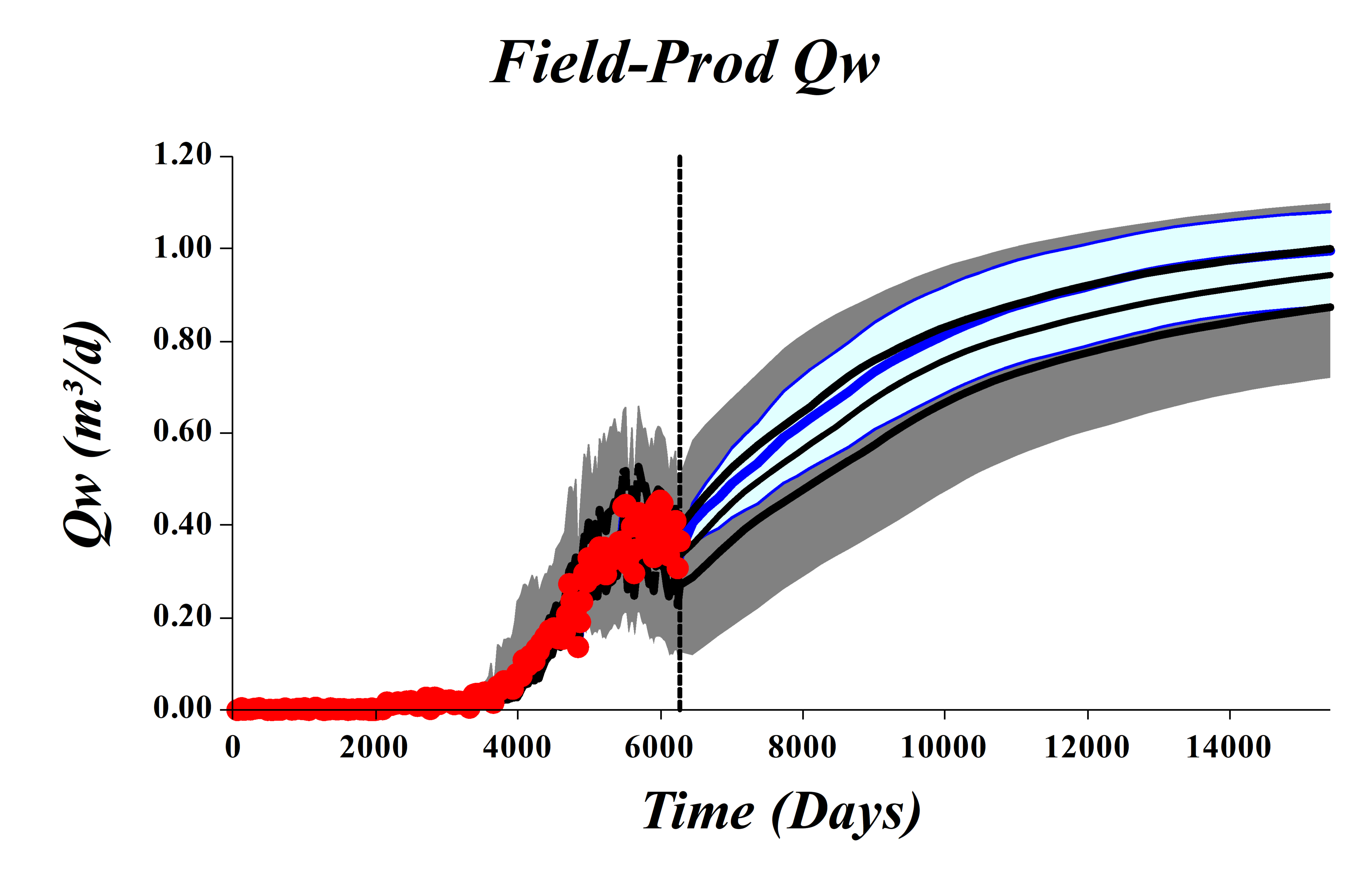}
    }
    \subfloat[]{
      \includegraphics[width=0.47\textwidth]{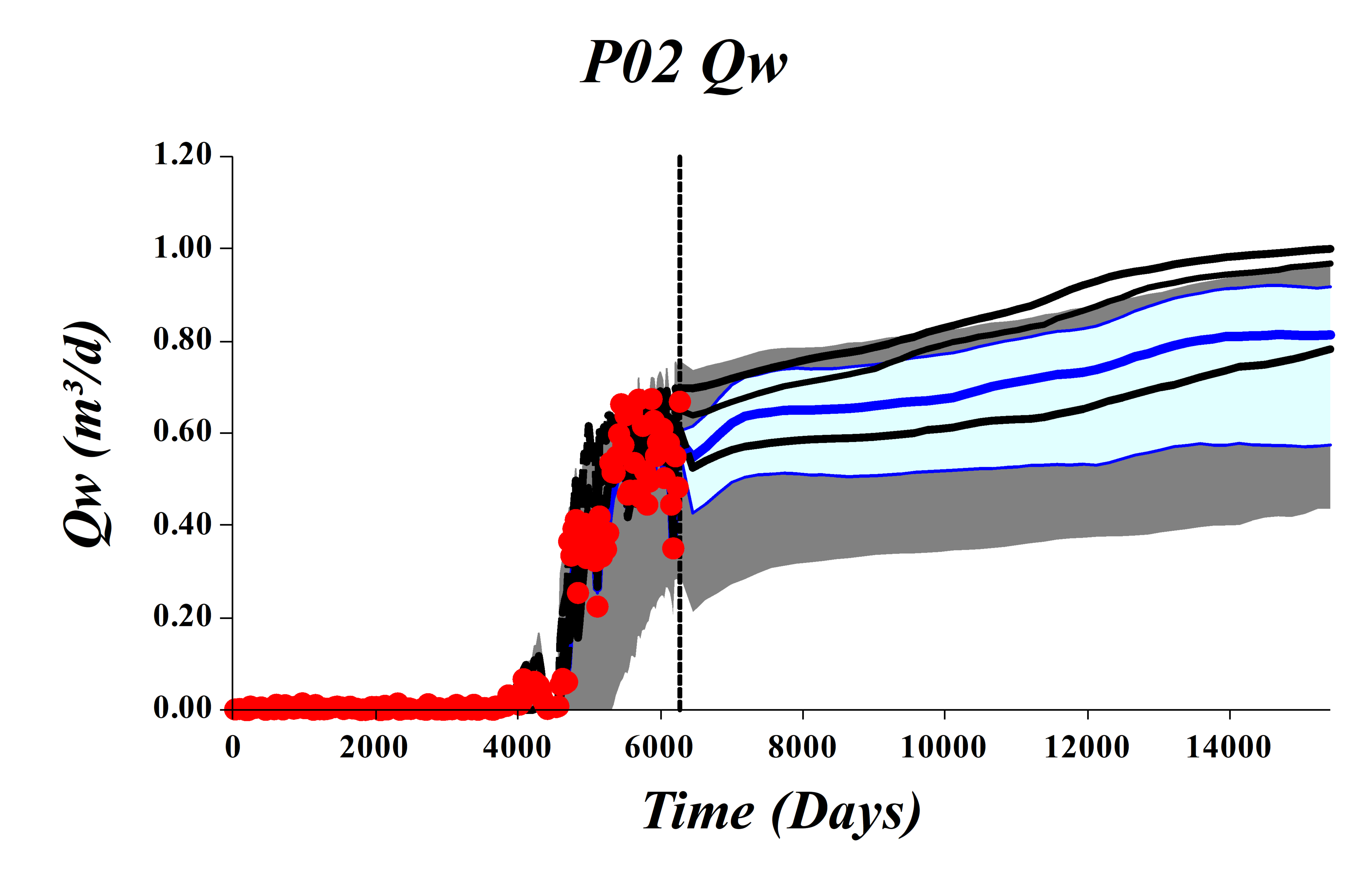}
    }
    \linebreak
    \subfloat[]{
      \includegraphics[width=0.47\textwidth]{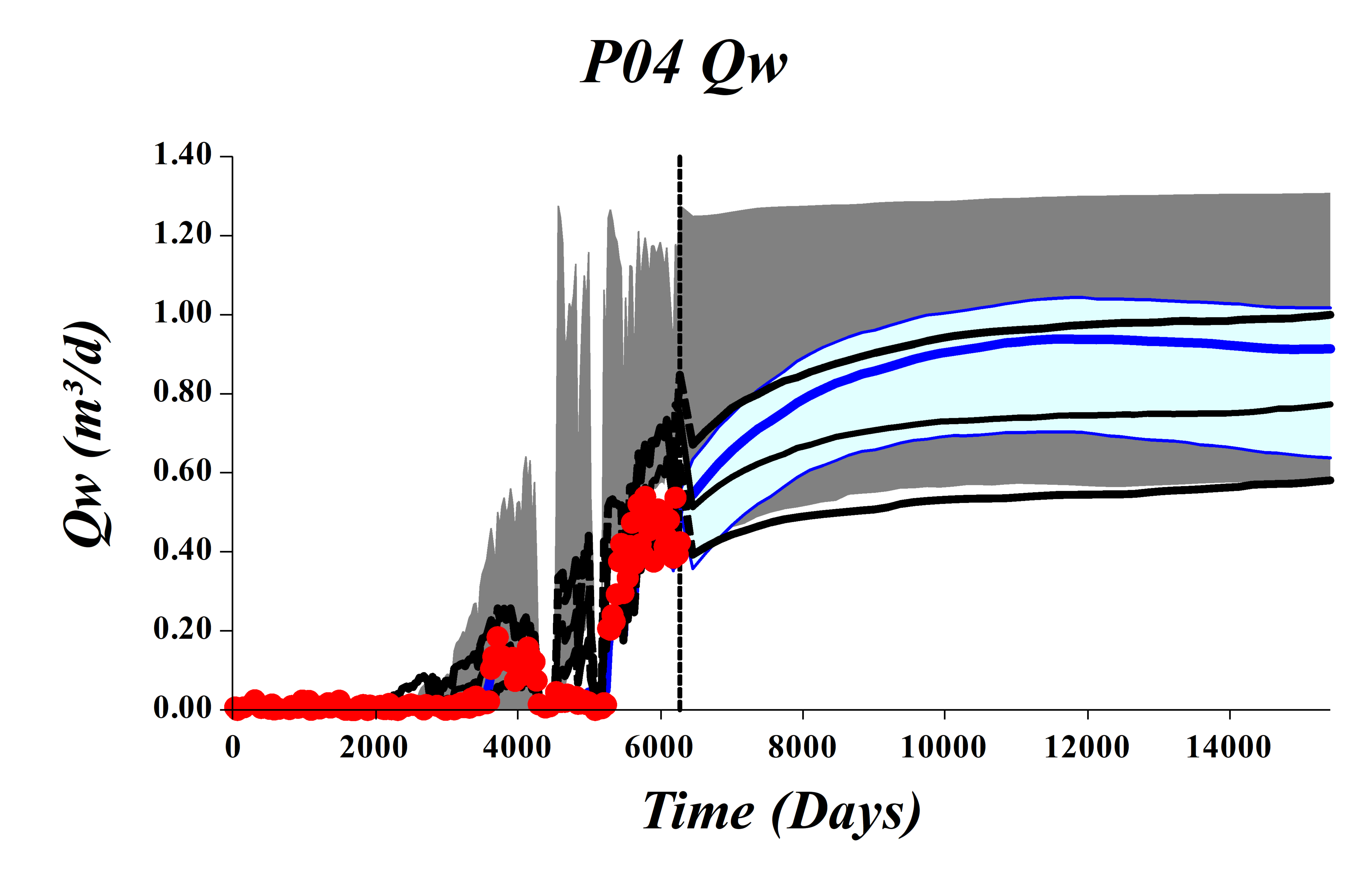}
    }
    \subfloat[]{
      \includegraphics[width=0.47\textwidth]{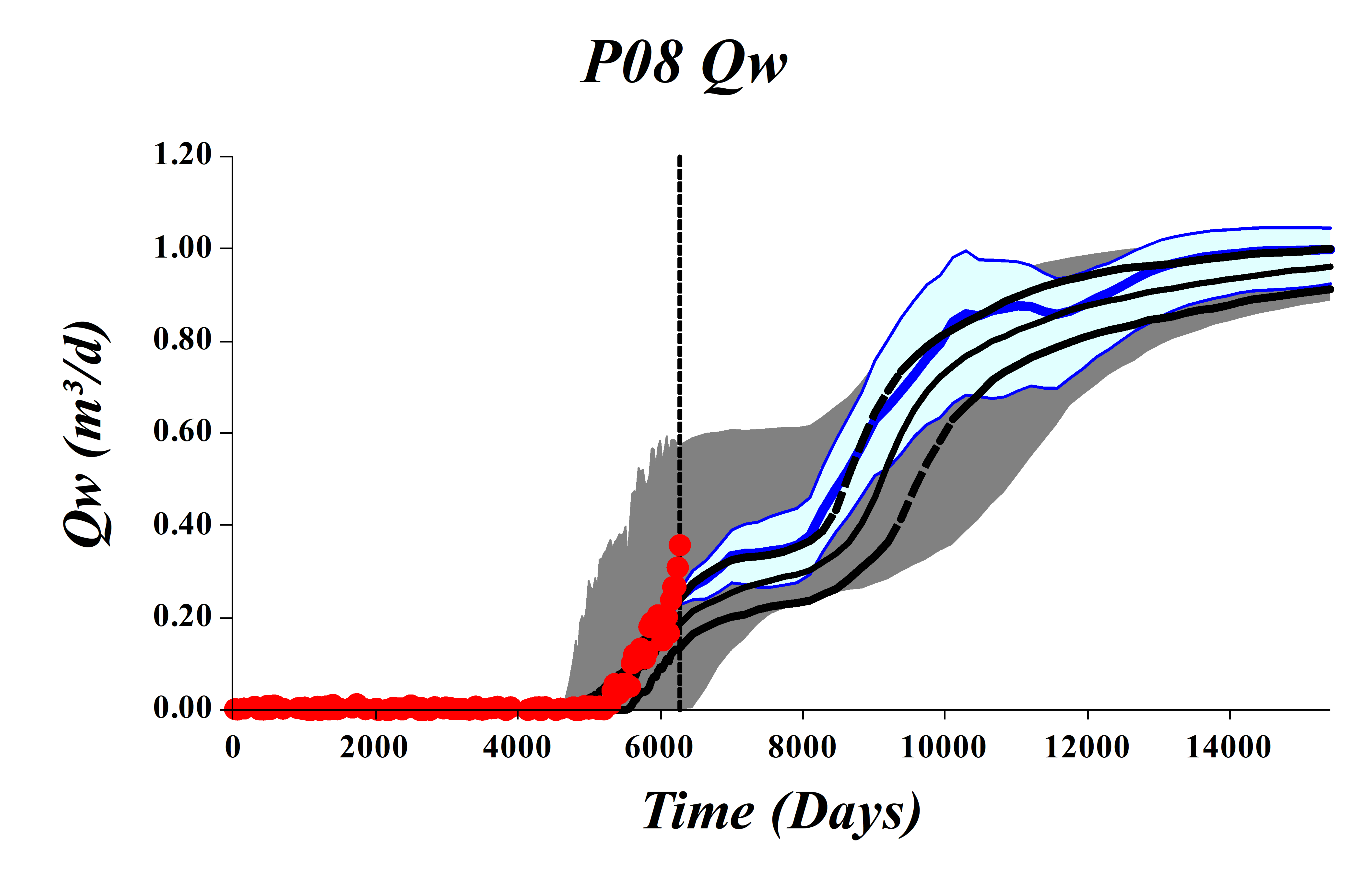}
    }
    \linebreak
    \subfloat[]{
      \includegraphics[width=0.47\textwidth]{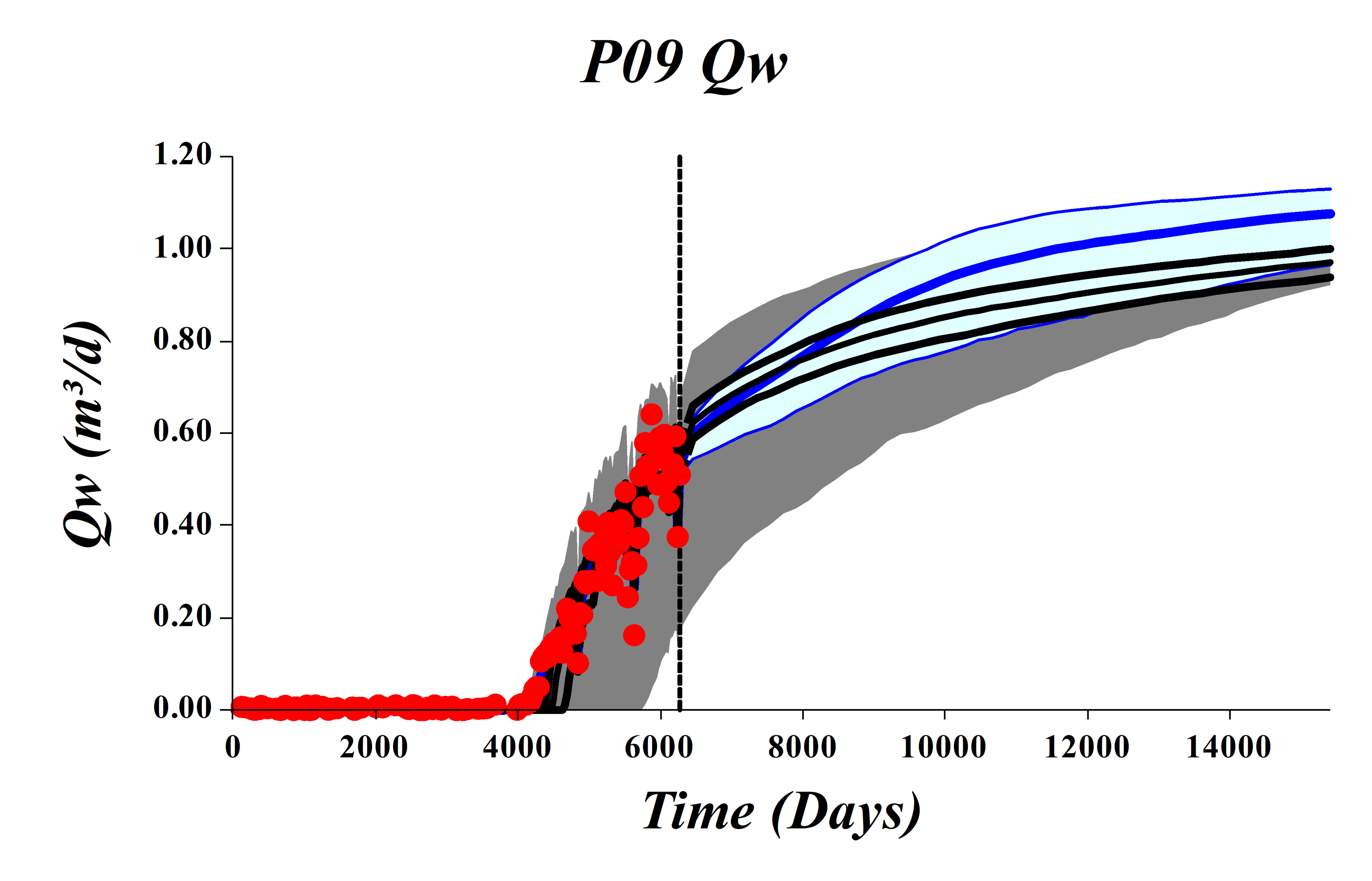}
    }
    \subfloat[]{
      \includegraphics[width=0.47\textwidth]{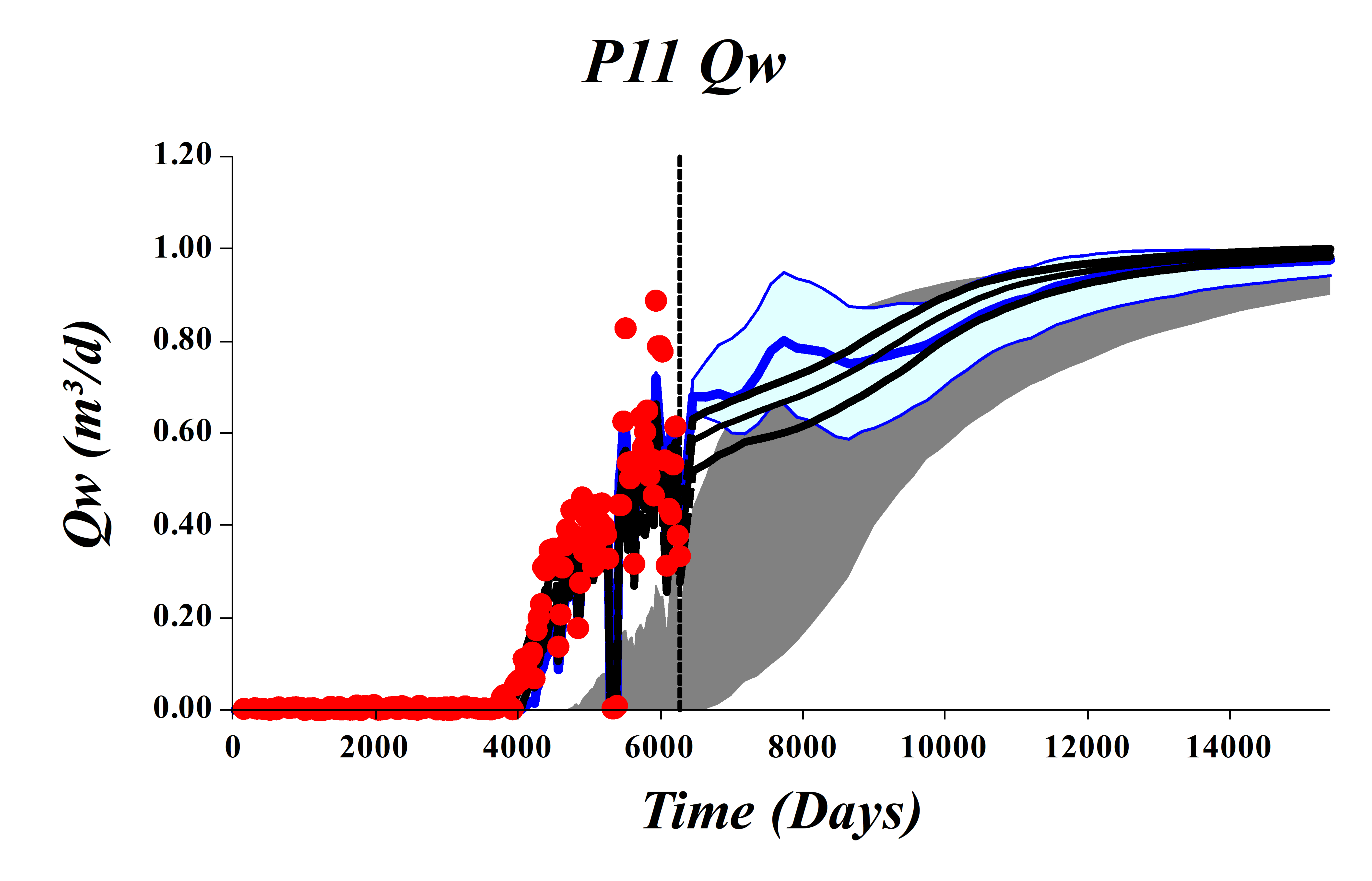}
    }
    \linebreak
    \subfloat[]{
      \includegraphics[width=0.47\textwidth]{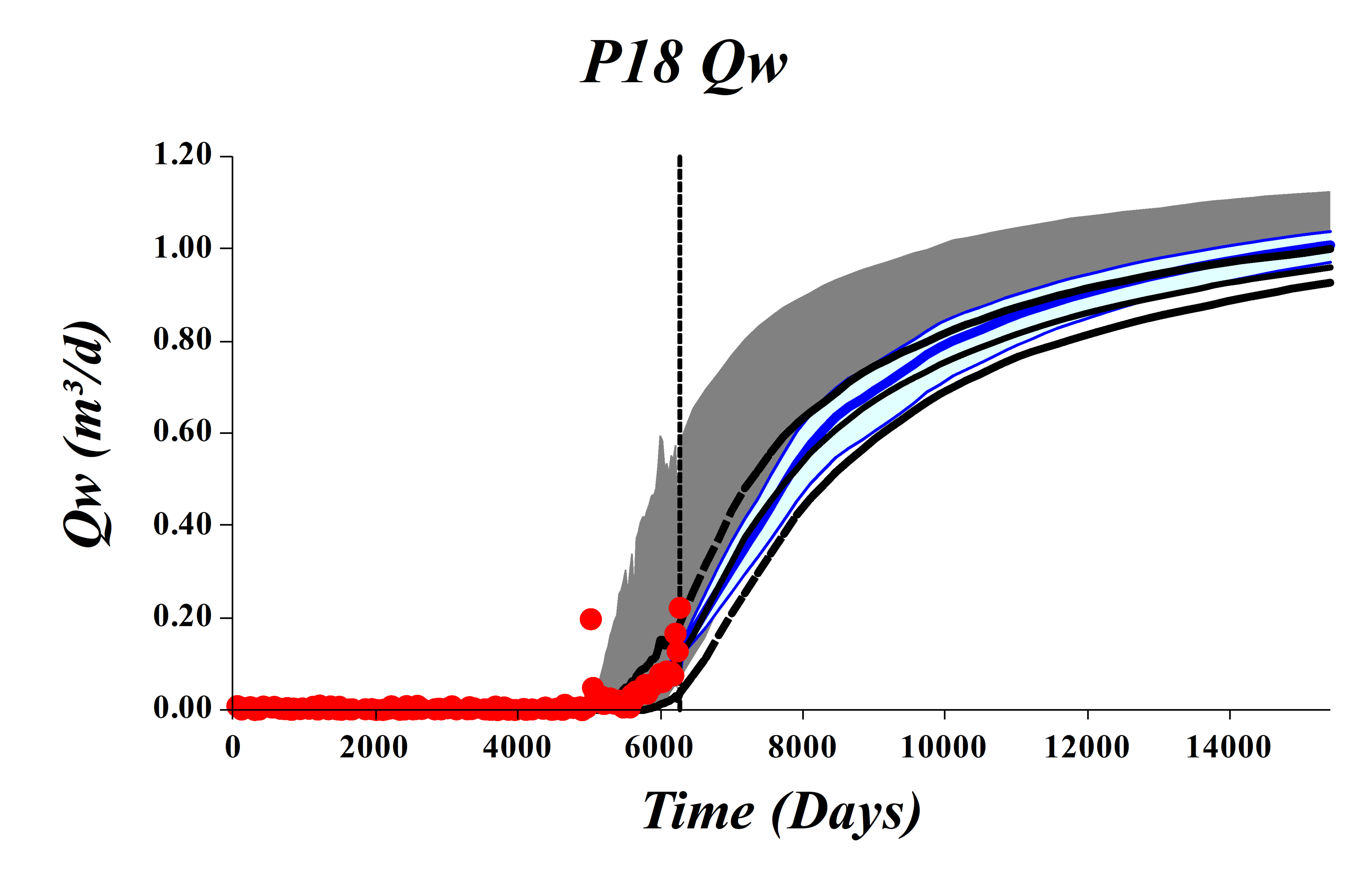}
    }
    \subfloat[]{
      \includegraphics[width=0.47\textwidth]{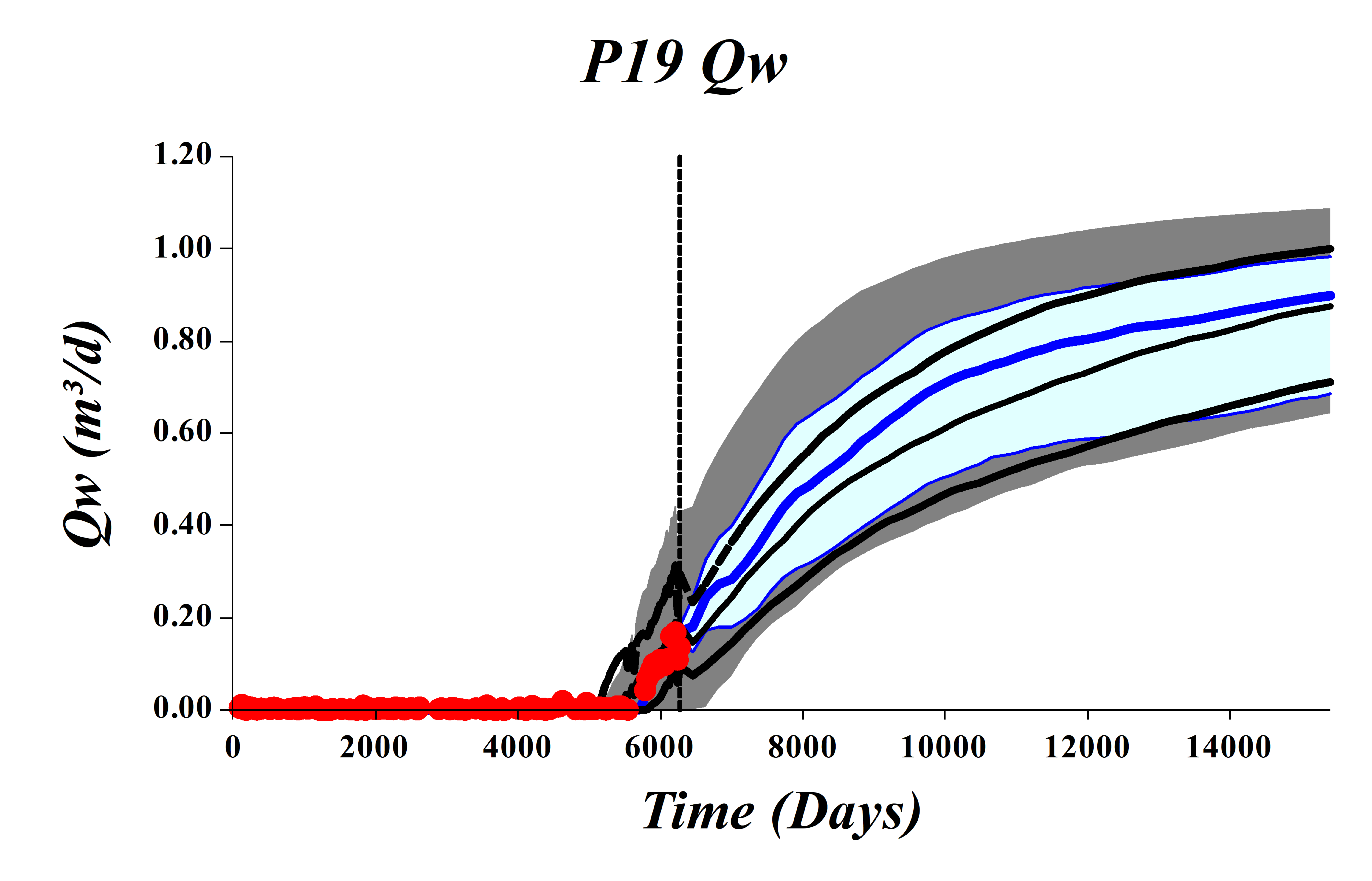}
    }
\captionsetup{justification=justified}
\caption{Normalized water production rate obtained with DSI-ESMDA with localization. (a) Total field production, (b)--(h) seven selected wells. Field case. See Fig.~\ref{Fig:WaterRateTestCase1} for description.}
\label{Fig:WaterRateFieldCase}
\end{figure}

\begin{figure}
\centering
    \captionsetup{justification=centering}
    \subfloat[]{
      \includegraphics[width=0.5\textwidth]{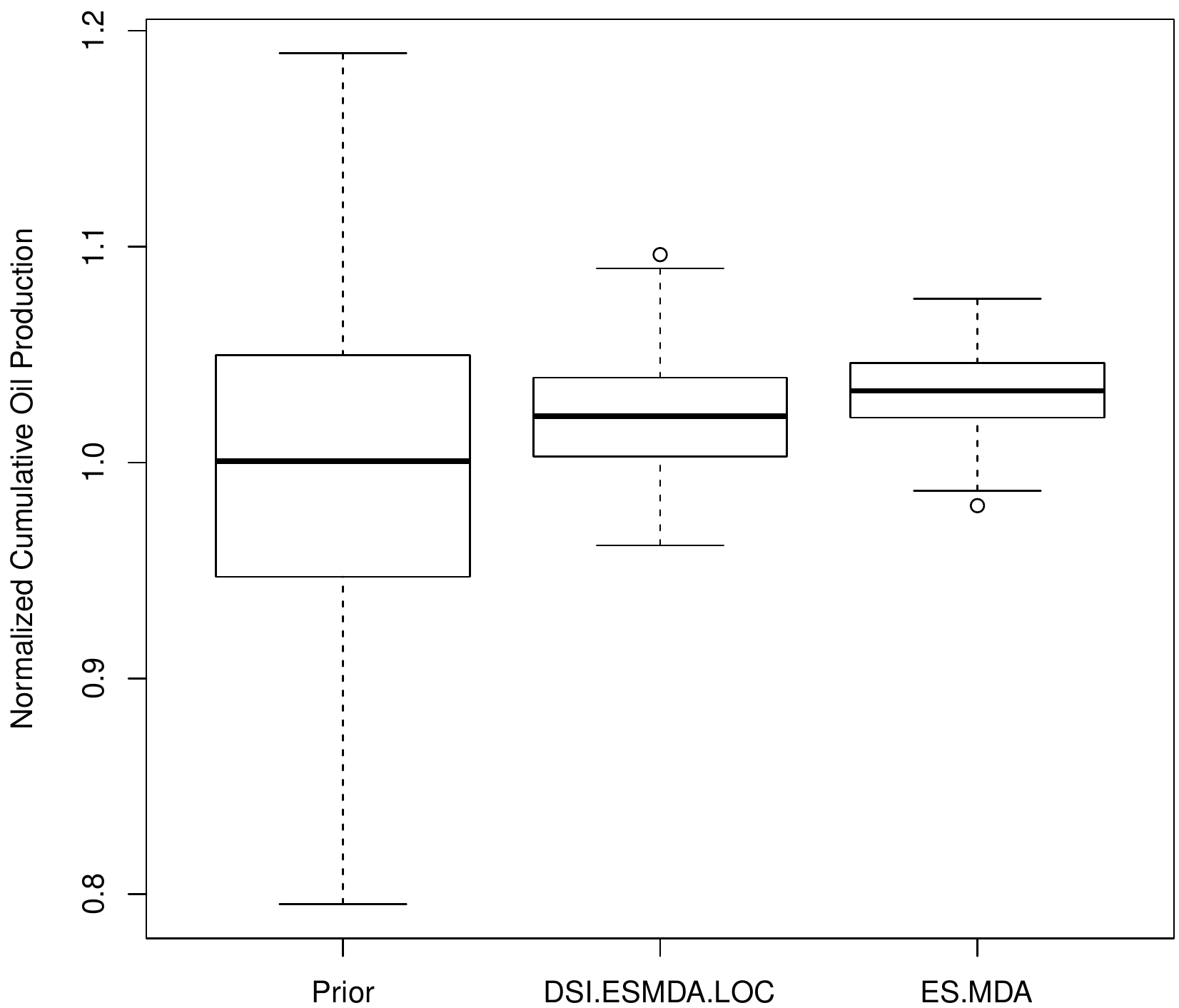}
    }
     \subfloat[]{
      \includegraphics[width=0.5\textwidth]{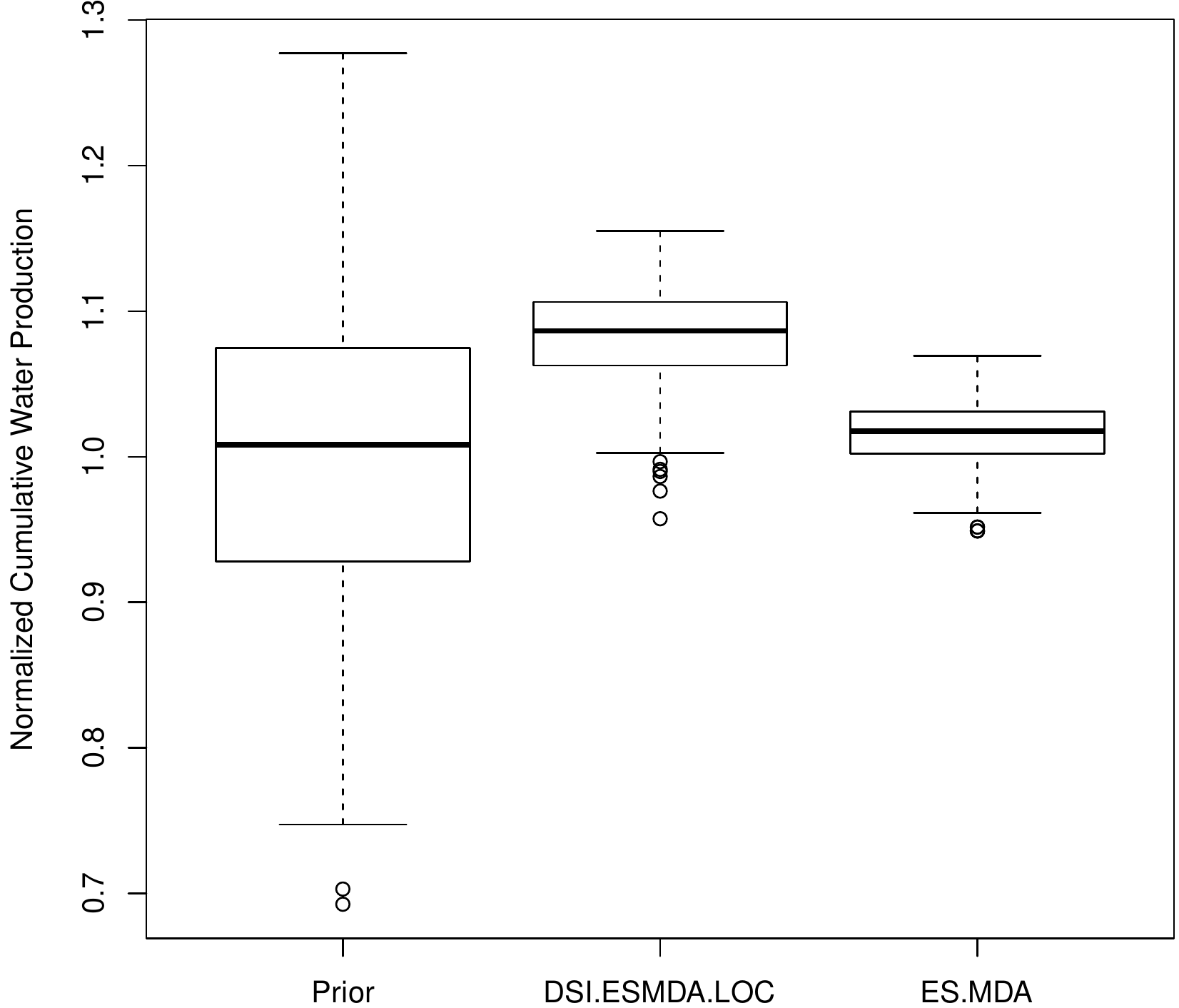}
    }
\captionsetup{justification=justified}
\caption{Normalized field cumulative production. (a) Oil and (b) water. Field case.}
\label{Fig:NpWpFieldCase}
\end{figure}

Unlike the UNISIM-I-H case, the predictions from the prior ensemble cover reasonably well the observed data for the large majority of the wells in the field, which creates a more favorable situation for the application of the methods. Once again, the optimizations in our DSI implementation failed to converge, providing unreasonable results for this field problem. DSI-ESMDA without localization showed severe reduction in the posterior variance, showing indications of ensemble colapse. Both results were not included in this paper. The only acceptable results obtained in this field are the ones from DSI-ESMDA with localization. In this case, we use the same configuration of the case UNISIM-I-H to define the critical length for localization, that is, $L_{x^\prime} = L_{y^\prime} = 2000$~m and $T = 6000$~days. For comparisons, we history matched the same ensemble of models using ES-MDA with $N_a = 4$ and localization with a critical length of 2000~m. Table~\ref{Tab:FieldCase-OF} presents the mean and standard deviation of the normalized data mismatch objective function. The results in this table indicate that DSI-ESMDA with localization obtained an excellent data match, with mean objective function of 0.55. The history matching with ES-MDA also improved significantly the data match, but with higher values of objective function. Figure~\ref{Fig:WaterRateFieldCase} shows the water production rate for eight representative wells of the field. The values in this figure were normalized to preserve the confidentiality of the information. Overall, the forecasts after DSI-ESMDA with localization are comparable with the forecasts from the history-matched models, despite some clear differences observed in some wells. Figure~\ref{Fig:NpWpFieldCase} shows boxplots of field cumulative oil and water production. The main difference between DSI-ESMDA and ES-MDA is for the cumulative water production, where DSI-ESMDA predicts a higher water production in the field.

\begin{table}
\caption{Normalized data-mismatch objective function. Field case}
\label{Tab:FieldCase-OF}
\begin{center}
\begin{tabular}{lcc}
\toprule
Case & Mean & Standard deviation\\
\midrule
Prior & 32.19 & 11.309 \\
DSI-ESMDA (with localization) & 0.55 & 0.004 \\
ES-MDA & 4.11 & 0.581 \\
\bottomrule
\end{tabular}
\end{center}
\end{table}

\section{Conclusions}
\label{Sec:Conclusions}

This paper introduced a new DSI implementation based on the data assimilation method ES-MDA. The new DSI-ESMDA was compared with the original DSI proposed in \citep{sun:17b,sun:17c}, which is based on PCA to reparameterize the predicted data from a prior ensemble combined to RML for sampling. The new implementation preserves the main advantage of the original DSI, namely, it is able to provide an ensemble of production forecasts requiring reservoir simulations only for a prior ensemble of models. We compared the DSI-ESMDA with the original DSI in two synthetic reservoir problems. We also applied DSI-ESMDA to a real field case with long production history and large number of wells. Based on the results for these test problems, we can summarize the following conclusions:

\begin{itemize}
  \item The proposed DSI-ESMDA is computationally faster than the original DSI. Even though the time to execute the reservoir simulations for the prior models tends to be dominant in both methods, the difference in the computational time for the inversion can be relevant for large problems. For example, for the UNISIM-I-H case, the original DSI required approximately 28 hours in a stand-alone computer while the DSI-ESMDA required only 45 minutes.
  \item The performance of both DSI implementations is highly dependent of the ability of the prior ensemble to provide reasonable estimates of the prior uncertainty. The same is also true for the more traditional model inversion with ES-MDA.
  \item The results indicate that the proposed DSI-ESMDA is more robust than the original DSI. In our tests, the optimizations required by DSI have not converged properly for the UNISIM-I-H case. As a result, DSI failed to obtained acceptable data matches and reasonable production forecasts. DSI-ESMDA, on the other hand, was able to match the observed data for all wells.
  \item The use of spatial and temporal localization improved significantly the results of the DSI-ESMDA method when the number of data point is large.
  \item The proposed DSI-ESMDA method with localization obtained forecasts of production comparable with the forecast provided by the posterior models generated with history matching using ES-MDA for a real field problem.
  \item The proposed method is very simple to implement. It does not require data transformations and it is straightforward integrate with different types of data and models.
\end{itemize}

\section*{Acknowledgement}
\label{Sec:Acknowledgement}
The authors would like to thank Petrobras for supporting this research and for the permission to publish this paper.


\begin{thebibliography}{52}
\providecommand{\natexlab}[1]{#1}
\providecommand{\url}[1]{\texttt{#1}}
\expandafter\ifx\csname urlstyle\endcsname\relax
  \providecommand{\doi}[1]{doi: #1}\else
  \providecommand{\doi}{doi: \begingroup \urlstyle{rm}\Url}\fi

\bibitem[Abadpour et~al.(2018)Abadpour, Adejare, Chugunova, Mathieu, and
  Haller]{abadpour:18a}
Abadpour, A., Adejare, M., Chugunova, T., Mathieu, H., and Haller, N.
\newblock Integrated geo-modeling and ensemble history matching of complex
  fractured carbonate and deep offshore turbidite fields, generation of several
  geologically coherent solutions using ensemble methods.
\newblock In \emph{Proceedings of the Abu Dhabi International Petroleum
  Exhibition \& Conference, Abu Dhabi, UAE, 12--15 November}, number
  {SPE}-193028-MS, 2018.
\newblock \doi{10.2118/193028-MS}.

\bibitem[Audet and Jr.(2006)]{audet:06a}
Audet, C. and Jr., J. E.~D.
\newblock Mesh adaptive direct search algorithms for constrained optimization.
\newblock \emph{{SIAM} Journal on Optimization}, 17\penalty0 (1):\penalty0
  188--217, 2006.
\newblock \doi{10.1137/040603371}.

\bibitem[Avansi and Schiozer(2015)]{avansi:15a}
Avansi, G.~D. and Schiozer, D.~J.
\newblock {UNISIM-I:} synthetic model for reservoir development and management
  applications.
\newblock \emph{International Journal of Modeling and Simulation for the
  Petroleum Industry}, 9\penalty0 (1):\penalty0 21--30, 2015.
\newblock URL
  \url{http://www.ijmspi.org/ojs/index.php/ijmspi/article/view/152}.

\bibitem[Chen and Oliver(2010)]{chen:10b}
Chen, Y. and Oliver, D.~S.
\newblock Cross-covariance and localization for {EnKF} in multiphase flow data
  assimilation.
\newblock \emph{Computational Geosciences}, 14\penalty0 (4):\penalty0 579--601,
  2010.
\newblock \doi{10.1007/s10596-009-9174-6}.

\bibitem[Chen and Oliver(2012)]{chen:12a}
Chen, Y. and Oliver, D.~S.
\newblock Ensemble randomized maximum likelihood method as an iterative
  ensemble smoother.
\newblock \emph{Mathematical Geosciences}, 44\penalty0 (1):\penalty0 1--26,
  2012.
\newblock \doi{10.1007/s11004-011-9376-z}.

\bibitem[Chen and Oliver(2013)]{chen:13a}
Chen, Y. and Oliver, D.~S.
\newblock {L}evenberg-{M}arquardt forms of the iterative ensemble smoother for
  efficient history matching and uncertainty quantification.
\newblock \emph{Computational Geosciences}, 17:\penalty0 689--703, 2013.
\newblock \doi{10.1007/s10596-013-9351-5}.

\bibitem[Chen and Oliver(2014)]{chen:14a}
Chen, Y. and Oliver, D.~S.
\newblock History matching of the {N}orne full-field model with an iterative
  ensemble smoother.
\newblock \emph{SPE Reservoir Evaluation \& Engineering}, 17\penalty0 (2),
  2014.
\newblock \doi{10.2118/164902-PA}.

\bibitem[Emerick(2014)]{emerick:14a}
Emerick, A.~A.
\newblock Estimation of pressure and saturation fields from time-lapse
  impedance data using the ensemble smoother.
\newblock \emph{Journal of Geophysics and Engineering}, 11\penalty0
  (3):\penalty0 035007, 2014.
\newblock \doi{10.1088/1742-2132/11/3/035007}.

\bibitem[Emerick(2016)]{emerick:16a}
Emerick, A.~A.
\newblock Analysis of the performance of ensemble-based assimilation of
  production and seismic data.
\newblock \emph{Journal of Petroleum Science and Engineering}, 139:\penalty0
  219--239, 2016.
\newblock \doi{10.1016/j.petrol.2016.01.029}.

\bibitem[Emerick(2018)]{emerick:18b}
Emerick, A.~A.
\newblock Analysis of geometric selection of the data-error covariance
  inflation for {ES-MDA}.
\newblock \emph{arXiv:1812.00924v1 [math.NA]}, 2018.
\newblock URL \url{https://arxiv.org/abs/1812.00924}.

\bibitem[Emerick and Reynolds(2011{\natexlab{a}})]{emerick:11c}
Emerick, A.~A. and Reynolds, A.~C.
\newblock History matching a field case using the ensemble {K}alman filter with
  covariance localization.
\newblock \emph{SPE Reservoir Evaluation \& Engineering}, 14\penalty0
  (4):\penalty0 423--432, 2011{\natexlab{a}}.
\newblock \doi{10.2118/141216-PA}.

\bibitem[Emerick and Reynolds(2011{\natexlab{b}})]{emerick:11e}
Emerick, A.~A. and Reynolds, A.~C.
\newblock Combining sensitivities and prior information for covariance
  localization in the ensemble {K}alman filter for petroleum reservoir
  applications.
\newblock \emph{Computational Geosciences}, 15\penalty0 (2):\penalty0 251--269,
  2011{\natexlab{b}}.
\newblock \doi{10.1007/s10596-010-9198-y}.

\bibitem[Emerick and Reynolds(2012)]{emerick:12a}
Emerick, A.~A. and Reynolds, A.~C.
\newblock History matching time-lapse seismic data using the ensemble {K}alman
  filter with multiple data assimilations.
\newblock \emph{Computational Geosciences}, 16\penalty0 (3):\penalty0 639--659,
  2012.
\newblock \doi{10.1007/s10596-012-9275-5}.

\bibitem[Emerick and Reynolds(2013{\natexlab{a}})]{emerick:13a}
Emerick, A.~A. and Reynolds, A.~C.
\newblock Investigation on the sampling performance of ensemble-based methods
  with a simple reservoir model.
\newblock \emph{Computational Geosciences}, 17\penalty0 (2):\penalty0 325--350,
  2013{\natexlab{a}}.
\newblock \doi{10.1007/s10596-012-9333-z}.

\bibitem[Emerick and Reynolds(2013{\natexlab{b}})]{emerick:13b}
Emerick, A.~A. and Reynolds, A.~C.
\newblock Ensemble smoother with multiple data assimilation.
\newblock \emph{Computers \& Geosciences}, 55:\penalty0 3--15,
  2013{\natexlab{b}}.
\newblock \doi{10.1016/j.cageo.2012.03.011}.

\bibitem[Emerick and Reynolds(2013{\natexlab{c}})]{emerick:13c}
Emerick, A.~A. and Reynolds, A.~C.
\newblock History matching of production and seismic data for a real field case
  using the ensemble smoother with multiple data assimilation.
\newblock In \emph{Proceedings of the SPE Reservoir Simulation Symposium, The
  Woodlands, Texas, USA, 18--20 February}, number {SPE}-163675-MS,
  2013{\natexlab{c}}.
\newblock \doi{10.2118/163675-MS}.

\bibitem[Evensen(1994)]{evensen:94}
Evensen, G.
\newblock Sequential data assimilation with a nonlinear quasi-geostrophic model
  using {M}onte {C}arlo methods to forecast error statistics.
\newblock \emph{Journal of Geophysical Research}, 99\penalty0 (C5):\penalty0
  10143--10162, 1994.
\newblock \doi{10.1029/94JC00572}.

\bibitem[Evensen(2004)]{evensen:04}
Evensen, G.
\newblock Sampling strategies and square root analysis schemes for the {EnKF}.
\newblock \emph{Ocean Dynamics}, 54\penalty0 (6):\penalty0 539--560, 2004.
\newblock \doi{10.1007/s10236-004-0099-2}.

\bibitem[Evensen(2007)]{evensen:07bk}
Evensen, G.
\newblock \emph{Data assimilation: the ensemble {K}alman filter}.
\newblock Springer, Berlin, 2007.

\bibitem[Evensen(2018)]{evensen:18a}
Evensen, G.
\newblock Analysis of iterative ensemble smoothers for solving inverse
  problems.
\newblock \emph{Computational Geosciences}, 2018.
\newblock \doi{10.1007/s10596-018-9731-y}.

\bibitem[Gaspari and Cohn(1999)]{gaspari:99}
Gaspari, G. and Cohn, S.~E.
\newblock Construction of correlation functions in two and three dimensions.
\newblock \emph{Quarterly Journal of the Royal Meteorological Society},
  125\penalty0 (554):\penalty0 723--757, 1999.
\newblock \doi{10.1002/qj.49712555417}.

\bibitem[He et~al.(2018)He, Sarma, Bhark, Tanaka, Chen, Wen, and
  Kamath]{he:18b}
He, J., Sarma, P., Bhark, E., Tanaka, S., Chen, B., Wen, X.-H., and Kamath, J.
\newblock Quantifying expected uncertainty reduction and value of information
  using ensemble-variance analysis.
\newblock \emph{SPE Journal}, 23\penalty0 (2), 2018.
\newblock \doi{10.2118/182609-PA}.

\bibitem[Houtekamer and Mitchell(2001)]{houtekamer:01}
Houtekamer, P.~L. and Mitchell, H.~L.
\newblock A sequential ensemble {K}alman filter for atmospheric data
  assimilation.
\newblock \emph{Monthly Weather Review}, 129\penalty0 (1):\penalty0 123--137,
  2001.
\newblock \doi{10.1175/1520-0493(2001)129<0123:ASEKFF>2.0.CO;2}.

\bibitem[Jeong et~al.(2018)Jeong, Sun, Lee, and Min]{jeong:18a}
Jeong, H., Sun, A.~Y., Lee, J., and Min, B.
\newblock A learning-based data-driven forecast approach for predicting future
  reservoir performance.
\newblock \emph{Advances in Water Resources}, 118:\penalty0 95--109, 2018.
\newblock \doi{10.1016/j.advwatres.2018.05.015}.

\bibitem[Jiang(2018)]{jiang:18a}
Jiang, S.
\newblock Data-space inversion with variable well controls in the prediction
  period.
\newblock Master's thesis, Stanford University, 2018.

\bibitem[Krishnamurti et~al.(2000)Krishnamurti, Kishtawal, Zhang, LaRow,
  Bachiochi, and Williford]{krishnamurti:00a}
Krishnamurti, T.~N., Kishtawal, C.~M., Zhang, Z., LaRow, T., Bachiochi, D., and
  Williford, E.
\newblock Multimodel ensemble forecasts for weather and seasonal climate.
\newblock \emph{Journal of Climate}, 13\penalty0 (23):\penalty0 4196--4216,
  2000.
\newblock \doi{10.1175/1520-0442(2000)013<4196:MEFFWA>2.0.CO;2}.

\bibitem[Le et~al.(2016)Le, Emerick, and Reynolds]{duc:16a}
Le, D.~H., Emerick, A.~A., and Reynolds, A.~C.
\newblock An adaptive ensemble smoother with multiple data assimilation for
  assisted history matching.
\newblock \emph{SPE Journal}, 21\penalty0 (6):\penalty0 2195--2207, 2016.
\newblock \doi{10.2118/173214-PA}.

\bibitem[Lorentzen et~al.(2019)Lorentzen, Luo, Bhakta, and
  Valestrand]{lorentzen:19a}
Lorentzen, R.~J., Luo, X., Bhakta, T., and Valestrand, R.
\newblock History matching the full norne field model using seismic and
  production data.
\newblock \emph{SPE Journal}, Preprint, 2019.
\newblock \doi{10.2118/194205-PA}.

\bibitem[Luo et~al.(2015)Luo, Stordal, Lorentzen, and N{\ae}vdal]{luo:15b}
Luo, X., Stordal, A.~S., Lorentzen, R.~J., and N{\ae}vdal, G.
\newblock Iterative ensemble smoother as an approximate solution to a
  regularized minimum-average-cost problem: Theory and applications.
\newblock \emph{SPE Journal}, 20\penalty0 (5), 2015.
\newblock \doi{10.2118/176023-PA}.

\bibitem[Ma et~al.(2017)Ma, Hetz, Wang, Bi, Stern, and Hoda]{ma:17a}
Ma, X., Hetz, G., Wang, X., Bi, L., Stern, D., and Hoda, N.
\newblock A robust iterative ensemble smoother method for efficient history
  matching and uncertainty quantification.
\newblock In \emph{Proceedings of the SPE Reservoir Simulation Conference,
  Montgomery, Texas, USA, 20--22 February}, number {SPE}-182693-MS, 2017.
\newblock \doi{10.2118/182693-MS}.

\bibitem[Mallet et~al.(2009)Mallet, Stoltz, and Mauricette]{mallet:09a}
Mallet, V., Stoltz, G., and Mauricette, B.
\newblock Ozone ensemble forecast with machine learning algorithms.
\newblock \emph{Journal of Geophysical Research: Atmospheres}, 114\penalty0
  (D5), 2009.
\newblock \doi{10.1029/2008JD009978}.

\bibitem[Maschio et~al.(2013)Maschio, Avansi, Santos, and
  Schiozer]{unisim-i-h:13}
Maschio, C., Avansi, G.~D., Santos, A.~A., and Schiozer, D.~J.
\newblock {UNISIM-I-H:} case study for history matching.
\newblock Dataset, 2013.
\newblock URL
  \url{www.unisim.cepetro.unicamp.br/benchmarks/br/unisim-i/unisim-i-h}.

\bibitem[Maucec et~al.(2016)Maucec, Ravanelli, Lyngra, Zhang, Alramadhan,
  Abdelhamid, and Al-Garni]{maucec:16a}
Maucec, M., Ravanelli, F. M. D.~M., Lyngra, S., Zhang, S.~J., Alramadhan,
  A.~A., Abdelhamid, O.~A., and Al-Garni, S.~A.
\newblock Ensemble-based assisted history matching with rigorous uncertainty
  quantification applied to a naturally fractured carbonate reservoir.
\newblock In \emph{Proceedings of the SPE Annual Technical Conference and
  Exhibition, Dubai, UAE, 26--28 September}, number {SPE}-181325-MS, 2016.
\newblock \doi{10.2118/181325-MS}.

\bibitem[Nocedal and Wright(2006)]{nocedal:06}
Nocedal, J. and Wright, S.~J.
\newblock \emph{Numerical Optimization}.
\newblock Springer, New York, 2006.

\bibitem[Oliver and Alfonzo(2018)]{oliver:18b}
Oliver, D.~S. and Alfonzo, M.
\newblock Calibration of imperfect models to biased observations.
\newblock \emph{Computational Geosciences}, 22:\penalty0 145--161, 2018.
\newblock \doi{10.1007/s10596-017-9678-4}.

\bibitem[Oliver and Chen(2011)]{oliver:11b}
Oliver, D.~S. and Chen, Y.
\newblock Recent progress on reservoir history matching: a review.
\newblock \emph{Computational Geosciences}, 15\penalty0 (1):\penalty0 185--221,
  2011.
\newblock \doi{10.1007/s10596-010-9194-2}.

\bibitem[Oliver et~al.(2008)Oliver, Reynolds, and Liu]{oliver:08bk}
Oliver, D.~S., Reynolds, A.~C., and Liu, N.
\newblock \emph{Inverse Theory for Petroleum Reservoir Characterization and
  History Matching}.
\newblock Cambridge University Press, Cambridge, UK, 2008.

\bibitem[Pagowski et~al.(2005)Pagowski, Grell, McKeen, D{\'e}v{\'e}nyi,
  Wilczak, Bouchet, Gong, Mchenry, Peckham, Mcqueen, Moffet, and
  Tang]{pagowski:05a}
Pagowski, M., Grell, G.~A., McKeen, S.~A., D{\'e}v{\'e}nyi, Wilczak, J.~M.,
  Bouchet, V.~S., Gong, W.~F., Mchenry, J.~N., Peckham, S., Mcqueen, J.~T.,
  Moffet, R., and Tang, Y.
\newblock A simple method to improve ensemble-based ozone forecasts.
\newblock \emph{Geophysical Research Letters}, 320\penalty0 (7), 2005.
\newblock \doi{10.1029/2004GL022305}.

\bibitem[Rafiee and Reynolds(2017)]{rafiee:17a}
Rafiee, J. and Reynolds, A.~C.
\newblock Theoretical and efficient practical procedures for the generation of
  inflation factors for {ES-MDA}.
\newblock \emph{Inverse Problems}, 33\penalty0 (11):\penalty0 115003, 2017.
\newblock \doi{10.1088/1361-6420/aa8cb2}.

\bibitem[Reynolds et~al.(2006)Reynolds, Zafari, and Li]{reynolds:06}
Reynolds, A.~C., Zafari, M., and Li, G.
\newblock Iterative forms of the ensemble {K}alman filter.
\newblock In \emph{Proceedings of 10th European Conference on the Mathematics
  of Oil Recovery, Amsterdam, 4--7 September}, 2006.
\newblock \doi{10.3997/2214-4609.201402496}.

\bibitem[Satija and Caers(2015)]{satija:15a}
Satija, A. and Caers, J.
\newblock Direct forecasting of subsurface flow response from non-linear
  dynamic data by linear least-squares in canonical functional principal
  component space.
\newblock \emph{Advances in Water Resources}, 77:\penalty0 69--81, 2015.
\newblock \doi{10.1016/j.advwatres.2015.01.002}.

\bibitem[Satija and Caers(2017)]{satija:17a}
Satija, A. and Caers, J.
\newblock Direct forecasting of reservoir performance using production data
  without history matching.
\newblock \emph{Computational Geosciences}, 21\penalty0 (2):\penalty0 315--333,
  2017.
\newblock \doi{10.1007/s10596-017-9614-7}.

\bibitem[Scheidt et~al.(2015)Scheidt, Renard, and Caers]{scheidt:15a}
Scheidt, C., Renard, P., and Caers, J.
\newblock Prediction-focused subsurface modeling: Investigating the need for
  accuracy in flow-based inverse modeling.
\newblock \emph{Mathematical Geosciences}, 47\penalty0 (2):\penalty0 173--191,
  2015.
\newblock \doi{10.1007/s11004-014-9521-6}.

\bibitem[Skjervheim et~al.(2011)Skjervheim, Evensen, Hove, and
  Vab{\o}]{skjervheim:11a}
Skjervheim, J.-A., Evensen, G., Hove, J., and Vab{\o}, J.~G.
\newblock An ensemble smoother for assisted history matching.
\newblock In \emph{Proceedings of the SPE Reservoir Simulation Symposium, The
  Woodlands, Texas, USA, 21--23 February}, number {SPE}-141929-MS, 2011.
\newblock \doi{10.2118/141929-MS}.

\bibitem[Souza(2017)]{souza:17a}
Souza, C.~R.
\newblock {A}ccord.{NET} {F}ramework {B}royden-{F}letcher-{G}oldfarb-{S}hanno
  class, 2017.
\newblock URL
  \url{http://accord-framework.net/docs/html/T_Accord_Math_Optimization_BroydenFletcherGoldfarbShanno.htm}.

\bibitem[Stordal(2014)]{stordal:14b}
Stordal, A.~S.
\newblock Iterative {B}ayesian inversion with {G}aussian mixtures: finite
  sample implementation and large sample asymptotics.
\newblock \emph{Computational Geosciences}, 19\penalty0 (1):\penalty0 1--15,
  2014.
\newblock \doi{10.1007/s10596-014-9444-9}.

\bibitem[Stordal and Elsheikh(2015)]{stordal:15a}
Stordal, A.~S. and Elsheikh, A.~H.
\newblock Iterative ensemble smoothers in the annealed importance sampling
  framework.
\newblock \emph{Advances in Water Resources}, 86:\penalty0 231--239, 2015.
\newblock \doi{10.1016/j.advwatres.2015.09.030}.

\bibitem[Sun and Durlofsky(2017)]{sun:17b}
Sun, W. and Durlofsky, L.~J.
\newblock A new data-space inversion procedure for efficient uncertainty
  quantification in subsurface flow problems.
\newblock \emph{Mathematical Geosciences}, Online, 2017.
\newblock \doi{10.1007/s11004-016-9672-8}.

\bibitem[Sun et~al.(2017{\natexlab{a}})Sun, Hui, and Durlofsky]{sun:17c}
Sun, W., Hui, M.-H., and Durlofsky, L.~J.
\newblock Production forecasting and uncertainty quantification for naturally
  fractured reservoirs using a new data-space inversion procedure.
\newblock \emph{Computational Geosciences}, Online, 2017{\natexlab{a}}.
\newblock \doi{10.1007/s10596-017-9633-4}.

\bibitem[Sun et~al.(2017{\natexlab{b}})Sun, Vink, and Gao]{sun:17a}
Sun, W., Vink, J.~C., and Gao, G.
\newblock A practical method to mitigate spurious uncertainty reduction in
  history matching workflows with imperfect reservoir models.
\newblock In \emph{Proceedings of the SPE Reservoir Simulation Conference,
  Montgomery, Texas, USA, 20-22 February}, number {SPE}-182599-MS,
  2017{\natexlab{b}}.
\newblock \doi{10.2118/182599-MS}.

\bibitem[Tarantola(2005)]{tarantola:05}
Tarantola, A.
\newblock \emph{Inverse Problem Theory and Methods for Model Parameter
  Estimation}.
\newblock {SIAM}, Philadelphia, {USA}, 2005.

\bibitem[van Leeuwen and Evensen(1996)]{vanleeuwen:96}
van Leeuwen, P.~J. and Evensen, G.
\newblock Data assimilation and inverse methods in terms of a probabilistic
  formulation.
\newblock \emph{Monthly Weather Review}, 124:\penalty0 2898--2913, 1996.
\newblock \doi{10.1175/1520-0493(1996)124<2898:DAAIMI>2.0.CO;2}.

\end{thebibliography}

\end{document}